\documentclass[12pt,twoside,leqno,openany]{amsart}

\pdfoutput=1

\usepackage{amssymb,amsbsy,amsmath,amsfonts,amssymb,amscd,times,
graphics,color,xypic,footmisc,fancyhdr,multicol,fancybox,
graphicx,mathrsfs,rotating,ifthen,wasysym}
\usepackage[all]{xy}

\usepackage[T1]{fontenc}
\newcommand{\C}{\mathbb{C}}

\newcommand{\N}{\mathbb{N}}
\newcommand{\R}{\mathbb{R}}

\newcommand{\HEAD}[2]{%
\pagestyle{fancy}
\fancyhead[RO]{\tiny\sf\thepage}
\fancyhead[CO]{{\tiny\sf #1}}
\fancyhead[LE]{\tiny\sf\thepage}
\fancyhead[CE]{{\tiny\sf #2}}
\fancyfoot{}}

\theoremstyle{definition}

\newcommand{\codim}{\text{\footnotesize\sf codim}}

\newcommand{\constant}{\text{\scriptsize\sf constant}}

\newcommand{\CRdim}{\text{\footnotesize\sf CRdim}}

\renewcommand{\dim}{\text{\footnotesize\sf dim}}

\newcommand{\isqrt}{{\scriptstyle{\sqrt{-1}}}}

\renewcommand{\lim}{\text{\footnotesize\sf lim}}

\newcommand{\mathmotsf}[1]{\text{\footnotesize\sf #1}}

\renewcommand{\mod}{\text{\footnotesize\sf mod}}

\newcommand{\NN}{\text{\scriptsize\sc n}}

\newcommand{\NNN}{\text{\sc n}}

\newcommand{\rank}{\text{\footnotesize\sf rank}}

\renewcommand{\Re}{\text{\footnotesize\sf Re}}

\newcommand{\remainder}{\text{\scriptsize\sf remainder}}

\newcommand{\smallbullet}{{\scriptscriptstyle{\bullet}}}

\newcommand{\vf}{\vfill\end{document}}

\newcommand{\weight}{\text{\footnotesize\sf weight}}

\let\mathcal\mathscr

\begin{document}

$\:$

\bigskip\bigskip

\begin{center}

{\Large\bf Equivalences of $5$-dimensional CR-manifolds}

\medskip

{\Large\bf III: Six models and (very) elementary}

\medskip

{\Large\bf normalizations}

\end{center}

\medskip

\begin{center}
Jo\"el {\sc Merker}
\end{center}

\bigskip

\begin{center}
\begin{minipage}[t]{10.25cm}
\baselineskip =0.32cm 
{\scriptsize
{\bf Abstract.}
The six nondegeneracy conditions of geometric nature that are
satisfied by the only six possibly existing nondegenerate general classes
$\text{\sf I}$,
$\text{\sf II}$, 
$\text{\sf III}_{\text{\sf 1}}$,
$\text{\sf III}_{\text{\sf 2}}$,
$\text{\sf IV}_{\text{\sf 1}}$,
$\text{\sf IV}_{\text{\sf 2}}$
of 5-dimensional CR manifolds are shown to be readable
instantaneously from their elementarily normalized respective defining
graphed equations, without advanced Moser theory.
}
\end{minipage}
\end{center}

\medskip

\begin{center}
\begin{minipage}[t]{11.75cm}
\baselineskip =0.35cm {\scriptsize

\centerline{\bf Table of contents}

\medskip

{\bf \ref{introduction}.~Introduction
\dotfill~\pageref{introduction}.}

{\bf \ref{removing-pluriharmonic}.~Removing pluriharmonic monomials
\dotfill~\pageref{removing-pluriharmonic}.}

{\bf \ref{general-class-I}.~General class $\text{\sf I}$
\dotfill~\pageref{general-class-I}.}

{\bf \ref{general-class-II}.~General class $\text{\sf II}$
\dotfill~\pageref{general-class-II}.}

{\bf \ref{general-class-III-1}.~General class 
$\text{\sf III}_{\text{\sf 1}}$
\dotfill~\pageref{general-class-III-1}.}

{\bf \ref{general-class-III-2}.~General class 
$\text{\sf III}_{\text{\sf 2}}$
\dotfill~\pageref{general-class-III-2}.}

{\bf \ref{general-class-IV-1}.~General class 
$\text{\sf IV}_{\text{\sf 1}}$
\dotfill~\pageref{general-class-IV-1}.}

{\bf \ref{general-class-IV-2}.~General class 
$\text{\sf IV}_{\text{\sf 2}}$
\dotfill~\pageref{general-class-IV-2}.}

{\bf \ref{smooth-categories}.~Smooth categories
\dotfill~\pageref{smooth-categories}.}

{\bf \ref{graphing-ontologies}.~Graphing ontologies
\dotfill~\pageref{graphing-ontologies}.}

}\end{minipage}
\end{center}


\bigskip

\section{\sf Introduction}
\label{introduction}
\HEAD{\ref{introduction}.~Introduction}{
Jo\"el {\sc Merker} (Paris-Sud)}

\medskip

Consider a CR-generic submanifold:
\[
\aligned
&
M^{2n+c}
\,\subset\,
\C^{n+c}
\\
&
{\scriptstyle{(c\,=\,{\sf codim}\,M,\,\,\,
n\,=\,{\sf CRdim}\,M)}},
\endaligned
\]
of smoothness:
\[
\mathcal{C}^\kappa\,
(\kappa\geqslant 1),
\ \ \ \ \ \ \ \ \ \ \ \ \ \ \
\text{\rm or}\,\,\,\,
\mathcal{C}^\infty,
\ \ \ \ \ \ \ \ \ \ \ \ \ \ \
\text{\rm or}\,\,\,\,
\mathcal{C}^\omega.
\]

According to~\cite{ Merker-Pocchiola-Sabzevari-5-CR-II}, there are
precisely {\em six} general classes of nondegenerate $M^{ 2n +c}
\subset \C^{ n+c}$ having dimension:
\[
2n+c
\,\leqslant\,
{\bf 5},
\]
hence having CR dimension: 
\[
n\,=\,
\left\{
\aligned
&
{\bf 1},
\\
&
{\bf 2}.
\endaligned\right.
\]
namely if one denotes by:
\[
\aligned
&
\big\{
\mathcal{L}
\big\},
\\
&
\big\{
\mathcal{L}_1,\mathcal{L}_2
\big\},
\endaligned
\]
any local frame for $T^{ 1, 0} M$, firstly the well known:
\[
\aligned
&
\boxed{\text{\sf General Class $\text{\sf I}$:}}
\\
&
\boxed{\,\,
\aligned
M^3\subset\C^2
\ \
&
\text{\rm with}\ \
\Big\{\mathcal{L},\,\overline{\mathcal{L}},\,\,
\big[\mathcal{L},\overline{\mathcal{L}}\big]\Big\}\ \
\\
&
\text{\rm constituting a frame for}\ \
\C\otimes_\R TM,\,\,
\endaligned}
\endaligned
\]
secondly the known:
\[
\aligned
&
\boxed{\text{\sf General Class $\text{\sf II}$:}}
\\
&
\boxed{\,\,
\aligned
M^4\subset\C^3
\ \
&
\text{\rm with}\ \
\Big\{\mathcal{L},\,\overline{\mathcal{L}},\,\,
\big[\mathcal{L},\overline{\mathcal{L}}\big],\,\,
\big[\mathcal{L},\,\big[\mathcal{L},\overline{\mathcal{L}}\big]\big]
\Big\}\ \
\\
&
\text{\rm constituting a frame for}\ \
\C\otimes_\R TM,\,\,
\endaligned}
\endaligned
\]
thirdly the known:
\[
\aligned
&
\boxed{\text{\sf General Class $\text{\sf III}_{\text{\sf 1}}$:}}
\\
&
\boxed{\,\,
\aligned
M^5\subset\C^4
\ \
&
\text{\rm with}\ \
\Big\{\mathcal{L},\,\overline{\mathcal{L}},\,\,
\big[\mathcal{L},\overline{\mathcal{L}}\big],\,\,
\big[\mathcal{L},\,
\big[\mathcal{L},\overline{\mathcal{L}}\big]\big],\,\,
\big[\overline{\mathcal{L}},\,
\big[\mathcal{L},\overline{\mathcal{L}}\big]\big]
\Big\}\,\,
\\
&
\text{\rm constituting a frame for}\ \
\C\otimes_\R TM,
\endaligned}
\endaligned
\]
fourthly the {\em new:}
\[
\aligned
&
\boxed{\text{\sf General Class $\text{\sf III}_{\text{\sf 2}}$:}}
\\
&
\boxed{\,\,
\aligned
M^5\subset\C^4
\ \
&
\text{\rm with}\ \
\Big\{\mathcal{L},\,\overline{\mathcal{L}},\,\,
\big[\mathcal{L},\overline{\mathcal{L}}\big],\,\,
\big[\mathcal{L},\,
\big[\mathcal{L},\overline{\mathcal{L}}\big]\big],\,\,
\big[\mathcal{L},\,\big[\mathcal{L},\,
\big[\mathcal{L},\overline{\mathcal{L}}\big]\big]\big]
\Big\}\,\,
\\
&
\text{\rm constituting a frame for}\ \
\C\otimes_\R TM,
\\
&
\text{\rm while}\ \
{\bf 4}
=
\rank_\C
\Big(\mathcal{L},\overline{\mathcal{L}},\,
\big[\mathcal{L},\overline{\mathcal{L}}\big],\,\,
\big[\mathcal{L},\,\big[\mathcal{L},\overline{\mathcal{L}}\big]\big],\,\,
\big[\overline{\mathcal{L}},\,
\big[\mathcal{L},\overline{\mathcal{L}}\big]\big]\Big),
\endaligned}
\endaligned
\]
fifthly the well known:
\[
\aligned
&
\boxed{\text{\sf General Class $\text{\sf IV}_{\text{\sf 1}}$:}}
\\
&
\boxed{\,\,
\aligned
M^5\subset\C^3
\ \
&
\text{\rm with}\ \
\Big\{\mathcal{L}_1,\,\mathcal{L}_2,\,
\overline{\mathcal{L}}_1,\,\overline{\mathcal{L}}_2,\,\,
\big[\mathcal{L}_1,\overline{\mathcal{L}}_1\big]
\Big\}\,\,
\\
&
\text{\rm constituting a frame for}\ \
\C\otimes_\R TM,
\\
&
\text{\rm and with the Levi-Form:}\ \
\\
&
\ \ \ \ \ \ \ \ \ \ \ \ \ \ \
\mathmotsf{Levi-Form}^M(p)
\\
&
\text{\rm being of rank}\,\,{\bf 2}\,\,
\text{\rm at every point}\,\,p\in M,
\endaligned}
\endaligned
\]
sixthly and lastly the known:

\[
\aligned
&
\boxed{\text{\sf General Class $\text{\sf IV}_{\text{\sf 2}}$:}}
\\
&
\boxed{\,\,
\aligned
M^5\subset\C^3
\ \
&
\text{\rm with}\ \
\Big\{\mathcal{L}_1,\,\mathcal{L}_2,\,
\overline{\mathcal{L}}_1,\,\overline{\mathcal{L}}_2,\,\,
\big[\mathcal{L}_1,\overline{\mathcal{L}}_1\big]
\Big\}\,\,
\\
&
\text{\rm constituting a frame for}\ \
\C\otimes_\R TM,
\\
&
\text{\rm with the Levi-Form:}\ \
\\
&
\ \ \ \ \ \ \ \ \ \ \ \ \ \ \
\mathmotsf{Levi-Form}^M(p)
\\
&
\text{\rm being of rank}\,\,{\bf 1}\,\,
\text{\rm at every point}\,\,p\in M,
\\
&
\text{\rm while the Freeman-Form:}\ \
\\
&
\ \ \ \ \ \ \ \ \ \ \ \ \ \ \
\mathmotsf{Freeman-Form}^M(p)
\\
&
\text{\rm is nondegenerate at every point.}\ \
\endaligned}
\endaligned
\]

\medskip

The objective is to find elementary, preliminary normalizing
coordinates in which the graphing equations would immediately
and visibly show
that these six nondegeneracy conditions are satisfied.

\medskip\noindent{\bf Remarks.}
For background foundational material, the reader is referred to Part
II of this memoir (\cite{ Merker-Pocchiola-Sabzevari-5-CR-II}),
considered to be conceptually synthetized.

Part~I, scheduled at the end, will include introduction,
references, credit.

Part~IV (forthcoming, brief) will set up
the six related matrix ambiguity groups (initial $G$-structures).

Part~V (forthcoming, {\em door to the core}), will set up 
the six initial coframes and their Darboux structures.

These thankless systematic re-foundations being achieved, 
the final launch of central computations will be firmly
grounded.


\bigskip

\section{\sf Removing pluriharmonic monomials}
\label{removing-pluriharmonic}
\HEAD{\ref{removing-pluriharmonic}.~Removing pluriharmonic monomials}{
Jo\"el {\sc Merker} (Paris-Sud)}

\medskip

Pick two integers:
\[
n\geqslant 1,
\ \ \ \ \ \ \ \ \ \ \ \ \ \ \ \ \ \
c\geqslant 1,
\]
denote coordinates on $\C^{n+c}$ by:
\[
(z_1,\dots,z_n,w_1,\dots,w_c)
\,=\,
\big(
x_1+\isqrt\,y_1,\dots,x_n+\isqrt\,y_n,\,
u_1+\isqrt\,v_1,\dots,u_c+\isqrt\,v_c
\big),
\]
and notationally contract them as:
\[
\big(z_\smallbullet,w_\smallbullet\big)
=
\big(
x_\smallbullet+\isqrt\,y_\smallbullet,\,\,
u_\smallbullet+\isqrt\,v_\smallbullet
\big).
\]

Given therefore a {\em local} CR-generic {\em real analytic} 
($\mathcal{ C}^\omega$) submanifold:
\[
M^{2n+c}
\subset
\C^{n+c}
\]
passing through the origin: 
\[
0
\in
M,
\]
with of course:
\[
\aligned
c
&
=
\codim_\R\,M,
\\
n
&
=
\CRdim\,M,
\endaligned
\]
then provided only that:
\[
\frac{\partial}{\partial z_1}
\bigg\vert_0
\,\not\in\,
T_0^{1,0}M,
\,\,\dots\dots,\,\,
\frac{\partial}{\partial z_n}
\bigg\vert_0
\,\not\in\,
T_0^{1,0}M,
\]
one easily convinces oneself that the analytic implicit function
theorem shows the existence of $c$ real analytic graphing functions:
\[
\varphi_1\big(z_\smallbullet,\overline{z}_\smallbullet,u_\smallbullet\big),
\,\dots\dots\dots,\,
\varphi_c\big(z_\smallbullet,\overline{z}_\smallbullet,u_\smallbullet\big),
\]
having convergent power series in some neighborhood of:
\[
\big(0_\smallbullet,0_\smallbullet,0_\smallbullet\big)
\,\in\,
\C^n\times\C^n\times\R^c,
\]
and vanishing at the origin:
\[
\aligned
0
&
=
\varphi_1\big(0_\smallbullet,0_\smallbullet,0_\smallbullet\big),
\\
\cdot\cdot
&
\cdots\cdots\cdots\cdots\cdots
\\
0
&
=
\varphi_c\big(0_\smallbullet,0_\smallbullet,0_\smallbullet\big),
\endaligned
\]
such that $M$ is locally represented by:
\[
\left[
\aligned
v_1
&
=
\varphi_1\big(z_1,\dots,z_n,\overline{z}_1,\dots,
\overline{z}_n,u_1,\dots,u_c\big),
\\
\cdots
&
\cdots\cdots\cdots\cdots\cdots\cdots\cdots\cdots\cdots\cdots\cdots\cdots
\cdot
\\
v_c
&
=
\varphi_c\big(z_1,\dots,z_n,\overline{z}_1,\dots,
\overline{z}_n,u_1,\dots,u_c\big).
\endaligned\right.
\]

After a $\C$-linear biholomorphism of $\C^{ n+c}$, one
can even assume that:
\[
T_0M
=
\big\{
0
=
v_1
=\cdots=
v_c
\big\},
\]
so that the $c$ graphing functions satisfy in addition:
\[
\aligned
0
&
=
\varphi_{1,z_k}\big(0_\smallbullet,0_\smallbullet,0_\smallbullet\big)
=
\varphi_{1,\overline{z}_k}\big(0_\smallbullet,0_\smallbullet,0_\smallbullet\big)
=
\varphi_{1,u_l}\big(0_\smallbullet,0_\smallbullet,0_\smallbullet\big)
\\
\cdot\,\cdot
&
\cdots\cdots\cdots\cdots\cdots\cdots\cdots\cdots\cdots\cdots
\cdots\cdots\cdots\cdots\cdots\cdots\cdots
\\
0
&
=
\varphi_{c,z_k}\big(0_\smallbullet,0_\smallbullet,0_\smallbullet\big)
=
\varphi_{c,\overline{z}_k}\big(0_\smallbullet,0_\smallbullet,0_\smallbullet\big)
=
\varphi_{c,u_l}\big(0_\smallbullet,0_\smallbullet,0_\smallbullet\big)
\\
&
\ \ \ \ \ \ \ \ \ \ \ \ \ \ \ \ \ \ \ \ \ \ \ \ \ \
{\scriptstyle{(1\,\leqslant\,k\,\leqslant\,n;\,\,\,
1\,\leqslant\,l\,\leqslant\,c)}}.
\endaligned
\]

Precisely, the assumption of
{\em local real analyticity} of the $\varphi_j$ means that
there exists a positive number:
\[
\rho_0
>
0
\]
for which the power series expansion:
\[
\aligned
\varphi_j
\big(z_\smallbullet,\overline{z}_\smallbullet,u_\smallbullet\big)
&
=
\sum_{\alpha_\smallbullet\in\N^n}\,
\sum_{\beta_\smallbullet\in\N^n}\,
\sum_{\gamma_\smallbullet\in\N^c}\,
\underbrace{\varphi_{j,\alpha_\smallbullet,\beta_\smallbullet,\gamma_\smallbullet}
}_{\in\,\R}
\big(z_\smallbullet\big)^{\alpha_\smallbullet}\,
\big(\overline{z}_\smallbullet\big)^{\beta_\smallbullet}\,
\big(u_\smallbullet\big)^{\gamma_\smallbullet}
\\
&
\ \ \ \ \ \ \ \ \ \ \ \ \ \ \ \ \ \ \ \ \ \ \ \ \ \
{\scriptstyle{(1\,\leqslant\,j\,\leqslant\,c)}},
\endaligned
\]
converges when:
\[
\vert z_\smallbullet\vert
<
\rho_0,
\ \ \ \ \ \ \ \ \ \ \ \ \ \ \ \ \ \ \ \ \ \ \ \ \ \
\vert u_\smallbullet\vert
<
\rho_0,
\]
in the sense that its coefficients satisfy
a Cauchy-type estimate:
\[
\big\vert
\varphi_{j,\alpha_\smallbullet,\beta_\smallbullet,\gamma\smallbullet}
\big\vert
\,\leqslant\,
\constant\,
\bigg(
\frac{1}{\rho_0}
\bigg)^{\vert\alpha_\smallbullet\vert
+
\vert\beta_\smallbullet\vert
+
\vert\gamma_\smallbullet\vert}
\ \ \ \ \ \ \ \ \ \ \ \ \
{\scriptstyle{(1\,\leqslant\,j\,\leqslant\,c)}},
\]
for some positive $\constant > 0$.

\medskip\noindent{\bf Convention.}
Given a (formal or convergent) 
power series in $\NNN$ variables $({\sf t}_1, \dots, {\sf
t}_\NN)$ having complex coefficients:
\[
\Phi\big({\sf t}_1,\dots,{\sf t}_\NN\big)
=
\sum_{\gamma_1=0}^\infty\,
\cdots\,
\sum_{\gamma_\NN=0}^\infty\,
\underbrace{\Phi_{\gamma_1,\dots,\gamma_\NN}}_{\in\,\C}\,
\big({\sf t}_1\big)^{\gamma_1}\,
\cdots\,
\big({\sf t}_\NN\big)^{\gamma_\NN},
\]
its {\sl conjugate series}:
\[
\overline{\Phi}
\big({\sf t}_1,\dots,{\sf t}_\NN\big)
:=
\sum_{\gamma_1=0}^\infty\,
\cdots\,
\sum_{\gamma_\NN=0}^\infty\,
\overline{\Phi_{\gamma_1,\dots,\gamma_\NN}}
\big({\sf t}_1\big)^{\gamma_1}\,
\cdots\,
\big({\sf t}_\NN\big)^{\gamma_\NN}\,
\]
is defined by conjugating only its coefficients, so that
an overall conjugation operator:
\[
\overline{\Phi({\sf t}_1,\dots,{\sf t}_\NN)}
\equiv
\overline{\Phi}
\big(\overline{\sf t}_1,\dots,\overline{\sf t}_\NN\big)
\]
distributes {\em simultaneously on the series symbol and on the
variables}.

\medskip\noindent{\bf Notation.}
The {\sl complexifications} of the antiholomorphic variables:
\[
\big(\overline{z}_\smallbullet,\overline{w}_\smallbullet\big)
\]
will be denoted using just some underlining:
\[
\big(\underline{z}_\smallbullet,\underline{w}_\smallbullet\big),
\]
and these are $(n+c)$ new variables totally independent of $(z_\smallbullet, 
w_\smallbullet)$.

\medskip

If one also introduce new complex variables:
\[
\nu_1,\dots,\nu_c
\,\in\,\C,
\]
with:
\[
\Re\,\nu_1
=
u_1,
\,\dots\dots,\,
\Re\,\nu_c
=
u_c,
\]
the {\sl complexification} of $\varphi_\smallbullet$ is the power series:
\[
\varphi_\smallbullet\big(z_\smallbullet,\underline{z}_\smallbullet,\nu_\smallbullet\big)
:=
\sum_{\alpha_\smallbullet\in\N^n}\,
\sum_{\beta_\smallbullet\in\N^n}\,
\sum_{\gamma_\smallbullet\in\N^c}\,
\underbrace{\varphi_{\smallbullet,\alpha_\smallbullet,\beta_\smallbullet,\gamma_\smallbullet}
}_{\in\,\R}
\big(z_\smallbullet\big)^{\alpha_\smallbullet}\,
\big(\underline{z}_\smallbullet\big)^{\beta_\smallbullet}\,
\big(\nu_\smallbullet\big)^{\gamma_\smallbullet}
\]
in the $2n+c$ {\em complex} variables:
\[
\big(z_\smallbullet,\underline{z}_\smallbullet,\nu_\smallbullet\big),
\]
and it yet converges when
\[
\vert z_\smallbullet\vert
<
\rho_0,
\ \ \ \ \ \ \ \ \ \ \ \ \ \ \ \ \ \ \ \ \ \ \ \ \ \
\vert \underline{z}_\smallbullet\vert
<
\rho_0,
\ \ \ \ \ \ \ \ \ \ \ \ \ \ \ \ \ \ \ \ \ \ \ \ \ \
\vert \nu_\smallbullet\vert
<
\rho_0,
\]
because:
\[
\big\vert
\varphi_{\smallbullet,\alpha_\smallbullet,\beta_\smallbullet,\gamma\smallbullet}
\big\vert
\,\leqslant\,
\constant\,
\bigg(
\frac{1}{\rho_0}
\bigg)^{\vert\alpha_\smallbullet\vert
+
\vert\beta_\smallbullet\vert
+
\vert\gamma_\smallbullet\vert}.
\]

\medskip\noindent{\bf Complexification-Identity Principle.}
{\em Two converging power series:}
\[
\aligned
G\big(z_\smallbullet,w_\smallbullet,\underline{z}_\smallbullet,
\underline{w}_\smallbullet
\big)
&
=
\sum_{\alpha_\smallbullet\in\N^n}\,
\sum_{\beta_\smallbullet\in\N^c}\,
\sum_{\gamma_\smallbullet\in\N^n}\,
\sum_{\delta_\smallbullet\in\N^c}\,
\underbrace{G_{\alpha_\smallbullet,
\beta_\smallbullet,\gamma_\smallbullet,
\delta_\smallbullet}}_{\in\,\C}\,
\big(z_\smallbullet\big)^{\alpha_\smallbullet}\,
\big(w_\smallbullet\big)^{\beta_\smallbullet}\,
\big(\underline{z}_\smallbullet\big)^{\gamma_\smallbullet}\,
\big(\underline{w}_\smallbullet\big)^{\delta_\smallbullet}
\\
H\big(z_\smallbullet,w_\smallbullet,\underline{z}_\smallbullet,
\underline{w}_\smallbullet
\big)
&
=
\sum_{\alpha_\smallbullet\in\N^n}\,
\sum_{\beta_\smallbullet\in\N^c}\,
\sum_{\gamma_\smallbullet\in\N^n}\,
\sum_{\delta_\smallbullet\in\N^c}\,
\underbrace{H_{\alpha_\smallbullet,
\beta_\smallbullet,\gamma_\smallbullet,
\delta_\smallbullet}}_{\in\,\C}\,
\big(z_\smallbullet\big)^{\alpha_\smallbullet}\,
\big(w_\smallbullet\big)^{\beta_\smallbullet}\,
\big(\underline{z}_\smallbullet\big)^{\gamma_\smallbullet}\,
\big(\underline{w}_\smallbullet\big)^{\delta_\smallbullet}
\endaligned
\]
{\em in the $2n + 2c$ independent complex variables:}
\[
\big(z_\smallbullet,w_\smallbullet, 
\underline{z}_\smallbullet
\underline{w}_\smallbullet
\big)
\,\in\,\C^{n+c+n+c}
\]
{\em are identically equal:}
\[
G_{\alpha_\smallbullet,
\beta_\smallbullet,\gamma_\smallbullet,
\delta_\smallbullet}
=
H_{\alpha_\smallbullet,
\beta_\smallbullet,\gamma_\smallbullet,
\delta_\smallbullet}
\ \ \ \ \ \ \ \ \ \ \ \ \
{\scriptstyle{(\forall\,\alpha_\smallbullet,\,
\beta_\smallbullet,\,\gamma_\smallbullet,\,\delta_\smallbullet)}}
\]
{\em if and only if their restrictions to the antiholomorphic
diagonal:}
\[
\overline{\Lambda}
:=
\big\{
\big(z_\smallbullet,w_\smallbullet,\underline{z}_\smallbullet,
\underline{w}_\smallbullet\big)
\in\C^{n+c+n+c}\,
\colon\,
\underline{z}_\smallbullet
=
\overline{z}_\smallbullet,\,\,
\underline{w}_\smallbullet
=
\overline{w}_\smallbullet
\big\},
\]
{\em coincide as functions:}
\[
G\big(
z_\smallbullet,w_\smallbullet,
\overline{z}_\smallbullet,
\overline{w}_\smallbullet
\big)
\,\equiv\,
H\big(
z_\smallbullet,
w_\smallbullet,
\overline{z}_\smallbullet,
\overline{w}_\smallbullet
\big),
\]
{\em that is more precisely, if and only if:}
\[
\aligned
&
G\Big(
x_\smallbullet+\isqrt\,y_\smallbullet,\,
u_\smallbullet+\isqrt\,v_\smallbullet,\,
x_\smallbullet-\isqrt\,y_\smallbullet,\,
u_\smallbullet-\isqrt\,v_\smallbullet
\Big)
\,\equiv
\\
&
\ \ \ \ \ \ \ \ \ \ \ \ \
\equiv\,
H\Big(
x_\smallbullet+\isqrt\,y_\smallbullet,\,
u_\smallbullet+\isqrt\,v_\smallbullet,\,
x_\smallbullet-\isqrt\,y_\smallbullet,\,
u_\smallbullet-\isqrt\,v_\smallbullet
\Big)
\endaligned
\]
{\em identically in:}
\[
\C\big\{x_\smallbullet,y_\smallbullet,u_\smallbullet,v_\smallbullet\big\},
\]
{\em {\em i.e.} as functions of:}
\[
\big(x_\smallbullet,y_\smallbullet,u_\smallbullet,v_\smallbullet\big)
\,\in\,
\R^{n+n+c+c}.
\]

\proof
One direction being trivial, starting then after subtraction from:
\[
0
\,\equiv\,
\sum_{\alpha_\smallbullet\in\N^n}\,
\sum_{\beta_\smallbullet\in\N^c}\,
\sum_{\gamma_\smallbullet\in\N^n}\,
\sum_{\delta_\smallbullet\in\N^c}\,
\Big[
G_{\alpha_\smallbullet,\beta_\smallbullet,
\gamma_\smallbullet,\delta_\smallbullet}
\,-\,
H_{\alpha_\smallbullet,\beta_\smallbullet,
\gamma_\smallbullet,\delta_\smallbullet}
\Big]\,
\big(z_\smallbullet\big)^{\alpha_\smallbullet}\,
\big(w_\smallbullet\big)^{\beta_\smallbullet}\,
\big(\overline{z}_\smallbullet\big)^{\gamma_\smallbullet}\,
\big(\overline{w}_\smallbullet\big)^{\delta_\smallbullet},
\]
one easily convinces oneself by suitable induction on:
\[
\big(
\alpha_\smallbullet,\beta_\smallbullet,
\gamma_\smallbullet,\delta_\smallbullet
\big)
\,\in\,
\N^{n+c+n+c}
\]
that all Taylor coefficients match up: 
\[
G_{\alpha_\smallbullet,\beta_\smallbullet,
\gamma_\smallbullet,\delta_\smallbullet}
=
H_{\alpha_\smallbullet,\beta_\smallbullet,
\gamma_\smallbullet,\delta_\smallbullet}
\ \ \ \ \ \ \ \ \ \ \ \ \
{\scriptstyle{(\forall\,\alpha_\smallbullet,\,
\beta_\smallbullet,\,\gamma_\smallbullet,\,\delta_\smallbullet)}},
\]
by successively differentiating
this identity and by setting in it:
\[
z_\smallbullet 
= 
w_\smallbullet
= 
0, 
\]
taking account of the elementary relations:
\[
\aligned
\delta_{k_1,k_2}
&
=
\frac{\partial}{\partial z_{k_1}}
\big(z_{k_2}\big)
=
\frac{\partial}{\partial\overline{z}_{k_1}}
\big(\overline{z}_{k_2}\big),
\\
\delta_{l_1,l_2}
&
=
\frac{\partial}{\partial w_{l_1}}
\big(w_{l_2}\big)
=
\frac{\partial}{\partial\overline{w}_{l_1}}
\big(\overline{w}_{l_2}\big),
\\
0
&
=
\frac{\partial}{\partial z_k}
\big(w_l\big)
=
\frac{\partial}{\partial\overline{z}_k}
\big(w_l\big),
\\
0
&
=
\frac{\partial}{\partial w_l}
\big(z_k\big)
=
\frac{\partial}{\partial\overline{w}_l}
\big(z_k\big),
\endaligned
\]
term by term differentiation being justified by standard
normal convergence arguments.
\endproof

\noindent{\bf Reality feature.}
{\em The $c$ scalar functions $\varphi_j$ being real-valued:}
\[
\overline{\varphi_j\big(z_\smallbullet,\overline{z}_\smallbullet,
u_\smallbullet\big)}
\,\equiv\,
\varphi_j\big(z_\smallbullet,
\overline{z}_\smallbullet,u_\smallbullet\big)
\ \ \ \ \ \ \ \ \ \ \ \ \
{\scriptstyle{(j\,=\,1\,\cdots\,c)}},
\]
{\em their power series complex coefficients satisfy:}
\[
\overline{\varphi_{j,\alpha_\smallbullet,\beta_\smallbullet,\gamma_\smallbullet}}
\,=\,
\varphi_{j,\beta_\smallbullet,\alpha_\smallbullet,\gamma_\smallbullet},
\]
{\em and conversely.}

\proof
Conjugating term by term the power series expansion 
is justified by normal convergence:
\[
\aligned
\overline{\varphi_j\big(z_\smallbullet,\overline{z}_\smallbullet,
u_\smallbullet\big)}
&
=
\overline{
\sum_{\alpha_\smallbullet\in\N^n}\,
\sum_{\beta_\smallbullet\in\N^n}\,
\sum_{\gamma_\smallbullet\in\N^c}\,
\varphi_{j,\alpha_\smallbullet,\beta_\smallbullet,\gamma_\smallbullet}\,
\big(z_\smallbullet\big)^{\alpha_\smallbullet}\,
\big(\overline{z}_\smallbullet\big)^{\beta_\smallbullet}\,
\big(u_\smallbullet\big)^{\gamma_\smallbullet}
}
\\
&
=
\sum_{\alpha_\smallbullet\in\N^n}\,
\sum_{\beta_\smallbullet\in\N^n}\,
\sum_{\gamma_\smallbullet\in\N^c}\,
\overline{\varphi_{j,\alpha_\smallbullet,\beta_\smallbullet,
\gamma_\smallbullet}}\,
\big(\overline{z}_\smallbullet\big)^{\alpha_\smallbullet}\,
\big(z_\smallbullet\big)^{\beta_\smallbullet}\,
\big(u_\smallbullet\big)^{\gamma_\smallbullet}
\\
&
=
\varphi_j\big(z_\smallbullet,
\overline{z}_\smallbullet,u_\smallbullet\big)
\\
&
=
\sum_{\alpha_\smallbullet\in\N^n}\,
\sum_{\beta_\smallbullet\in\N^n}\,
\sum_{\gamma_\smallbullet\in\N^c}\,
\varphi_{j,\alpha_\smallbullet,\beta_\smallbullet,\gamma_\smallbullet}\,
\big(z_\smallbullet\big)^{\alpha_\smallbullet}\,
\big(\overline{z}_\smallbullet\big)^{\beta_\smallbullet}\,
\big(u_\smallbullet\big)^{\gamma_\smallbullet},
\endaligned
\]
and the assumed reality yields the stated coefficient equalities,
after reorganizing the index summation:
\[
\alpha_\smallbullet
\,\longleftrightarrow\, 
\beta_\smallbullet
\]
for comparison-identification thanks to the complexification-identity
principle.
\endproof

\noindent{\bf Scholium.}
{\em For $\nu_\smallbullet \in \C^c$, on has:}
\[
\overline{\varphi_j\big(z_\smallbullet,\overline{z}_\smallbullet,
\nu_\smallbullet\big)}
\,\equiv\,
\varphi_j
\big(z_\smallbullet,\overline{z}_\smallbullet,\overline{\nu}_\smallbullet\big).
\]

\proof
At the level of converging power series, this follows (exercise)
from the above symmetry relations between coefficients.
\endproof

Now, writing the $c$ scalar equations of $M$ as:
\[
\frac{w_\smallbullet-\overline{w}_\smallbullet}{2\,\isqrt}
=
\varphi_\smallbullet
\bigg(
z_\smallbullet,\,
\overline{z}_\smallbullet,\,
\frac{w_\smallbullet+\overline{w}_\smallbullet}{2}
\bigg),
\]
and reminding:
\[
\varphi_\smallbullet
=
{\rm O}(2),
\]
the implicit function theorem enables one to solve for either $w_\smallbullet$ 
or $\overline{ w}_\smallbullet$, but to be appropriately rigorous,
it is better in fact to {\sl complexify} 
the equation that one wants to solve, namely to
look at the $c$ scalar equations:
\[
\frac{w_\smallbullet-\underline{w}_\smallbullet}{2\,\isqrt}
=
\varphi_\smallbullet
\bigg(
z_\smallbullet,\,
\underline{z}_\smallbullet,\,
\frac{w_\smallbullet+\underline{w}_\smallbullet}{2}
\bigg),
\]
that are now {\em holomorphic} with respect to the $2n + 2c$ variables:
\[
\big(z_\smallbullet,w_\smallbullet,
\underline{z}_\smallbullet,\underline{w}_\smallbullet\big)
\,\in\,\C^{n+c+n+c}.
\]

Using the analytic implicit function theorem, two options
indeed present themselves:

\medskip\noindent$\square$\,
either solve with respect to $w_\smallbullet$;

\medskip\noindent$\square$\,
or solve with respect to $\underline{ w}_\smallbullet$.

\medskip
The first solving provides:
\[
w_\smallbullet
=
\Theta_\smallbullet
\big(z_\smallbullet,\underline{z}_\smallbullet,\underline{w}_\smallbullet\big),
\]
for some certain $c$ {\em complex-analytic} (holomorphic) functions:
\[
\aligned
\Theta_j
\big(z_\smallbullet,\underline{z}_\smallbullet,\underline{w}_\smallbullet\big)
&
=
\sum_{\alpha_\smallbullet\in\N^n}\,
\sum_{\beta_\smallbullet\in\N^n}\,
\sum_{\gamma_\smallbullet\in\N^c}\,
\Theta_{j,\alpha_\smallbullet,\beta_\smallbullet,\gamma_\smallbullet}
\big(z_\smallbullet\big)^{\alpha_\smallbullet}\,
\big(\underline{z}_\smallbullet\big)^{\beta_\smallbullet}\,
\big(\underline{w}_\smallbullet\big)^{\gamma_\smallbullet}
\\
&
\ \ \ \ \ \ \ \ \ \ \ \ \ \ \ \ \ \ \ \ \ \ \ \ \ \
{\scriptstyle{(1\,\leqslant\,j\,\leqslant\,c)}},
\endaligned
\]
these $c$ power series being also convergent,
namely having coefficients which also satisfy\,\,---\,\,shrinking
$\rho_0 > 0$ if necessary\,\,---\,\,a Cauchy-type estimate:
\[
\big\vert
\Theta_{j,\alpha_\smallbullet,\beta_\smallbullet,\gamma_\smallbullet}
\big\vert
\,\leqslant\,
\constant\,
\bigg(
\frac{1}{\rho_0}
\bigg)^{
\vert\alpha_\smallbullet\vert+\vert\beta_\smallbullet\vert
+\vert\gamma_\smallbullet\vert}
\ \ \ \ \ \ \ \ \ \ \ \ \ \ \ 
{\scriptstyle{(1\,\leqslant\,j\,\leqslant\,c)}}.
\]

By definition, the solution $\Theta_\smallbullet$ satisfies:
\[
\frac{\Theta_\smallbullet\big(z_\smallbullet,\underline{z}_\smallbullet,
\underline{w}_\smallbullet\big)-\underline{w}_\smallbullet}{
2\,\isqrt}
\,\,\equiv\,\,
\varphi_\smallbullet
\bigg(
z_\smallbullet,\,
\underline{z}_\smallbullet,\,\,
\frac{\Theta_\smallbullet\big(z_\smallbullet,\underline{z}_\smallbullet,
\underline{w}_\smallbullet\big)+\underline{w}_\smallbullet}{2}
\bigg),
\]
identically in the ring of $c$ converging power series:
\[
\C\big\{z_\smallbullet,\underline{z}_\smallbullet,\underline{w}_\smallbullet\big\}^c.
\]

On restriction to the antiholomorphic diagonal:
\[
\big\{
\underline{z}_\smallbullet
=
\overline{z}_\smallbullet,\,
\underline{w}_\smallbullet
=
\overline{w}_\smallbullet
\big\},
\]
one then obtains:
\[
\frac{\Theta_\smallbullet\big(z_\smallbullet,\overline{z}_\smallbullet,
\overline{w}_\smallbullet\big)-\overline{w}_\smallbullet}{
2\,\isqrt}
\,\,\equiv\,\,
\varphi_\smallbullet
\bigg(
z_\smallbullet,\,
\overline{z}_\smallbullet,\,\,
\frac{\Theta_\smallbullet\big(z_\smallbullet,\overline{z}_\smallbullet,
\overline{w}_\smallbullet\big)+\overline{w}_\smallbullet}{2}
\bigg),
\]
so that the $c$ complex scalar equations:
\[
w_\smallbullet
=
\Theta_\smallbullet
\big(z_\smallbullet,\overline{z}_\smallbullet,\overline{w}_\smallbullet\big)
\]
also constitute defining equations for $M^{ 2n + c}
\subset \C^{ n+c}$, {\em the $c$ real part
equations and the $c$ imaginary
part equations being equivalent} (exercise of understanding,
or {\em see} below).

Furthermore, using:
\[
\overline{\varphi_\smallbullet
\big(z_\smallbullet,\overline{z}_\smallbullet,\nu_\smallbullet\big)}
\,\equiv\,
\varphi_\smallbullet
\big(z_\smallbullet,\overline{z}_\smallbullet,\overline{\nu}_\smallbullet\big),
\]
a plain conjugation of what precedes yields:
\[
\frac{-\,\overline{\Theta}_\smallbullet\big(
\overline{z}_\smallbullet,z_\smallbullet,w_\smallbullet\big)
+w_\smallbullet}{2\,\isqrt}
\,\equiv\,
\varphi_\smallbullet
\bigg(
z_\smallbullet,\,
\overline{z}_\smallbullet,\,\,
\frac{
\overline{\Theta}_\smallbullet\big(\overline{z}_\smallbullet,z_\smallbullet,
w_\smallbullet\big)+w_\smallbullet}{2}
\bigg),
\]
and a last complexification yields:
\[
\frac{-\,\overline{\Theta}_\smallbullet\big(
\underline{z}_\smallbullet,z_\smallbullet,w_\smallbullet\big)
+w_\smallbullet}{2\,\isqrt}
\,\equiv\,
\varphi_\smallbullet
\bigg(
z_\smallbullet,\,
\underline{z}_\smallbullet,\,\,
\frac{
\overline{\Theta}_\smallbullet\big(\underline{z}_\smallbullet,z_\smallbullet,
w_\smallbullet\big)+w_\smallbullet}{2}
\bigg).
\]

On the other hand, if one solves secondly 
for $\underline{ w}_\smallbullet$ the same $c$ equations:
\[
\frac{w_\smallbullet-\underline{w}_\smallbullet}{2\,\isqrt}
\,=\,
\varphi_\smallbullet
\bigg(
z_\smallbullet,\underline{z}_\smallbullet,
\frac{w_\smallbullet+\underline{w}_\smallbullet}{2}
\bigg),
\]
getting:
\[
\underline{w}_\smallbullet
=
\widetilde{\Theta}_\smallbullet
\big(
\underline{z}_\smallbullet,z_\smallbullet,w_\smallbullet
\big),
\]
for a certain $\C^c$-valued holomorphic map $\widetilde{ \Theta}_\smallbullet$
satisfying by definition:
\[
\frac{w_\smallbullet
-
\widetilde{\Theta}_\smallbullet\big(\underline{z}_\smallbullet,
z_\smallbullet,w_\smallbullet\big)}{2\,\isqrt}
\,\equiv\,
\varphi_\smallbullet
\bigg(
z_\smallbullet,\underline{z}_\smallbullet,\,
\frac{w_\smallbullet+\widetilde{\Theta}_\smallbullet
\big(\underline{z}_\smallbullet,z_\smallbullet,w_\smallbullet\big)}{2}
\bigg),
\]
then because this $\widetilde{ \Theta}_\smallbullet$ provided by the analytic
implicit function theorem is {\em unique}, a comparison
with what has been written at the moment yields
the coincidence:
\[
\widetilde{\Theta}_\smallbullet
=
\overline{\Theta}_\smallbullet.
\]
In summary:

\medskip\noindent{\bf Proposition.}
{\em
Solving either with respect to $w_\smallbullet$ or to $\underline{ w}_\smallbullet$
some $c$ real analytic local defining equations:}
\[
\frac{w_\smallbullet-\overline{w}_\smallbullet}{2\,\isqrt}
=
v_\smallbullet
=
\varphi_\smallbullet
\big(z_\smallbullet,\overline{z}_\smallbullet,u_\smallbullet\big)
=
\varphi_\smallbullet
\bigg(
z_\smallbullet,\overline{z}_\smallbullet,\,
\frac{w_\smallbullet+\overline{w}_\smallbullet}{2}
\bigg)
\]
{\em for a $\mathcal{ C}^\omega$ CR-generic submanifold
$M^{ 2n + c} \subset \C^{ n+c}$, there exists a
single $\C^c$-valued local holomorphic function:}
\[
\Theta_\smallbullet
=
\Theta_\smallbullet\big(z_\smallbullet,\overline{z}_\smallbullet,\overline{w}_\smallbullet\big)
\]
{\em such that the two solutions are conjugate to each other:}
\[
\aligned
w_\smallbullet
&
=
\Theta_\smallbullet
\big(z_\smallbullet,\overline{z}_\smallbullet,\overline{w}_\smallbullet\big),
\\
\overline{w}_\smallbullet
&
=
\overline{\Theta}_\smallbullet
\big(\overline{z}_\smallbullet,z_\smallbullet,w_\smallbullet\big),
\endaligned
\]
{\em and a point of coordinates:}
\[
\big(z_\smallbullet,w_\smallbullet\big)
\,\in\,M
\]
{\em belongs to $M$ if and only if it satisfies either
one of these three systems of $c$ equations (hence the
other two as well):}
\[
\aligned
v_\smallbullet
&
=
\varphi_\smallbullet
\big(z_\smallbullet,\overline{z}_\smallbullet,u_\smallbullet\big),
\\
w_\smallbullet
&
=
\Theta_\smallbullet
\big(z_\smallbullet,\overline{z}_\smallbullet,\overline{w}_\smallbullet\big),
\\
\overline{w}_\smallbullet
&
=
\overline{\Theta}_\smallbullet
\big(\overline{z}_\smallbullet,z_\smallbullet,w_\smallbullet\big).
\qed
\endaligned
\]

\medskip
Passing to complexified variables:
\[
\big(z_\smallbullet,w_\smallbullet,\underline{z}_\smallbullet,
\underline{w}_\smallbullet\big)
\,\in\,
\C^{n+c+n+c},
\]
the three systems of $c$ scalar equations:
\[
\aligned
{\textstyle{\frac{w_\smallbullet-\underline{w}_\smallbullet}{2\,\isqrt}}}
&
=
\varphi_\smallbullet
\big(z_\smallbullet,\underline{z}_\smallbullet,
{\textstyle{\frac{w_\smallbullet+\underline{w}_\smallbullet}{2}}}
\big),
\\
w_\smallbullet
&
=
\Theta_\smallbullet
\big(z_\smallbullet,\underline{z}_\smallbullet,\underline{w}_\smallbullet\big),
\\
\underline{w}_\smallbullet
&
=
\overline{\Theta}_\smallbullet
\big(\underline{z}_\smallbullet,z_\smallbullet,w_\smallbullet\big)
\endaligned
\]
are therefore all equivalent by pairs, and since their differentials
at the origin are all of maximal rank equal to $c$, 
if follows from the so-called {\sl Hadamard lemma} 
that there exist two $c \times c$ invertible matrices of 
converging power series:
\[
\aligned
{\tt b}_{\smallbullet,\smallbullet}
\big(
z_\smallbullet,w_\smallbullet,\underline{z}_\smallbullet,\underline{w}_\smallbullet
\big)
&
\,\in\,
{\sf GL}_{c\times c}
\Big(
\C\big\{
z_\smallbullet,w_\smallbullet,\underline{z}_\smallbullet,\underline{w}_\smallbullet
\big\}
\Big),
\\
{\tt a}_{\smallbullet,\smallbullet}
\big(
z_\smallbullet,w_\smallbullet,\underline{z}_\smallbullet,\underline{w}_\smallbullet
\big)
&
\,\in\,
{\sf GL}_{c\times c}
\Big(
\C\big\{
z_\smallbullet,w_\smallbullet,\underline{z}_\smallbullet,\underline{w}_\smallbullet
\big\}
\Big),
\endaligned
\]
such that:
\[
\aligned
{\textstyle{\frac{w_\smallbullet-\underline{w}_\smallbullet}{2\,\isqrt}}}
-
\varphi_\smallbullet
\big(z_\smallbullet,\underline{z}_\smallbullet,
{\textstyle{\frac{w_\smallbullet+\underline{w}_\smallbullet}{2}}}
\big)
&
\,\equiv\,
{\tt b}_{\smallbullet,\smallbullet}\,
\big[
w_\smallbullet
-
\Theta_\smallbullet
\big(z_\smallbullet,\underline{z}_\smallbullet,\underline{w}_\smallbullet\big)
\big],
\\
\underline{w}_\smallbullet
-
\overline{\Theta}_\smallbullet
\big(\underline{z}_\smallbullet,z_\smallbullet,w_\smallbullet\big)
&
\,\equiv\,
{\tt a}_{\smallbullet,\smallbullet}\,
\big[
w_\smallbullet
-
\Theta_\smallbullet
\big(z_\smallbullet,\underline{z}_\smallbullet,\underline{w}_\smallbullet\big)
\big].
\endaligned
\]

Next, the {\sl extrinsic complexification} of $M$ is
the complex submanifold:
\[
M^{e_c}
\,\subset\,
\C^{n+c+n+c}
\]
of complex codimension $c$ which is 
defined by either one of the two (equivalent (mental exercise) 
systems of $c$ equations:
\[
\aligned
w_\smallbullet
&
=
\Theta_\smallbullet
\big(z_\smallbullet,\underline{z}_\smallbullet,\underline{w}_\smallbullet\big),
\\
\underline{w}_\smallbullet
&
=
\overline{\Theta}_\smallbullet
\big(\underline{z}_\smallbullet,z_\smallbullet,w_\smallbullet\big).
\qed
\endaligned
\]
It follows (exercise of understanding) that the two functional
equations:
\[
\boxed{\,\,
\aligned
w_\smallbullet
&
\equiv
\Theta_\smallbullet
\big(z_\smallbullet,\underline{z}_\smallbullet,\,
\overline{\Theta}_\smallbullet(\underline{z}_\smallbullet,z_\smallbullet,
w_\smallbullet)\big),\,\,
\\
\underline{w}_\smallbullet
&
\equiv
\overline{\Theta}_\smallbullet
\big(\underline{z}_\smallbullet,z_\smallbullet,
\Theta_\smallbullet(z_\smallbullet,\underline{z}_\smallbullet,
\underline{w}_\smallbullet)\big),\,\,
\endaligned}
\]
are identically satisfied, respectively, in:
\[
\aligned
&
\C\big\{\underline{z}_\smallbullet,z_\smallbullet,w_\smallbullet\big\}^c,
\\
&
\C\big\{z_\smallbullet,\underline{z}_\smallbullet,\underline{w}_\smallbullet\big\}^c.
\endaligned
\]

\medskip

Consider now a local biholomorphic map:
\[
\aligned
\big(z_\smallbullet,w_\smallbullet\big)
&
\,\longmapsto\,
\big(
z_\smallbullet'(z_\smallbullet,w_\smallbullet),\,
w_\smallbullet'(z_\smallbullet,w_\smallbullet)
\big)
\\
&\ \ \,
=:
\big(z_\smallbullet',w_\smallbullet'\big),
\endaligned
\]
which, to fix ideas, 
sends $(0_\smallbullet, 0_\smallbullet)$ to
$(0_\smallbullet', 0_\smallbullet')$, and moreover\,\,---\,\,after
a possible renumbering of coordinates\,\,---, which sends
$M$ onto another CR-generic:
\[
{M'}^{2n+c}
\,\subset\,
{\C'}^{n+c}
\]
that is also representable under similar equivalent graphed forms:
\[
\aligned
v_\smallbullet'
&
=
\varphi_\smallbullet'
\big(z_\smallbullet',\overline{z}_\smallbullet',u_\smallbullet'\big),
\\
w_\smallbullet'
&
=
\Theta_\smallbullet'
\big(z_\smallbullet',\overline{z}_\smallbullet',\overline{w}_\smallbullet'\big),
\\
\overline{w}_\smallbullet'
&
=
\overline{\Theta}_\smallbullet'
\big(\overline{z}_\smallbullet',z_\smallbullet',w_\smallbullet'\big).
\endaligned
\]

Since the differentials of each among these three collections
of $c$ equations are of rank $c$ at the origin (hence
at every nearby point), the so-called Hadamard lemma
provides a $c \times c$ invertible matrix of
converging power series:
\[
{\tt c}_{\smallbullet,\smallbullet}
\big(
z_\smallbullet,w_\smallbullet,\underline{z}_\smallbullet,\underline{w}_\smallbullet
\big)
\,\in\,
{\sf GL}_{c\times c}
\Big(
\C\big\{
z_\smallbullet,w_\smallbullet,\underline{z}_\smallbullet,\underline{w}_\smallbullet
\big\}
\Big),
\]
such that one has:
\[
\aligned
&
w_\smallbullet'\big(z_\smallbullet,w_\smallbullet\big)
-
\Theta_\smallbullet'
\big(
z_\smallbullet'
(z_\smallbullet,w_\smallbullet),\,
\overline{z}_\smallbullet'
(\underline{z}_\smallbullet,\underline{w}_\smallbullet),\,
\overline{w}_\smallbullet'
(\underline{z}_\smallbullet,\underline{w}_\smallbullet)
\big)
\,\equiv
\\
&
\ \ \ \ \ \ \ \ \ \ \ 
\equiv\,
{\tt c}_{\smallbullet,\smallbullet}
\big(z_\smallbullet,w_\smallbullet,\underline{z}_\smallbullet,\underline{w}_\smallbullet\big)\,
\big[
w_\smallbullet
-
\Theta_\smallbullet
\big(z_\smallbullet,\underline{z}_\smallbullet,\underline{w}_\smallbullet\big)
\big],
\endaligned
\]
identically in:
\[
\C\big\{
z_\smallbullet,w_\smallbullet,\underline{z}_\smallbullet,\underline{w}_\smallbullet
\big\}^c.
\]

In other (interpretational) words, the {\sl complexified mapping:}
\[
\aligned
\big(z_\smallbullet,w_\smallbullet,\underline{z}_\smallbullet,\underline{w}_\smallbullet\big)
&
\,\longmapsto\,
\big(
z_\smallbullet'(z_\smallbullet,w_\smallbullet),\,
w_\smallbullet'(z_\smallbullet,w_\smallbullet),\,
\overline{z}_\smallbullet'(\underline{z}_\smallbullet,\underline{w}_\smallbullet),\,
\overline{w}_\smallbullet'(\underline{z}_\smallbullet,\underline{w}_\smallbullet)
\big)
\\
&\ \ \,
=:
\big(z_\smallbullet',w_\smallbullet',\underline{z}_\smallbullet',
\underline{w}_\smallbullet'\big),
\endaligned
\]
sends the complexification:
\[
M^{e_c}
=
\big\{
(z_\smallbullet,w_\smallbullet,\underline{z}_\smallbullet,\underline{w}_\smallbullet)
\colon\,
w_\smallbullet
=
\Theta_\smallbullet
\big(z_\smallbullet,\underline{z}_\smallbullet,\underline{w}_\smallbullet\big)
\big\}
\]
to the complexification:
\[
{M'}^{e_c}
=
\big\{
(z_\smallbullet',w_\smallbullet',\underline{z}_\smallbullet',\underline{w}_\smallbullet')
\colon\,
w_\smallbullet'
=
\Theta_\smallbullet'
\big(z_\smallbullet',\underline{z}_\smallbullet',\underline{w}_\smallbullet'\big)
\big\}.
\]

\medskip\noindent{\bf Proposition.}
{\em After a local biholomorphism fixing the origin of the form:}
\[
\aligned
z_k
&
\,\longmapsto\,
z_k
\ \ \ \ \ \ \ \ \ \ \ \ \ \ \ \ \ \ \ \ \ \ \ \ \ \ \ \ \ \ \ 
{\scriptstyle{(1\,\leqslant\,k\,\leqslant\,n)}},
\\
w_l
&
\,\longmapsto\,
w_l\big(z_\smallbullet,w_\smallbullet\big)
\ \ \ \ \ \ \ \ \ \ \ \ \ \ \ \ \ \ 
{\scriptstyle{(1\,\leqslant\,l\,\leqslant\,c)}},
\endaligned
\]
{\em the new real analytic local graphing functions satisfy:}
\[
\aligned
0
&
\equiv
\varphi_1\big(0,\overline{z}_\smallbullet,u_\smallbullet\big)
\equiv
\varphi_1\big(z_\smallbullet,0,u_\smallbullet\big),
\\
\cdot\,\cdot
&
\cdots\cdots\cdots\cdots\cdots\cdots\cdots\cdots\cdots\cdots
\\
0
&
\equiv
\varphi_c\big(0,\overline{z}_\smallbullet,u_\smallbullet\big)
\equiv
\varphi_c\big(z_\smallbullet,0,u_\smallbullet\big),
\endaligned
\]
{\em that is to say after notational contraction:}
\[
0
\equiv
\varphi_\smallbullet
\big(z_\smallbullet,0_\smallbullet,u_\smallbullet\big)
\equiv
\varphi_\smallbullet\big(0_\smallbullet,\overline{z}_\smallbullet,
u_\smallbullet),
\]
{\em and simultaneously also:}
\[
\aligned
\underline{w}_\smallbullet
&
\equiv
\Theta_\smallbullet\big(0_\smallbullet,\underline{z}_\smallbullet,
\underline{w}_\smallbullet\big)
\equiv\,
\Theta_\smallbullet\big(z_\smallbullet,0_\smallbullet,
\underline{w}_\smallbullet\big),
\\
w_\smallbullet
&
\equiv
\overline{\Theta}
\big(0_\smallbullet,z_\smallbullet,w_\smallbullet\big)\,\,\,
\equiv\,
\overline{\Theta}_\smallbullet\big(
\underline{z}_\smallbullet,0_\smallbullet,w_\smallbullet\big).
\endaligned
\]

\proof
A composition of {\em two} appropriate
local biholomorphisms will do the job.

As a {\sl Step I}, introduce the local holomorphic transformation defined by:
\[
\aligned
z_\smallbullet
&
=
z_\smallbullet',
\\
w_\smallbullet
&
=
w_\smallbullet'
+
\isqrt\,\varphi_\smallbullet\big(0_\smallbullet,0_\smallbullet,w_\smallbullet'\big).
\endaligned
\] 
It has an $(n + c) \times (n+c)$ Jacobian matrix at the
origin equal to the identity, since:
\[
\varphi_\smallbullet
=
{\rm O}(2)
\]
after the (assumed) preliminary affine normalization:
\[
T_0M
=
\big\{v_\smallbullet=0_\smallbullet\big\},
\]
hence this transformation is a local {\em biholomorphism}.

Through it, $M$ is tranformed to a certain CR-generic
submanifold:
\[
{M'}^{2n+c}
\,\subset\,
{\C'}^{n+c}
\]
still satisfying:
\[
T_{0'}M'
=
\big\{v_\smallbullet'=0_\smallbullet'\big\},
\]
so that the analytic implicit function theorem solves:
\[
v_\smallbullet'
=
\varphi_\smallbullet'
\big(z_\smallbullet',\overline{z}_\smallbullet',u_\smallbullet'\big).
\]

Visibly:
\[
\big\{
z_\smallbullet=0_\smallbullet
\big\}
\,=\,
\big\{
z_\smallbullet'=0_\smallbullet'
\big\},
\]
whence:
\[
\big\{
z_\smallbullet=0_\smallbullet
\big\}
\cap
M
\,=\,
\big\{
z_\smallbullet'=0_\smallbullet'
\big\}
\cap
M',
\]
that is to say:
\[
\big\{
\big(0_\smallbullet,\,u_\smallbullet
+
\isqrt\,\varphi_\smallbullet(0_\smallbullet,0_\smallbullet,u_\smallbullet)\big)
\big\}
=
\big\{
\big(0_\smallbullet',\,u_\smallbullet'
+
\isqrt\,\varphi_\smallbullet'(0_\smallbullet',0_\smallbullet',u_\smallbullet')\big)
\big\},
\]
this coincidence being through the related restriction of
the biholomorphism.

But through the biholomorphism in question, when:
\[
\aligned
w_\smallbullet'
&
=
u_\smallbullet'
+
\isqrt\,0_\smallbullet'
\\
&
=
u_\smallbullet',
\endaligned
\]
and when:
\[
z_\smallbullet'
=
0_\smallbullet',
\]
one visibly covers arbitrary points:
\[
\aligned
\big(z_\smallbullet,w_\smallbullet\big)
&
=
\big(0_\smallbullet,\,
u_\smallbullet'
+
\isqrt\,
\varphi_\smallbullet\big(0_\smallbullet,0_\smallbullet,u_\smallbullet'\big)
\big)
\\
&
\in\,
\big\{z_\smallbullet=0_\smallbullet\big\}
\cap
M,
\endaligned
\]
so that necessarily:
\[
\big\{z_\smallbullet'=0_\smallbullet'\big\}
\cap
M'
=
\big\{
\big(0_\smallbullet',u_\smallbullet'\big)
\big\},
\]
and one concludes that:
\[
\varphi_\smallbullet'
\big(
0_\smallbullet',0_\smallbullet',u_\smallbullet'
\big)
\,\equiv\,
0_\smallbullet'.
\]

Dropping the primes, one may therefore now assume in order
to pursue the proof that:
\[
\varphi_\smallbullet
\big(0_\smallbullet,0_\smallbullet,u_\smallbullet\big)
\equiv
0_\smallbullet.
\]
It follows (mental exercise) that:
\[
\aligned
\Theta_\smallbullet
\big(0_\smallbullet,0_\smallbullet,\overline{w}_\smallbullet\big)
&
\,\equiv\,
\overline{w}_\smallbullet,
\\
\overline{\Theta}_\smallbullet
\big(0_\smallbullet,0_\smallbullet,w_\smallbullet\big)
&
\,\equiv\,
w_\smallbullet.
\endaligned
\]

As a {\sl Step II}, introduce the (second) local biholomorphism:
\[
\aligned
z_\smallbullet'
&
=
z_\smallbullet,
\\
w_\smallbullet'
&
=
\overline{\Theta}_\smallbullet
\big(0_\smallbullet,z_\smallbullet,w_\smallbullet\big).
\endaligned
\]
It transforms $M$ to a certain generic submanifold:
\[
{M'}^{2n+c}
\,\subset\,
{\C'}^{n+c},
\]
still satisfying:
\[
T_{0'}M'
=
\big\{v_\smallbullet'=0_\smallbullet'\big\},
\]
whence $M'$ is defined by:
\[
\aligned
v_\smallbullet'
&
=
\varphi_\smallbullet'
\big(z_\smallbullet',\overline{z}_\smallbullet',u_\smallbullet'\big),
\\
w_\smallbullet'
&
=
\Theta_\smallbullet'
\big(z_\smallbullet',\overline{z}_\smallbullet',\overline{w}_\smallbullet'\big),
\\
\overline{w}_\smallbullet'
&
=
\overline{\Theta}_\smallbullet'
\big(\overline{z}_\smallbullet',
z_\smallbullet',w_\smallbullet'\big).
\endaligned
\]

\medskip

So the {\sl complexification} is the second biholomorphism:
\[
\aligned
z_\smallbullet'
&
=
z_\smallbullet,
\ \ \ \ \ \ \ \ \ \ \ \ \ \ \ \ \ \ \ \ \ \ \ \ \ \ \ \ \ \ \ \ \ 
\ \ \ \ \ \ \ \ \ \ \ \ \ \ \ \ \ \ \ \ 
\underline{z}_\smallbullet'
=
\underline{z}_\smallbullet,
\\
w_\smallbullet'
&
=
\overline{\Theta}_\smallbullet
\big(0_\smallbullet,z_\smallbullet,w_\smallbullet\big),
\ \ \ \ \ \ \ \ \ \ \ \ \ \ \ \ \ \ \ \ \ \ \ \ \ \ \ \ \ \ \ \ \ 
\underline{w}_\smallbullet'
=
\Theta\big(0_\smallbullet,\underline{z}_\smallbullet,
\underline{w}_\smallbullet\big),
\endaligned
\]
sends biholomorphically the extrinsic complexification $M^{e_c}$ onto
the extrinsic complexification ${M'}^{e_c}$, namely:
\[
\aligned
\big(z_\smallbullet,w_\smallbullet,\underline{z}_\smallbullet,
\underline{w}_\smallbullet\big)
\,\in\,M^{e_c}
&
\,\,\,\,\,
\Longleftrightarrow
\,\,\,\,\,
\big(z_\smallbullet',w_\smallbullet',\underline{z}_\smallbullet',
\underline{w}_\smallbullet'\big)
\,\in\,{M'}^{e_c}
\\
&
\,\,\,\,\,
\Longleftrightarrow
\,\,\,\,\,
\big(
z_\smallbullet,\,\overline{\Theta}_\smallbullet
(0_\smallbullet,z_\smallbullet,w_\smallbullet),\,
\underline{z}_\smallbullet,\,
\Theta_\smallbullet
(0_\smallbullet,\underline{z}_\smallbullet,\underline{w}_\smallbullet)
\big)
\,\in\,M',
\endaligned
\]
that is to say, if one takes as equations for ${M'}^{e_c}$:
\[
\underline{w}_\smallbullet'
=
\overline{\Theta}'
\big(\underline{z}_\smallbullet',z_\smallbullet',w_\smallbullet'\big),
\]
one has:
\[
\Theta_\smallbullet
\big(0_\smallbullet,\underline{z}_\smallbullet,\underline{w}_\smallbullet\big)
=
\overline{\Theta}_\smallbullet'
\big(
\underline{z}_\smallbullet,z_\smallbullet,
\overline{\Theta}_\smallbullet
(0_\smallbullet,z_\smallbullet,w_\smallbullet)\big),
\]
still for:
\[
\big(z_\smallbullet,w_\smallbullet,\underline{z}_\smallbullet,\underline{w}_\smallbullet\big)
\,\in\,
M^{e_c}.
\]
But this means that after replacing $w_\smallbullet$ occuring
at the very last place by:
\[
w_\smallbullet
=
\Theta_\smallbullet
\big(
z_\smallbullet,\underline{z}_\smallbullet,\underline{w}_\smallbullet
\big),
\]
one has the $c$ equations:
\[
\Theta_\smallbullet
\big(0_\smallbullet,\underline{z}_\smallbullet,\underline{w}_\smallbullet\big)
\,\equiv\,
\overline{\Theta}_\smallbullet'
\Big(
\underline{z}_\smallbullet,z_\smallbullet,
\overline{\Theta}_\smallbullet
\big(0_\smallbullet,z_\smallbullet,
\Theta_\smallbullet
\big(
z_\smallbullet,\underline{z}_\smallbullet,\underline{w}_\smallbullet
\big)
\big)
\Big),
\]
identically satisfied in:
\[
\C\big\{z_\smallbullet,\underline{z}_\smallbullet,\underline{w}_\smallbullet\big\}^c.
\]

Set in these:
\[
\underline{z}_\smallbullet
:=
0_\smallbullet,
\]
and get:
\[
\underbrace{\Theta_\smallbullet
\big(0_\smallbullet,0_\smallbullet,\underline{w}_\smallbullet\big)}_{
\equiv\,\underline{w}_\smallbullet
\atop
\text{\sf after Step I}}
\,\equiv\,
\overline{\Theta}_\smallbullet'
\Big(
0_\smallbullet',z_\smallbullet,
\underbrace{
\overline{\Theta}_\smallbullet
\big(0_\smallbullet,z_\smallbullet,
\Theta_\smallbullet
\big(
z_\smallbullet,0_\smallbullet,\underline{w}_\smallbullet
\big)
\big)}_{
\equiv\,\underline{w}_\smallbullet\,\,
\text{\sf by one of the functional equations}}
\Big),
\]
that is to say:
\[
\underline{w}_\smallbullet
\equiv
\overline{\Theta}_\smallbullet'
\big(0_\smallbullet',z_\smallbullet,
\underline{w}_\smallbullet\big).
\]
Changing the name of variables to improve notational harmony,
one therefore arrives at one of the announced identities:
\[
w_\smallbullet'
\equiv
\overline{\Theta}_\smallbullet'
\big(0_\smallbullet',z_\smallbullet',w_\smallbullet'\big).
\]
Conjugating and renaming again variables, one also
obtains:
\[
\underline{w}_\smallbullet'
\equiv
\Theta_\smallbullet'
\big(0_\smallbullet',\underline{z}_\smallbullet',\underline{w}_\smallbullet'\big).
\]

In the functional equation:
\[
w_\smallbullet'
\equiv
\Theta'
\big(
z_\smallbullet',\overline{z}_\smallbullet',\,
\overline{\Theta}'(\underline{z}_\smallbullet',
z_\smallbullet',w_\smallbullet')
\big),
\]
setting:
\[
\underline{z}_\smallbullet'
:=
0_\smallbullet',
\]
taking account of what has been just obtained, one gets:
\[
w_\smallbullet'
\equiv
\Theta_\smallbullet'
\big(z_\smallbullet',0_\smallbullet',
\underbrace{\overline{\Theta}_\smallbullet'
(0_\smallbullet',z_\smallbullet',w_\smallbullet')}_{
\equiv\,w_\smallbullet'}\big),
\]
that is to say after renaming variables:
\[
\underline{w}_\smallbullet'
\equiv
\Theta_\smallbullet'
\big(z_\smallbullet',0_\smallbullet',
\underline{w}_\smallbullet'\big).
\]
Conjugating and complexifying:
\[
w_\smallbullet'
\equiv
\overline{\Theta}_\smallbullet'
\big(\underline{z}_\smallbullet',0_\smallbullet',w_\smallbullet'\big).
\]

It only remains to treat $\varphi_\smallbullet'$. But from:
\[
\aligned
\frac{\Theta_\smallbullet'\big(z_\smallbullet',\underline{z}_\smallbullet',
\underline{w}_\smallbullet'\big)-\underline{w}_\smallbullet'}{2\,\isqrt}
&
\,\equiv\,
\varphi_\smallbullet'
\bigg(
z_\smallbullet',\underline{z}_\smallbullet',\,
\frac{\Theta_\smallbullet'\big(z_\smallbullet',\underline{z}_\smallbullet',
\underline{w}_\smallbullet'\big)+\underline{w}_\smallbullet'}{2}
\bigg),
\\
\frac{w_\smallbullet'-\overline{\Theta}_\smallbullet'\big(\underline{z}_\smallbullet',
z_\smallbullet',w_\smallbullet'\big)}{2\,\isqrt}
&
\,\equiv\,
\varphi_\smallbullet'
\bigg(
z_\smallbullet',\underline{z}_\smallbullet',\,
\frac{w_\smallbullet'+\overline{\Theta}_\smallbullet'
\big(\underline{z}_\smallbullet',z_\smallbullet',w_\smallbullet'\big)}{2}
\bigg),
\endaligned
\]
setting either $z_\smallbullet' := 0_\smallbullet'$ or (not simultaneously)
$\underline{ z}_\smallbullet' := 0_\smallbullet'$, one obtains:
\[
0_\smallbullet'
\,\equiv\,
\varphi_\smallbullet'
\big(0_\smallbullet',\overline{z}_\smallbullet',\nu_\smallbullet'\big)
\,\equiv\,
\varphi_\smallbullet'
\big(
z_\smallbullet',0_\smallbullet',\nu_\smallbullet'\big),
\]
which finishes.
\endproof


\bigskip

\section{\sf General class $\text{\sf I}$}
\label{general-class-I}
\HEAD{\ref{general-class-I}.~General class $\text{\sf I}$}{
Jo\"el {\sc Merker} (Paris-Sud)}

\medskip\noindent{\bf Proposition.}
{\em A local real analytic hypersurface passing through the origin:}
\[
0
\,\in\,
M^3
\,\subset\,
\C^2
\]
{\em which belongs to the general Class $\text{\sf I}$, namely such that,
for any local vector field generator $\mathcal{ L}$ of $T^{1, 0}M$:}
\[
\Big\{
\mathcal{L},\,\overline{\mathcal{L}},\,\,
\big[\mathcal{L},\overline{\mathcal{L}}\big]\Big\}
\]
{\em constitute a frame for $\C\otimes_\R TM$, 
may always be represented, in suitable local holomorphic coordinates:}
\[
(z,w)
\in
\C^2
\]
{\em by a specific real analytic equation of the form:}
\[
\boxed{\,\,
\underline{\text{\footnotesize\sf (I):}}
\ \ \ \ \
\aligned
v
=
z\overline{z}
+
z\overline{z}\,{\rm O}_1\big(z,\overline{z}\big)
+
z\overline{z}\,{\rm O}_1(u).\,\,
\rule[-5pt]{0pt}{20pt}
\endaligned}
\]

\proof
Choose first coordinates $(z, w)$ so that $M$ is:
\[
v
=
\varphi\big(z,\overline{z},u\big),
\]
with:
\[
0
=
\varphi_{z}(0)
=
\varphi_{\overline{z}}(0)
=
\varphi_u(0),
\]
and with:
\[
0
\equiv
\varphi\big(0,\overline{z},u\big)
\equiv
\varphi(z,0,u).
\]

Equivalently, the {\em converging} power series expansion:
\[
\varphi
=
\sum_{j,k,l\in\N
\atop
j+k+l\geqslant 2}\,
\underbrace{\varphi_{j,k,l}}_{\in\,\C}\,
z^j\,\overline{z}^k\,u^l
\ \ \ \ \ \ \ \ \ \ \ \ \
{\scriptstyle{(\overline{\varphi_{j,k,l}}\,=\,\varphi_{k,j,l})}},
\]
has vanishing coefficients:
\[
0
=
\varphi_{0,k,l}
=
\varphi_{j,0,l},
\]
so that:
\[
v
=
a\,z\overline{z}
+
z\overline{z}\,{\rm O}_1\big(z,\overline{z}\big)
+
z\overline{z}\,u\,{\rm O}_0\big(z,\overline{z},u\big),
\]
or shortly:
\[
v
=
a\,z\overline{z}
+
z\overline{z}\,{\rm O}_1\big(z,\overline{z}\big)
+
z\overline{z}\,{\rm O}_1\big(u\big),
\]
with of course:
\[
a
\,\in\,\R.
\]

The standard intrinsic generator for $T^{1, 0}M$
({\em see}~\cite{Merker-Pocchiola-Sabzevari-5-CR-II}):
\[
\mathcal{L}
=
\frac{\partial}{\partial z}
\underbrace{-
\frac{\varphi_z}{\isqrt+\varphi_u}}_{=:\,A}\,
\frac{\partial}{\partial u}
\]
has coefficient:
\[
\aligned
A
&
=
-\,\frac{a\,\overline{z}+{\rm O}(2)}{
\isqrt+{\rm O}(2)}
\\
&
=
\isqrt\,a\,\overline{z}
+
{\rm O}(2).
\endaligned
\]

Thus starting from:
\[
\aligned
\mathcal{L}
&
=
\frac{\partial}{\partial z}
+
\big(\isqrt\,a\,\overline{z}
+
{\rm O}(2)\big)\,
\frac{\partial}{\partial u},
\\
\overline{\mathcal{L}}
&
=
\frac{\partial}{\partial\overline{z}}
-
\big(\isqrt\,a\,z+{\rm O}(2)\big)\,
\frac{\partial}{\partial u},
\endaligned
\]
one computes the Lie bracket:
\[
\big[\mathcal{L},\overline{\mathcal{L}}\big]
=
\big(2\,a+{\rm O}(1)\big)\,
\frac{\partial}{\partial u}.
\]

By hypothesis, the three vectors at the origin:
\[
\aligned
\mathcal{L}\big\vert_0
&
=
\frac{\partial}{\partial z}
\bigg\vert_0,
\\
\overline{\mathcal{L}}\big\vert_0
&
=
\frac{\partial}{\partial\overline{z}}
\bigg\vert_0,
\\
\isqrt\,\big[\mathcal{L},\overline{\mathcal{L}}\big]\Big\vert_0
&
=
2a\,
\frac{\partial}{\partial u},
\endaligned
\]
should make a basis for $\C\otimes_\R T_0M$ having:
\[
{\bf 3}
=
\dim_\C\big(\C\otimes_\R T_0M\big),
\]
and this forces:
\[
a\neq 0.
\]

If $a < 0$, changing:
\[
w
\longmapsto
-\,w,
\]
one makes:
\[
a>0,
\]
and lastly changing:
\[
z\,\longmapsto\,
\sqrt{a}\,z
\]
one makes $a = 1$, which finishes.
\endproof

\noindent{\bf Scholium.}
{\em In such elementarily normalized coordinates, one has
the diagonal normalization at the origin:}
\[
\aligned
\mathcal{L}\big\vert_0
&
=
\frac{\partial}{\partial z}
\bigg\vert_0,
\\
\overline{\mathcal{L}}\big\vert_0
&
=
\frac{\partial}{\partial\overline{z}}
\bigg\vert_0,
\\
\big[\mathcal{L},\overline{\mathcal{L}}\big]
\Big\vert_0
&
=
-\,2\,\isqrt\,\,
\frac{\partial}{\partial u}
\bigg\vert_0,
\endaligned
\]
{\em which conveniently fixes ideas when performing explicitly 
the Cartan equivalence procedure \cite{
Merker-Sabzevari-Cartan-M3-C2}.\qed}

\medskip

Taking zero remainders, one obtains the:
\[
\boxed{\,\,
{\text{\footnotesize\sf Model (I):}}
\ \ \ \ \ \ \ \ \ 
v
=
z\overline{z}.\,\,
\rule[-4pt]{0pt}{15pt}}
\]


\bigskip

\section{\sf General class $\text{\sf II}$}
\label{general-class-II}
\HEAD{\ref{general-class-II}.~General class $\text{\sf II}$}{
Jo\"el {\sc Merker} (Paris-Sud)}

\medskip

\medskip\noindent{\bf Proposition.}
{\em A local real analytic $4$-dimensional
CR-generic submanifold passing through the origin:}
\[
0
\,\in\,
M^4
\,\subset\,
\C^3
\]
{\em which belongs to the general Class $\text{\sf II}$, namely such that,
for any local vector field generator $\mathcal{ L}$ of $T^{1, 0}M$:}
\[
\Big\{
\mathcal{L},\,\overline{\mathcal{L}},\,\,
\big[\mathcal{L},\overline{\mathcal{L}}\big],\,\,
\big[\mathcal{L},\,\big[\mathcal{L},\overline{\mathcal{L}}\big]\big]
\Big\}
\]
{\em constitute a frame for $\C\otimes_\R TM$, 
may always be represented, in suitable local holomorphic coordinates:}
\[
(z,w_1,w_2)
\in
\C^3
\]
{\em by two specific real analytic equations of the form:}
\[
\boxed{\,\,
\underline{\text{\footnotesize\sf (II):}}
\ \ \ \ \
\aligned
v_1
&
=
z\overline{z}
\ \ \ \ \ \ \ \ \ \ \ \ 
+
z\overline{z}\,{\rm O}_2\big(z,\overline{z}\big)
+
z\overline{z}\,{\rm O}_1(u_1)
+
z\overline{z}\,{\rm O}_1(u_2),\,\,
\\
v_2
&
=
z^2\overline{z}
+
z\overline{z}^2
+
z\overline{z}\,{\rm O}_2\big(z,\overline{z}\big)
+
z\overline{z}\,{\rm O}_1(u_1)
+
z\overline{z}\,{\rm O}_1(u_2).\,\,
\endaligned
\rule[-5pt]{0pt}{27pt}}
\]

\proof
Choose first coordinates $(z, w_1, w_2)$ so that $M$ is:
\[
\aligned
v_1
&
=
\varphi_1\big(z,\overline{z},u_1,u_2\big),
\\
v_2
&
=
\varphi_2\big(z,\overline{z},u_1,u_2\big),
\endaligned
\]
with:
\[
\aligned
0
&
=
\varphi_{1,z}(0)
=
\varphi_{1,\overline{z}}(0)
=
\varphi_{1,u_1}(0)
=
\varphi_{1,u_2}(0),
\\
0
&
=
\varphi_{2,z}(0)
=
\varphi_{2,\overline{z}}(0)
=
\varphi_{2,u_1}(0)
=
\varphi_{2,u_2}(0),
\endaligned
\]
and with:
\[
\aligned
0
&
\equiv
\varphi_1\big(0,\overline{z},u_1,u_2)
\equiv
\varphi_1\big(z,0,u_1,u_2\big),
\\
0
&
\equiv
\varphi_2\big(0,\overline{z},u_1,u_2)
\equiv
\varphi_2\big(z,0,u_1,u_2\big).
\endaligned
\]

Equivalently, the {\em converging} power series expansions:
\[
\aligned
\varphi_1
&
=
\sum_{j,k,l_1,l_2\in\N
\atop
j+k+l_1+l_2\geqslant 2}\,
\underbrace{\varphi_{1,j,k,l_1,l_2}}_{\in\,\C}\,
z^j\,\overline{z}^k\,u_1^{l_1}\,u_2^{l_2}
\ \ \ \ \ \ \ \ \ \ \ \ \
{\scriptstyle{(\overline{\varphi_{1,j,k,l_1,l_2}}\,=\,
\varphi_{1,k,j,l_1,l_2})}},
\\
\varphi_2
&
=
\sum_{j,k,l_1,l_2\in\N
\atop
j+k+l_1+l_2\geqslant 2}\,
\underbrace{\varphi_{2,j,k,l_1,l_2}}_{\in\,\C}\,
z^j\,\overline{z}^k\,u_1^{l_1}\,u_2^{l_2}
\ \ \ \ \ \ \ \ \ \ \ \ \
{\scriptstyle{(\overline{\varphi_{2,j,k,l_1,l_2}}\,=\,
\varphi_{2,k,j,l_1,l_2})}},
\endaligned
\]
have vanishing coefficients:
\[
\aligned
0
&
=
\varphi_{1,0,k,l_1,l_2}
=
\varphi_{1,j,0,l_1,l_2},
\\
0
&
=
\varphi_{2,0,k,l_1,l_2}
=
\varphi_{2,j,0,l_1,l_2},
\endaligned
\]
so that:
\[
\aligned
v_1
&
=
a_1\,z\overline{z}
+
z\overline{z}\,{\rm O}_1\big(z,\overline{z}\big)
+
z\overline{z}\,u_1\,{\rm O}_0
\big(z,\overline{z},u_1,u_2\big)
+
z\overline{z}\,u_2\,{\rm O}_0
\big(z,\overline{z},u_1,u_2\big),
\\
v_2
&
=
a_2\,z\overline{z}
+
z\overline{z}\,{\rm O}_1\big(z,\overline{z}\big)
+
z\overline{z}\,u_1\,{\rm O}_0
\big(z,\overline{z},u_1,u_2\big)
+
z\overline{z}\,u_2\,{\rm O}_0
\big(z,\overline{z},u_1,u_2\big),
\endaligned
\]
or shortly:
\[
\aligned
v_1
&
=
a_1\,z\overline{z}
+
z\overline{z}\,{\rm O}_1\big(z,\overline{z}\big)
+
z\overline{z}\,{\rm O}_1(u_1)
+
z\overline{z}\,{\rm O}_1(u_2)
\\
v_2
&
=
a_2\,z\overline{z}
+
z\overline{z}\,{\rm O}_1\big(z,\overline{z}\big)
+
z\overline{z}\,{\rm O}_1(u_1)
+
z\overline{z}\,{\rm O}_1(u_2),
\endaligned
\]
with of course:
\[
\aligned
&
a_1\,\in\,\R,
\\
&
a_2\,\in\,\R.
\endaligned
\]

The standard intrinsic generator for $T^{1, 0}M$
({\em see}~\cite{ Merker-Pocchiola-Sabzevari-5-CR-II}):
\[
\mathcal{L}
=
\frac{\partial}{\partial z}
+
A_1\,\frac{\partial}{\partial u_1}
+
A_2\,\frac{\partial}{\partial u_2}
\]
has its two coefficients given by the formulas:
\[
\aligned
A_1
&
=
\left\vert\!
\begin{array}{cc}
-\varphi_{1,z} & \varphi_{1,u_2}
\\
-\varphi_{2,z} & \isqrt+\varphi_{2,u_2}
\end{array}
\!\right\vert
\bigg/
\left\vert\!
\begin{array}{cc}
\isqrt+\varphi_{1,u_1} & \varphi_{1,u_2}
\\
\varphi_{2,u_1} & \isqrt+\varphi_{2,u_2}
\end{array}
\!\right\vert,
\\
A_2
&
=
\left\vert\!
\begin{array}{cc}
\isqrt+\varphi_{1,u_1} & -\varphi_{1,z}
\\
\varphi_{2,u_1} & -\varphi_{2,z}
\end{array}
\!\right\vert
\bigg/
\left\vert\!
\begin{array}{cc}
\isqrt+\varphi_{1,u_1} & \varphi_{1,u_2}
\\
\varphi_{2,u_1} & \isqrt+\varphi_{2,u_2}
\end{array}
\!\right\vert.
\endaligned
\]

An approximant for $A_1$ is:
\[
\aligned
A_1
&
=
\left\vert\!
\begin{array}{cc}
-a_1\overline{z}+{\rm O}(2) & {\rm O}(2)
\\
-a_2\overline{z}+{\rm O}(2) & \isqrt+{\rm O}(2)
\end{array}
\!\right\vert
\bigg/
\left\vert\!
\begin{array}{cc}
\isqrt+{\rm O}(2) & {\rm O}(2)
\\
{\rm O}(2) & \isqrt+{\rm O}(2)
\end{array}
\!\right\vert
\\
&
=
\isqrt\,a_1\,\overline{z}
+
{\rm O}(2),
\endaligned
\]
and similarly:
\[
A_2
=
\isqrt\,a_2\,\overline{z}
+
{\rm O}(2),
\]
so that:
\[
\aligned
\mathcal{L}
&
=
\frac{\partial}{\partial z}
+
\isqrt\,
\big(a_1\,\overline{z}+{\rm O}(2)\big)\,
\frac{\partial}{\partial u_1}
+
\isqrt\,
\big(a_2\,\overline{z}+{\rm O}(2)\big)\,
\frac{\partial}{\partial u_2},
\\
\overline{\mathcal{L}}
&
=
\frac{\partial}{\partial\overline{z}}
-
\isqrt\,
\big(a_1\,z+{\rm O}(2)\big)\,
\frac{\partial}{\partial u_1}
-
\isqrt\,
\big(a_2\,z+{\rm O}(2)\big)\,
\frac{\partial}{\partial u_2},
\endaligned
\]
whence:
\[
\isqrt\,
\big[\mathcal{L},\overline{\mathcal{L}}\big]
=
\big(2\,a_1+{\rm O}(1)\big)\,
\frac{\partial}{\partial u_1}
+
\big(2\,a_2+{\rm O}(1)\big)\,
\frac{\partial}{\partial u_2},
\]
and at the origin:
\[
\isqrt\,
\big[\mathcal{L},\overline{\mathcal{L}}\big]
\Big\vert_0
=
2a_1\,\frac{\partial}{\partial u_1}
\bigg\vert_0
+
2a_2\,\frac{\partial}{\partial u_2}
\bigg\vert_0.
\]

Since the three vectors:
\[
\mathcal{L}\big\vert_0, 
\ \ \ \ \ \ \ \ \ \ \ \ \
\overline{\mathcal{L}}\big\vert_0,
\ \ \ \ \ \ \ \ \ \ \ \ \
\big[\mathcal{L},\overline{\mathcal{L}}\big] 
\Big\vert_0
\]
must in particular by hypothesis necessarily be $\C$-linearly independent
at the origin, one has:
\[
a_1\neq 0
\ \ \ \ \ \ \ \ \ \ \ \ \ \ \ \ \ \ \
\text{\rm or}
\ \ \ \ \ \ \ \ \ \ \ \ \ \ \ \ \ \ \
a_2\neq 0.
\]

Say:
\[
a_1\neq 0.
\]
Then do if necessary:
\[
w_1
\,\longmapsto\,
-\,w_1
\]
to make:
\[
a_1>0,
\]
and do:
\[
z
\,\longmapsto\,
\sqrt{a_1}\,z,
\]
to make:
\[
a_1
=
1.
\]

Writing out now the third-order $(z, \overline{z})$-terms, one
therefore comes to:
\[
\aligned
v_1
&
=
z\overline{z}
+
\alpha_1\,z^2\overline{z}
+
\overline{\alpha}_1\,z\overline{z}^2
+
z\overline{z}\,
{\rm O}_2\big(z,\overline{z}\big)
+
z\overline{z}\,{\rm O}_1(u_1)
+
z\overline{z}\,{\rm O}_1(u_2),
\\
v_2
&
=
a_2\,z\overline{z}
+
\alpha_2\,z^2\overline{z}
+
\overline{\alpha}_2\,z\overline{z}^2
+
z\overline{z}\,
{\rm O}_2\big(z,\overline{z}\big)
+
z\overline{z}\,{\rm O}_1(u_1)
+
z\overline{z}\,{\rm O}_1(u_2),
\endaligned
\]
with a new $\alpha_2 \in \R$ and with of course:
\[
\alpha_1
\in
\C
\ \ \ \ \ \ \ \ \ \ \ \ \ \ \ \ \ \ \
\text{\rm and}
\ \ \ \ \ \ \ \ \ \ \ \ \ \ \ \ \ \ \
\alpha_2
\in
\C.
\]

Changing:
\[
w_2
\,\longmapsto\,
w_2-a_2\,w_1,
\]
one makes:
\[
a_2
=
0.
\]
In the process, $\alpha_2$ is modified.

Changing:
\[
z
\,\longmapsto\,
z+\alpha_1\,z^2,
\]
one makes:
\[
\alpha_1
=
0.
\]

One thus arrives at:
\[
\aligned
v_1
&
=
z\overline{z}
\ \ \ \ \ \ \ \ \ \ \ \ \ \ \ \ \ \ \ \ \ \,
+
z\overline{z}\,
{\rm O}_2\big(z,\overline{z}\big)
+
z\overline{z}\,{\rm O}_1(u_1)
+
z\overline{z}\,{\rm O}_1(u_2),
\\
v_2
&
=
\alpha_2\,z^2\overline{z}
+
\overline{\alpha}_2\,z\overline{z}^2
+
z\overline{z}\,
{\rm O}_2\big(z,\overline{z}\big)
+
z\overline{z}\,{\rm O}_1(u_1)
+
z\overline{z}\,{\rm O}_1(u_2).
\endaligned
\]

Next, a deeper approximant for $A_1$ is:
\[
\!\!\!\!\!\!\!\!\!\!\!\!\!\!\!\!\!\!\!\!
\!\!\!\!\!\!\!\!\!\!\!\!\!\!\!\!\!\!\!\!
\aligned
A_1
&
=
\left\vert\!
\begin{array}{cc}
-\overline{z}
+
{\rm O}_3(z,\overline{z})+\overline{z}{\rm O}_1(u_1,u_2)
& 
z\overline{z}{\rm O}(0)
\\
-2\alpha_2z\overline{z}-\overline{\alpha}_2\overline{z}^2
+
{\rm O}_3(z,\overline{z})+\overline{z}{\rm O}_1(u_1,u_2)
& 
\isqrt+z\overline{z}{\rm O}(0)
\end{array}
\!\right\vert
\bigg/
\left\vert\!
\begin{array}{cc}
\isqrt+z\overline{z}{\rm O}(0) & z\overline{z}{\rm O}(0)
\\
z\overline{z}{\rm O}(0) & \isqrt+z\overline{z}{\rm O}(0)
\end{array}
\!\right\vert
\\
&
=
\isqrt\,\overline{z}
+
\overline{z}\,{\rm O}_1(u_1,u_2)
+
{\rm O}(3),
\endaligned
\]
while an approximant for $A_2$ is:
\[
\!\!\!\!\!\!\!\!\!\!\!\!\!\!\!\!\!\!\!\!
\!\!\!\!\!\!\!\!\!\!\!\!\!\!\!\!\!\!\!\!
\aligned
A_2
&
=
\left\vert\!
\begin{array}{cc}
\isqrt+z\overline{z}{\rm O}(0)
&
-\overline{z}
+
{\rm O}_3(z,\overline{z})+\overline{z}{\rm O}_1(u_1,u_2)
\\
z\overline{z}{\rm O}(0)
&
-2\alpha_2z\overline{z}-\overline{\alpha}_2\overline{z}^2
+
{\rm O}_3(z,\overline{z})+\overline{z}{\rm O}_1(u_1,u_2)
\end{array}
\!\right\vert
\bigg/
\left\vert\!
\begin{array}{cc}
\isqrt+z\overline{z}{\rm O}(0) & z\overline{z}{\rm O}(0)
\\
z\overline{z}{\rm O}(0) & \isqrt+z\overline{z}{\rm O}(0)
\end{array}
\!\right\vert
\\
&
=
2\,\isqrt\,\alpha_2\,z\overline{z}
+
\isqrt\,\overline{\alpha}_2\,\overline{z}^2
+
\overline{z}\,{\rm O}_1(u_1,u_2)
+
{\rm O}(3),
\endaligned
\]
so that:
\[
\!\!\!\!\!\!\!\!\!\!\!\!\!\!\!\!\!\!\!\!
\!\!\!\!\!\!\!\!\!\!\!\!\!\!\!\!\!\!\!\!
\aligned
\mathcal{L}
&
=
\frac{\partial}{\partial z}
+
\isqrt\,
\big(\overline{z}
+
\overline{z}\,{\rm O}_1(u_1,u_2)
+
{\rm O}(3)\big)\,
\frac{\partial}{\partial u_1}
+
\isqrt\,
\big(
2\alpha_2\,z\overline{z}
+
\overline{\alpha}_2\,\overline{z}^2
+
\overline{z}\,{\rm O}_1(u_1,u_2)
+
{\rm O}(3)\big)\,
\frac{\partial}{\partial u_2},
\\
\overline{\mathcal{L}}
&
=
\frac{\partial}{\partial\overline{z}}
-
\isqrt\,
\big(
z
+
z\,{\rm O}_1(u_1,u_2)
+
{\rm O}(3)\big)\,
\frac{\partial}{\partial u_1}
-
\isqrt\,
\big(
2\overline{\alpha}_2\,z\overline{z}
+
\alpha_2\,z^2
+
z\,{\rm O}_1(u_1,u_2)
+
{\rm O}(3)\big)\,
\frac{\partial}{\partial u_2},
\endaligned
\]
whence:
\[
\big[\mathcal{L},\overline{\mathcal{L}}\big]
=
-\,\isqrt\,
\big(
2
+
{\rm O}_1(u_1,u_2)
+
{\rm O}(2)\big)\,
\frac{\partial}{\partial u_1}
-
\isqrt\,
\big(
4\alpha_2\,z+4\overline{\alpha}_2\,\overline{z}
+
{\rm O}_1(u_1,u_2)
+
{\rm O}(2)\big)\,
\frac{\partial}{\partial u_2},
\]
whence further:
\[
\aligned
\big[\mathcal{L},\,\big[\mathcal{L},\overline{\mathcal{L}}\big]\big]
&
=
{\rm O}(1)\,\frac{\partial}{\partial u_1}
-
\isqrt\,\big(4\alpha_2+{\rm O}(1)\big)\,
\frac{\partial}{\partial u_2},
\endaligned
\]
and lastly, at the origin:
\[
\aligned
\mathcal{L}\big\vert_0
&
=
\frac{\partial}{\partial z}
\bigg\vert_0,
\\
\overline{\mathcal{L}}\big\vert_0
&
=
\frac{\partial}{\partial\overline{z}}
\bigg\vert_0,
\\
\big[\mathcal{L},\overline{\mathcal{L}}\big]
\Big\vert_0
&
=
-\,2\,\isqrt\,
\frac{\partial}{\partial u_1}
\bigg\vert_0,
\\
\big[\mathcal{L},\,\big[\mathcal{L},\overline{\mathcal{L}}\big]\big]
\Big\vert_0
&
=
-\,4\,\isqrt\,\alpha_2\,
\frac{\partial}{\partial u_2}
\bigg\vert_0.
\endaligned
\]

Since these four vectors should by hypothesis constitute
a basis for: 
\[
\C\otimes_\R T_0M
\,=\,
\C\frac{\partial}{\partial z}
\bigg\vert_0
\oplus
\C\frac{\partial}{\partial\overline{z}}
\bigg\vert_0
\oplus
\C\frac{\partial}{\partial u_1}
\bigg\vert_0
\oplus
\C\frac{\partial}{\partial u_2}
\bigg\vert_0,
\] 
one necessarily has:
\[
\alpha_2
\neq
0.
\]

Using a complex dilation:
\[
z
\,\longmapsto\,
\lambda\,z,
\]
with $\lambda$ solution of (exercise):
\[
1
=
\alpha_2\,\lambda^2\,\overline{\lambda},
\]
one makes:
\[
\alpha_2
=
1,
\]
so that the second equation receives the form announced:
\[
v_2
=
z^2\overline{z}
+
z\overline{z}^2
+
z\overline{z}\,{\rm O}_2\big(z,\overline{z}\big)
+
z\overline{z}\,{\rm O}_1(u_1)
+
z\overline{z}\,{\rm O}_1(u_2).
\]

In the process, the first equation is changed to:
\[
v_1
=
\lambda\overline{\lambda}\,z\overline{z}
+
z\overline{z}\,{\rm O}_2\big(z,\overline{z}\big)
+
z\overline{z}\,{\rm O}_1(u_1)
+
z\overline{z}\,{\rm O}_1(u_2),
\]
but one can yet replace:
\[
w_1
\,\longmapsto\,
\frac{w_1}{\lambda\overline{\lambda}},
\]
to keep $1\cdot z \overline{ z}$, which concludes.
\endproof

\noindent{\bf Scholium.}
{\em In such elementarily normalized coordinates, one has
the diagonal normalization at the origin:}
\[
\aligned
\mathcal{L}\big\vert_0
&
=
\frac{\partial}{\partial z}
\bigg\vert_0,
\\
\overline{\mathcal{L}}\big\vert_0
&
=
\frac{\partial}{\partial\overline{z}}
\bigg\vert_0,
\\
\big[\mathcal{L},\overline{\mathcal{L}}\big]
\Big\vert_0
&
=
-\,2\,\isqrt\,
\frac{\partial}{\partial u_1}
\bigg\vert_0,
\\
\big[\mathcal{L},\,\big[\mathcal{L},\overline{\mathcal{L}}\big]\big]
\Big\vert_0
&
=
-\,4\,\isqrt\,
\frac{\partial}{\partial u_2}
\bigg\vert_0,
\endaligned
\]
{\em which conveniently fixes ideas 
when performing explicitly the Cartan equivalence
procedure.\qed}

\medskip

Taking zero remainders, one obtains Beloshapka's cubic:
\[
\boxed{\,\,
{\text{\footnotesize\sf Model (II):}}
\ \ \ \ \ \ \ \ \ 
\aligned
v_1
&
=
z\overline{z},\,\,
\\
v_2
&
=
z^2\overline{z}+z\overline{z}^2.
\endaligned
\rule[-4pt]{0pt}{27pt}}
\]


\bigskip

\section{\sf General class $\text{\sf III}_{\text{\sf 1}}$}
\label{general-class-III-1}
\HEAD{\ref{general-class-III-1}.~General class 
$\text{\sf III}_{\text{\sf 1}}$}{
Jo\"el {\sc Merker} (Paris-Sud)}

\medskip

\medskip\noindent{\bf Proposition.}
{\em A local real analytic $5$-dimensional
CR-generic submanifold passing through the origin:}
\[
0
\,\in\,
M^5
\,\subset\,
\C^4
\]
{\em which belongs to the general Class $\text{\sf III}_{
\text{\sf 1}}$, namely such that,
for any local vector field generator $\mathcal{ L}$ of $T^{1, 0}M$:}
\[
\Big\{
\mathcal{L},\,\overline{\mathcal{L}},\,\,
\big[\mathcal{L},\overline{\mathcal{L}}\big],\,\,
\big[\mathcal{L},\,\big[\mathcal{L},\overline{\mathcal{L}}\big]\big],\,\,
\big[\overline{\mathcal{L}},\,
\big[\mathcal{L},\overline{\mathcal{L}}\big]\big]
\Big\}
\]
{\em constitute a frame for $\C\otimes_\R TM$, 
may always be represented, in suitable local holomorphic coordinates:}
\[
(z,w_1,w_2,w_3)
\in
\C^3
\]
{\em by three specific real analytic equations of the form:}
\[
\boxed{\,\,
\underline{\text{\footnotesize\sf 
$\text{\sf (III)}_{\text{\sf 1}}$:}}
\ \ \ \ \
\aligned
v_1
&
=
z\overline{z}
\ \ \ \ \ \ \ \ \ \ \ \ \ \ \ \ \ \ \ \ \ \ \
+
z\overline{z}\,{\rm O}_2\big(z,\overline{z}\big)
+
z\overline{z}\,{\rm O}_1(u_1)
+
z\overline{z}\,{\rm O}_1(u_2)
+
z\overline{z}\,{\rm O}_1(u_3),\,\,
\\
v_2
&
=
z^2\overline{z}
+
z\overline{z}^2
\ \ \ \ \ \ \ \ \ \ \,
+
z\overline{z}\,{\rm O}_2\big(z,\overline{z}\big)
+
z\overline{z}\,{\rm O}_1(u_1)
+
z\overline{z}\,{\rm O}_1(u_2)
+
z\overline{z}\,{\rm O}_1(u_3),\,\,
\\
v_3
&
=
\isqrt\,\big(z^2\overline{z}
-
z\overline{z}^2\big)
+
z\overline{z}\,{\rm O}_2\big(z,\overline{z}\big)
+
z\overline{z}\,{\rm O}_1(u_1)
+
z\overline{z}\,{\rm O}_1(u_2)
+
z\overline{z}\,{\rm O}_1(u_3).\,\,
\endaligned
\rule[-5pt]{0pt}{27pt}}
\]

\proof
Choose first coordinates $(z, w_1, w_2,w_3)$ so that $M$ is:
\[
\aligned
v_1
&
=
\varphi_1\big(z,\overline{z},u_1,u_2,u_3\big),
\\
v_2
&
=
\varphi_2\big(z,\overline{z},u_1,u_2,u_3\big),
\\
v_3
&
=
\varphi_3\big(z,\overline{z},u_1,u_2,u_3\big),
\endaligned
\]
with:
\[
\aligned
0
&
=
\varphi_{1,z}(0)
=
\varphi_{1,\overline{z}}(0)
=
\varphi_{1,u_1}(0)
=
\varphi_{1,u_2}(0)
=
\varphi_{1,u_3}(0),
\\
0
&
=
\varphi_{2,z}(0)
=
\varphi_{2,\overline{z}}(0)
=
\varphi_{2,u_1}(0)
=
\varphi_{2,u_2}(0)
=
\varphi_{2,u_3}(0),
\\
0
&
=
\varphi_{3,z}(0)
=
\varphi_{3,\overline{z}}(0)
=
\varphi_{3,u_1}(0)
=
\varphi_{3,u_2}(0)
=
\varphi_{3,u_3}(0),
\endaligned
\]
and with:
\[
\aligned
0
&
\equiv
\varphi_1\big(0,\overline{z},u_1,u_2,u_3)
\equiv
\varphi_1\big(z,0,u_1,u_2,u_3\big),
\\
0
&
\equiv
\varphi_2\big(0,\overline{z},u_1,u_2,u_3)
\equiv
\varphi_2\big(z,0,u_1,u_2,u_3\big),
\\
0
&
\equiv
\varphi_3\big(0,\overline{z},u_1,u_2,u_3)
\equiv
\varphi_3\big(z,0,u_1,u_2,u_3\big).
\endaligned
\]

Equivalently, the {\em converging} power series expansions:
\[
\aligned
\varphi_1
&
=
\sum_{j,k,l_1,l_2,l_3\in\N
\atop
j+k+l_1+l_2+l_3\geqslant 2}\,
\underbrace{\varphi_{1,j,k,l_1,l_2,l_3}}_{\in\,\C}\,
z^j\,\overline{z}^k\,u_1^{l_1}\,u_2^{l_2}\,u_3^{l_3}
\ \ \ \ \ \ \ \ \ \ \ \ \
{\scriptstyle{(\overline{\varphi_{1,j,k,l_1,l_2,l_3}}\,=\,
\varphi_{1,k,j,l_1,l_2,l_3})}},
\\
\varphi_2
&
=
\sum_{j,k,l_1,l_2,l_3\in\N
\atop
j+k+l_1+l_2+l_3\geqslant 2}\,
\underbrace{\varphi_{2,j,k,l_1,l_2,l_3}}_{\in\,\C}\,
z^j\,\overline{z}^k\,u_1^{l_1}\,u_2^{l_2}\,u_3^{l_3}
\ \ \ \ \ \ \ \ \ \ \ \ \
{\scriptstyle{(\overline{\varphi_{2,j,k,l_1,l_2,l_3}}\,=\,
\varphi_{2,k,j,l_1,l_2,l_3})}},
\\
\varphi_3
&
=
\sum_{j,k,l_1,l_2,l_3\in\N
\atop
j+k+l_1+l_2+l_3\geqslant 2}\,
\underbrace{\varphi_{3,j,k,l_1,l_2,l_3}}_{\in\,\C}\,
z^j\,\overline{z}^k\,u_1^{l_1}\,u_2^{l_2}\,u_3^{l_3}
\ \ \ \ \ \ \ \ \ \ \ \ \
{\scriptstyle{(\overline{\varphi_{3,j,k,l_1,l_2,l_3}}\,=\,
\varphi_{3,k,j,l_1,l_2,l_3})}},
\endaligned
\]
have vanishing coefficients:
\[
\aligned
0
&
=
\varphi_{1,0,k,l_1,l_2,l_3}
=
\varphi_{1,j,0,l_1,l_2,l_3},
\\
0
&
=
\varphi_{2,0,k,l_1,l_2,l_3}
=
\varphi_{2,j,0,l_1,l_2,l_3},
\\
0
&
=
\varphi_{3,0,k,l_1,l_2,l_3}
=
\varphi_{3,j,0,l_1,l_2,l_3},
\endaligned
\]
so that:
\[
\aligned
v_1
&
=
a_1\,z\overline{z}
+
z\overline{z}\,{\rm O}_1\big(z,\overline{z}\big)
+
z\overline{z}\,{\rm O}_1(u_1)
+
z\overline{z}\,{\rm O}_1(u_2)
+
z\overline{z}\,{\rm O}_1(u_3),
\\
v_2
&
=
a_2\,z\overline{z}
+
z\overline{z}\,{\rm O}_1\big(z,\overline{z}\big)
+
z\overline{z}\,{\rm O}_1(u_1)
+
z\overline{z}\,{\rm O}_1(u_2)
+
z\overline{z}\,{\rm O}_1(u_3),
\\
v_3
&
=
a_3\,z\overline{z}
+
z\overline{z}\,{\rm O}_1\big(z,\overline{z}\big)
+
z\overline{z}\,{\rm O}_1(u_1)
+
z\overline{z}\,{\rm O}_1(u_2)
+
z\overline{z}\,{\rm O}_1(u_3),
\endaligned
\]
with of course:
\[
\aligned
&
a_1\,\in\,\R,
\\
&
a_2\,\in\,\R,
\\
&
a_3\,\in\,\R.
\endaligned
\]

The standard intrinsic generator for $T^{1, 0}M$
({\em see}~\cite{ Merker-Pocchiola-Sabzevari-5-CR-II}):
\[
\mathcal{L}
=
\frac{\partial}{\partial z}
+
A_1\,\frac{\partial}{\partial u_1}
+
A_2\,\frac{\partial}{\partial u_2}
+
A_3\,\frac{\partial}{\partial u_3}
\]
has its three coefficients given by the formulas:
\[
\aligned
A_1
&
=
\left\vert\!
\begin{array}{ccc}
-\varphi_{1,z} & \varphi_{1,u_2} & \varphi_{1,u_3}
\\
-\varphi_{2,z} & \isqrt+\varphi_{2,u_2} & \varphi_{2,u_3}
\\
-\varphi_{3,z} & \varphi_{3,u_2} & \isqrt+\varphi_{3,u_3}
\end{array}
\!\right\vert
\bigg/
\left\vert\!
\begin{array}{ccc}
\isqrt+\varphi_{1,u_1} & \varphi_{1,u_2} & \varphi_{1,u_3}
\\
\varphi_{2,u_1} & \isqrt+\varphi_{2,u_2} & \varphi_{2,u_3}
\\
\varphi_{3,u_1} & \varphi_{3,u_2} & \isqrt+\varphi_{3,u_3}
\end{array}
\!\right\vert,
\endaligned
\]
\[
\aligned
A_2
&
=
\left\vert\!
\begin{array}{ccc}
\isqrt+\varphi_{1,u_1} & -\varphi_{1,z} & \varphi_{1,u_3}
\\
\varphi_{2,u_1} & -\varphi_{2,z} &  \varphi_{2,u_3}
\\
\varphi_{3,u_1} & -\varphi_{3,z} & \isqrt+\varphi_{3,u_3}
\end{array}
\!\right\vert
\bigg/
\left\vert\!
\begin{array}{ccc}
\isqrt+\varphi_{1,u_1} & \varphi_{1,u_2} & \varphi_{1,u_3}
\\
\varphi_{2,u_1} & \isqrt+\varphi_{2,u_2} & \varphi_{2,u_3}
\\
\varphi_{3,u_1} & \varphi_{3,u_2} & \isqrt+\varphi_{3,u_3}
\end{array}
\!\right\vert,
\endaligned
\]
\[
\aligned
A_3
&
=
\left\vert\!
\begin{array}{ccc}
\isqrt+\varphi_{1,u_1} & \varphi_{1,u_2} & -\varphi_{1,z}
\\
\varphi_{2,u_1} &  \isqrt+\varphi_{2,u_2} & -\varphi_{2,z}
\\
\varphi_{3,u_1} & \varphi_{3,u_2} & -\varphi_{3,z}
\end{array}
\!\right\vert
\bigg/
\left\vert\!
\begin{array}{ccc}
\isqrt+\varphi_{1,u_1} & \varphi_{1,u_2} & \varphi_{1,u_3}
\\
\varphi_{2,u_1} & \isqrt+\varphi_{2,u_2} & \varphi_{2,u_3}
\\
\varphi_{3,u_1} & \varphi_{3,u_2} & \isqrt+\varphi_{3,u_3}
\end{array}
\!\right\vert.
\endaligned
\]

An approximant for $A_1$ is:
\[
\!\!\!\!\!\!\!\!\!\!\!\!\!\!\!\!\!\!\!\!
\aligned
A_1
&
=
\left\vert\!
\begin{array}{ccc}
-a_1\overline{z}+{\rm O}(2) & {\rm O}(2) & {\rm O}(2)
\\
-a_2\overline{z}+{\rm O}(2) & \isqrt+{\rm O}(2) & {\rm O}(2)
\\
-a_3\overline{z}+{\rm O}(2) & {\rm O}(2) & \isqrt+{\rm O}(2)
\end{array}
\!\right\vert
\bigg/
\left\vert\!
\begin{array}{ccc}
\isqrt+{\rm O}(2) & {\rm O}(2) & {\rm O}(2)
\\
{\rm O}(2) & \isqrt+{\rm O}(2) & {\rm O}(2)
\\
{\rm O}(2) & {\rm O}(2) & \isqrt+{\rm O}(2)
\end{array}
\!\right\vert
\\
&
=
\isqrt\,a_1\,\overline{z}
+
{\rm O}(2),
\endaligned
\]
and similarly:
\[
\aligned
A_2
=
\isqrt\,a_2\,\overline{z}
+
{\rm O}(2),
\\
A_3
=
\isqrt\,a_3\,\overline{z}
+
{\rm O}(2),
\endaligned
\]
so that:
\[
\aligned
\mathcal{L}
&
=
\frac{\partial}{\partial z}
+
\isqrt\,
\big(a_1\,\overline{z}+{\rm O}(2)\big)\,
\frac{\partial}{\partial u_1}
+
\isqrt\,
\big(a_2\,\overline{z}+{\rm O}(2)\big)\,
\frac{\partial}{\partial u_2}
+
\isqrt\,
\big(a_3\,\overline{z}+{\rm O}(2)\big)\,
\frac{\partial}{\partial u_3},
\\
\overline{\mathcal{L}}
&
=
\frac{\partial}{\partial\overline{z}}
-
\isqrt\,
\big(a_1\,z+{\rm O}(2)\big)\,
\frac{\partial}{\partial u_1}
-
\isqrt\,
\big(a_2\,z+{\rm O}(2)\big)\,
\frac{\partial}{\partial u_2}
-
\isqrt\,
\big(a_3\,z+{\rm O}(2)\big)\,
\frac{\partial}{\partial u_3},
\endaligned
\]
whence:
\[
\isqrt\,
\big[\mathcal{L},\overline{\mathcal{L}}\big]
=
\big(2\,a_1+{\rm O}(1)\big)\,
\frac{\partial}{\partial u_1}
+
\big(2\,a_2+{\rm O}(1)\big)\,
\frac{\partial}{\partial u_2}
+
\big(2\,a_2+{\rm O}(1)\big)\,
\frac{\partial}{\partial u_3},
\]
and at the origin:
\[
\isqrt\,
\big[\mathcal{L},\overline{\mathcal{L}}\big]
\Big\vert_0
=
2a_1\,\frac{\partial}{\partial u_1}
\bigg\vert_0
+
2a_2\,\frac{\partial}{\partial u_2}
\bigg\vert_0
+
2a_3\,\frac{\partial}{\partial u_2}
\bigg\vert_0.
\]

Since the three vectors:
\[
\mathcal{L}\big\vert_0, 
\ \ \ \ \ \ \ \ \ \ \ \ \
\overline{\mathcal{L}}\big\vert_0,
\ \ \ \ \ \ \ \ \ \ \ \ \
\big[\mathcal{L},\overline{\mathcal{L}}\big] 
\Big\vert_0
\]
must in particular by hypothesis necessarily be $\C$-linearly independent
at the origin, one has:
\[
a_1\neq 0,
\ \ \ \ \ \ \ \ \ \ \ \ \ \ \ \ \ \ \
\text{\rm or}
\ \ \ \ \ \ \ \ \ \ \ \ \ \ \ \ \ \ \
a_2\neq 0,
\ \ \ \ \ \ \ \ \ \ \ \ \ \ \ \ \ \ \
\text{\rm or}
\ \ \ \ \ \ \ \ \ \ \ \ \ \ \ \ \ \ \
a_3\neq 0.
\]

Say:
\[
a_1\neq 0.
\]
Then do if necessary:
\[
w_1
\,\longmapsto\,
-\,w_1
\]
to make:
\[
a_1>0,
\]
and do:
\[
z
\,\longmapsto\,
\sqrt{a_1}\,z,
\]
to make:
\[
a_1
=
1.
\]

Writing out now the third-order $(z, \overline{z})$-terms, one
therefore comes to:
\[
\aligned
v_1
&
=
z\overline{z}
+
\alpha_1\,z^2\overline{z}
+
\overline{\alpha}_1\,z\overline{z}^2
+
z\overline{z}\,
{\rm O}_2\big(z,\overline{z}\big)
+
z\overline{z}\,{\rm O}_1(u_1)
+
z\overline{z}\,{\rm O}_1(u_2)
+
z\overline{z}\,{\rm O}_1(u_3),
\\
v_2
&
=
a_2\,z\overline{z}
+
\alpha_2\,z^2\overline{z}
+
\overline{\alpha}_2\,z\overline{z}^2
+
z\overline{z}\,
{\rm O}_2\big(z,\overline{z}\big)
+
z\overline{z}\,{\rm O}_1(u_1)
+
z\overline{z}\,{\rm O}_1(u_2)
+
z\overline{z}\,{\rm O}_1(u_3),
\\
v_3
&
=
a_3\,z\overline{z}
+
\alpha_3\,z^2\overline{z}
+
\overline{\alpha}_3\,z\overline{z}^2
+
z\overline{z}\,
{\rm O}_2\big(z,\overline{z}\big)
+
z\overline{z}\,{\rm O}_1(u_1)
+
z\overline{z}\,{\rm O}_1(u_2)
+
z\overline{z}\,{\rm O}_1(u_3),
\endaligned
\]
with of course:
\[
\alpha_1
\in
\C,
\ \ \ \ \ \ \ \ \ \ \ \ \ \ \ \ \ \ \
\alpha_2
\in
\C,
\ \ \ \ \ \ \ \ \ \ \ \ \ \ \ \ \ \ \
\alpha_3
\in
\C.
\]

Changing:
\[
\aligned
w_2
&
\,\longmapsto\,
w_2-a_2\,w_1,
\\
w_3
&
\,\longmapsto\,
w_3-a_3\,w_1,
\endaligned
\]
one makes:
\[
\aligned
a_2
&
=
0,
\\
a_3
&
=
0.
\endaligned
\]
In the process, $\alpha_2$, $\alpha_3$ are modified.

Changing:
\[
z
\,\longmapsto\,
z+\alpha_1\,z^2,
\]
one makes:
\[
\alpha_1
=
0.
\]

One thus arrives at:
\[
\aligned
v_1
&
=
z\overline{z}
\ \ \ \ \ \ \ \ \ \ \ \ \ \ \ \ \ \ \ \ \ \,
+
z\overline{z}\,
{\rm O}_2\big(z,\overline{z}\big)
+
z\overline{z}\,{\rm O}_1(u_1)
+
z\overline{z}\,{\rm O}_1(u_2)
+
z\overline{z}\,{\rm O}_1(u_3),
\\
v_2
&
=
\alpha_2\,z^2\overline{z}
+
\overline{\alpha}_2\,z\overline{z}^2
+
z\overline{z}\,
{\rm O}_2\big(z,\overline{z}\big)
+
z\overline{z}\,{\rm O}_1(u_1)
+
z\overline{z}\,{\rm O}_1(u_2)
+
z\overline{z}\,{\rm O}_1(u_3),
\\
v_3
&
=
\alpha_3\,z^2\overline{z}
+
\overline{\alpha}_3\,z\overline{z}^2
+
z\overline{z}\,
{\rm O}_2\big(z,\overline{z}\big)
+
z\overline{z}\,{\rm O}_1(u_1)
+
z\overline{z}\,{\rm O}_1(u_2)
+
z\overline{z}\,{\rm O}_1(u_3).
\endaligned
\]

Take again the approximant for the denominator
common to $A_1$, $A_2$, $A_3$:
\[
\aligned
\left\vert\!
\begin{array}{ccc}
\isqrt+\varphi_{1,u_1} & \varphi_{1,u_2} & \varphi_{1,u_3}
\\
\varphi_{2,u_1} & \isqrt+\varphi_{2,u_2} & \varphi_{2,u_3}
\\
\varphi_{3,u_1} & \varphi_{3,u_2} & \isqrt+\varphi_{3,u_3}
\end{array}
\!\right\vert
&
=
\left\vert\!
\begin{array}{ccc}
\isqrt+z\overline{z}{\rm O}(0) & z\overline{z}{\rm O}(0) & z\overline{z}{\rm O}(0)
\\
z\overline{z}{\rm O}(0) & \isqrt+z\overline{z}{\rm O}(0) & z\overline{z}{\rm O}(0)
\\
z\overline{z}{\rm O}(0) & z\overline{z}{\rm O}(0) & \isqrt+z\overline{z}{\rm O}(0)
\end{array}
\!\right\vert
\\
&
=
\isqrt^3
+
z\overline{z}{\rm O}(0).
\endaligned
\]

Next, a deeper approximant for the numerator of $A_1$ is:
\[
\aligned
\left\vert\!
\begin{array}{ccc}
-\varphi_{1,z} & \varphi_{1,u_2} & \varphi_{1,u_3}
\\
-\varphi_{2,z} & \isqrt+\varphi_{2,u_2} & \varphi_{2,u_3}
\\
-\varphi_{3,z} & \varphi_{3,u_2} & \isqrt+\varphi_{3,u_3}
\end{array}
\!\right\vert
&
=
\left\vert\!
\begin{array}{ccc}
\substack{
-\overline{z}+{\rm O}_3(z,\overline{z})\\
+\overline{z}{\rm O}_1(u_1)
+\overline{z}{\rm O}_1(u_2)+\overline{z}{\rm O}_1(u_3)}
&
\substack{z\overline{z}{\rm O}(0)}
&
\substack{z\overline{z}{\rm O}(0)}\bigskip
\\
\substack{
-2\alpha_2z\overline{z}-\overline{\alpha}_2\overline{z}^2
+{\rm O}_3(z,\overline{z})+\\
+\overline{z}{\rm O}_1(u_1)+\overline{z}{\rm O}_1(u_2)
+\overline{z}{\rm O}_1(u_3)}
&
\isqrt+\substack{z\overline{z}{\rm O}(0)}
&
\substack{z\overline{z}{\rm O}(0)}\bigskip
\\
\substack{
-2\alpha_3z\overline{z}-\overline{\alpha}_3\overline{z}^2
+{\rm O}_3(z,\overline{z})+\\
+\overline{z}{\rm O}_1(u_1)+\overline{z}{\rm O}_1(u_2)
+\overline{z}{\rm O}_1(u_3)}
&
\substack{z\overline{z}{\rm O}(0)}
&
\isqrt+\substack{z\overline{z}{\rm O}(0)}
\end{array}
\!\right\vert
\endaligned
\]
that is:
\[
\aligned
&
\isqrt^2\big(
-\,
\overline{z}+
{\rm O}_3(z,\overline{z})
+
\overline{z}\,{\rm O}_1(u_1)
+
\overline{z}\,{\rm O}_1(u_2)
+
\overline{z}\,{\rm O}_1(u_3)
\big),
\endaligned
\]
so that:
\[
\aligned
A_1
&
=
\frac{
\isqrt^2\big(
-\,\overline{z}+
{\rm O}_3(z,\overline{z})
+
\overline{z}\,{\rm O}_1(u_1)
+
\overline{z}\,{\rm O}_1(u_2)
+
\overline{z}\,{\rm O}_1(u_3)
\big)}{
\isqrt^3+z\overline{z}{\rm O}(0)}
\\
&
=
\isqrt\,\big(
\overline{z}
+
{\rm O}_3(z,\overline{z})
+
\overline{z}\,{\rm O}_1(u_1)
+
\overline{z}\,{\rm O}_1(u_2)
+
\overline{z}\,{\rm O}_1(u_3)
\big).
\endaligned
\]

Next, a deeper approximant for the numerator of $A_2$ is:
\[
\aligned
\left\vert\!
\begin{array}{ccc}
\isqrt+\varphi_{1,u_2} & -\varphi_{1,z} & \varphi_{1,u_3}
\\
\varphi_{2,u_1} & -\varphi_{2,z} &  \varphi_{2,u_3}
\\
\varphi_{3,u_1} & -\varphi_{3,z} & \isqrt+\varphi_{3,u_3}
\end{array}
\!\right\vert
&
=
\left\vert\!
\begin{array}{ccc}
\substack{\isqrt+z\overline{z}{\rm O}(0)}
&
\substack{
-\overline{z}+{\rm O}_3(z,\overline{z})\\
+\overline{z}{\rm O}_1(u_1)
+\overline{z}{\rm O}_1(u_2)+\overline{z}{\rm O}_1(u_3)}
&
\substack{z\overline{z}{\rm O}(0)}\bigskip
\\
\substack{z\overline{z}{\rm O}(0)}
&
\substack{
-2\alpha_2z\overline{z}-\overline{\alpha}_2\overline{z}^2
+{\rm O}_3(z,\overline{z})+\\
+\overline{z}{\rm O}_1(u_1)+\overline{z}{\rm O}_1(u_2)
+\overline{z}{\rm O}_1(u_3)}
&
\substack{z\overline{z}{\rm O}(0)}\bigskip
\\
\substack{z\overline{z}{\rm O}(0)}
&
\substack{
-2\alpha_3z\overline{z}-\overline{\alpha}_3\overline{z}^2
+{\rm O}_3(z,\overline{z})+\\
+\overline{z}{\rm O}_1(u_1)+\overline{z}{\rm O}_1(u_2)
+\overline{z}{\rm O}_1(u_3)}
&
\isqrt+\substack{z\overline{z}(0)}
\end{array}
\!\right\vert,
\endaligned
\]
that is:
\[
\aligned
&
\isqrt^2\big(
-\,2\alpha_2\,z\overline{z}-\overline{\alpha}_2\,\overline{z}^2
+
{\rm O}_3(z,\overline{z})
+
\overline{z}{\rm O}_1(u_1)
+
\overline{z}{\rm O}_1(u_2)
+
\overline{z}{\rm O}_1(u_3)
\big),
\endaligned
\]
so that:
\[
\aligned
A_2
&
=
\frac{
\isqrt^2\big(
-\,2\alpha_2\,z\overline{z}-\overline{\alpha}_2\,\overline{z}^2
+
{\rm O}_3(z,\overline{z})
+
\overline{z}{\rm O}_1(u_1)
+
\overline{z}{\rm O}_1(u_2)
+
\overline{z}{\rm O}_1(u_3)
\big)}{
\isqrt^3+z\overline{z}{\rm O}(0)}
\\
&
=
\isqrt\,\big(
2\alpha_2\,z\overline{z}+\overline{\alpha}_2\,\overline{z}^2
+
{\rm O}_3(z,\overline{z})
+
\overline{z}{\rm O}_1(u_1)
+
\overline{z}{\rm O}_1(u_2)
+
\overline{z}{\rm O}_1(u_3)
\big).
\endaligned
\]

Quite similarly:
\[
A_3
=
\isqrt\,\big(
2\alpha_3\,z\overline{z}+\overline{\alpha}_3\,\overline{z}^2
+
{\rm O}_3(z,\overline{z})
+
\overline{z}{\rm O}_1(u_1)
+
\overline{z}{\rm O}_1(u_2)
+
\overline{z}{\rm O}_1(u_3)
\big).
\]

It follows that:
\[
\aligned
\mathcal{L}
&
=
\frac{\partial}{\partial z}
+
\isqrt\,
\big(\overline{z}
+
\overline{z}\,{\rm O}_1(u_1,u_2,u_3)
+
{\rm O}(3)\big)\,
\frac{\partial}{\partial u_1}
+
\\
&
\ \ \ \ \ \ \ \ \ \ \,
+
\isqrt\,
\big(
2\alpha_2\,z\overline{z}
+
\overline{\alpha}_2\,\overline{z}^2
+
\overline{z}\,{\rm O}_1(u_1,u_2,u_3)
+
{\rm O}(3)\big)\,
\frac{\partial}{\partial u_2},
\\
&
\ \ \ \ \ \ \ \ \ \ \,
+
\isqrt\,
\big(
2\alpha_3\,z\overline{z}
+
\overline{\alpha}_3\,\overline{z}^2
+
\overline{z}\,{\rm O}_1(u_1,u_2,u_3)
+
{\rm O}(3)\big)\,
\frac{\partial}{\partial u_3},
\endaligned
\]
\[
\aligned
\overline{\mathcal{L}}
&
=
\frac{\partial}{\partial\overline{z}}
-
\isqrt\,
\big(
z
+
z\,{\rm O}_1(u_1,u_2,u_3)
+
{\rm O}(3)\big)\,
\frac{\partial}{\partial u_1}
-
\\
&
\ \ \ \ \ \ \ \ \ \ \,
-
\isqrt\,
\big(
2\overline{\alpha}_2\,z\overline{z}
+
\alpha_2\,z^2
+
z\,{\rm O}_1(u_1,u_2,u_3)
+
{\rm O}(3)\big)\,
\frac{\partial}{\partial u_2},
\\
&
\ \ \ \ \ \ \ \ \ \ \,
-
\isqrt\,
\big(
2\overline{\alpha}_3\,z\overline{z}
+
\alpha_3\,z^2
+
z\,{\rm O}_1(u_1,u_2,u_3)
+
{\rm O}(3)\big)\,
\frac{\partial}{\partial u_3},
\endaligned
\]
whence:
\[
\aligned
\big[\mathcal{L},\overline{\mathcal{L}}\big]
&
=
-\,\isqrt\,
\big(
2
+
{\rm O}_1(u_1,u_2,u_3)
+
{\rm O}(2)\big)\,
\frac{\partial}{\partial u_1}
-
\\
&
\ \ \ \ \
-
\isqrt\,
\big(
4\alpha_2\,z+4\overline{\alpha}_2\,\overline{z}
+
{\rm O}_1(u_1,u_2,u_3)
+
{\rm O}(2)\big)\,
\frac{\partial}{\partial u_2},
\\
&
\ \ \ \ \
-
\isqrt\,
\big(
4\alpha_3\,z+4\overline{\alpha}_3\,\overline{z}
+
{\rm O}_1(u_1,u_2,u_3)
+
{\rm O}(2)\big)\,
\frac{\partial}{\partial u_3},
\endaligned
\]
whence further:
\[
\aligned
\big[\mathcal{L},\,\big[\mathcal{L},\overline{\mathcal{L}}\big]\big]
&
=
{\rm O}(1)\,\frac{\partial}{\partial u_1}
-
\isqrt\,\big(4\alpha_2+{\rm O}(1)\big)\,
\frac{\partial}{\partial u_2}
-
\isqrt\,\big(4\alpha_3+{\rm O}(1)\big)\,
\frac{\partial}{\partial u_3},
\\
\big[\overline{\mathcal{L}},\,
\big[\mathcal{L},\overline{\mathcal{L}}\big]\big]
&
=
{\rm O}(1)\,\frac{\partial}{\partial u_1}
-
\isqrt\,\big(4\overline{\alpha}_2+{\rm O}(1)\big)\,
\frac{\partial}{\partial u_2}
-
\isqrt\,\big(4\overline{\alpha}_3+{\rm O}(1)\big)\,
\frac{\partial}{\partial u_3},
\endaligned
\]
and lastly, at the origin:
\[
\aligned
\mathcal{L}\big\vert_0
&
=
\frac{\partial}{\partial z}
\bigg\vert_0,
\\
\overline{\mathcal{L}}\big\vert_0
&
=
\frac{\partial}{\partial\overline{z}}
\bigg\vert_0,
\\
\big[\mathcal{L},\overline{\mathcal{L}}\big]
\Big\vert_0
&
=
-\,2\,\isqrt\,
\frac{\partial}{\partial u_1}
\bigg\vert_0,
\\
\big[\mathcal{L},\,\big[\mathcal{L},\overline{\mathcal{L}}\big]\big]
\Big\vert_0
&
=
-\,4\,\isqrt\,\alpha_2\,
\frac{\partial}{\partial u_2}
\bigg\vert_0
-4\,\isqrt\,\alpha_3\,
\frac{\partial}{\partial u_3}
\bigg\vert_0,
\\
\big[\overline{\mathcal{L}},\,
\big[\mathcal{L},\overline{\mathcal{L}}\big]\big]
\Big\vert_0
&
=
-\,4\,\isqrt\,\overline{\alpha}_2\,
\frac{\partial}{\partial u_2}
\bigg\vert_0
-4\,\isqrt\,\overline{\alpha}_3\,
\frac{\partial}{\partial u_3}
\bigg\vert_0.
\endaligned
\]

These four vectors should by hypothesis constitute a basis for:
\[
\C\otimes_\R T_0M
\,=\,
\C\frac{\partial}{\partial z}
\bigg\vert_0
\oplus
\C\frac{\partial}{\partial\overline{z}}
\bigg\vert_0
\oplus
\C\frac{\partial}{\partial u_1}
\bigg\vert_0
\oplus
\C\frac{\partial}{\partial u_2}
\bigg\vert_0
\oplus
\C\frac{\partial}{\partial u_3}
\bigg\vert_0,
\] 
having dimension ${\bf 5}$, which necessitates the nonzeroness:
\[
0
\neq
\left\vert\!
\begin{array}{ccc}
2 & 0 & 0
\\
0 & 4\alpha_2 & 4\alpha_3
\\
0 & 4\overline{\alpha}_2 & 4\overline{\alpha}_3
\end{array}
\!\right\vert.
\]

At least:
\[
\alpha_2
\neq
0.
\]

Using a complex dilation:
\[
z
\,\longmapsto\,
\lambda\,z,
\]
with $\lambda$ solution of:
\[
1
=
\alpha_2\,\lambda^2\,\overline{\lambda},
\]
one makes:
\[
\alpha_2
=
1,
\]
so that the second equation receives the form announced:
\[
v_2
=
z^2\overline{z}
+
z\overline{z}^2
+
z\overline{z}\,{\rm O}_2\big(z,\overline{z}\big)
+
z\overline{z}\,{\rm O}_1(u_1)
+
z\overline{z}\,{\rm O}_1(u_2)
+
z\overline{z}\,{\rm O}_1(u_3).
\]

In the process, the first equation is changed to:
\[
v_1
=
\lambda\overline{\lambda}\,z\overline{z}
+
z\overline{z}\,{\rm O}_2\big(z,\overline{z}\big)
+
z\overline{z}\,{\rm O}_1(u_1)
+
z\overline{z}\,{\rm O}_1(u_2)
+
z\overline{z}\,{\rm O}_1(u_3),
\]
but one can yet replace:
\[
w_1
\,\longmapsto\,
\frac{w_1}{\lambda\overline{\lambda}},
\]
to keep $1\cdot z \overline{ z}$.

Writing:
\[
\alpha_3
=
a_3
+
\isqrt\,b_3,
\]
replacing:
\[
w_3
\,\longmapsto\,
w_3
-
a_3\,w_2,
\]
one makes:
\[
a_3
=
0.
\]

Lastly:
\[
0
\neq
\left\vert\!
\begin{array}{ccc}
2 & 0 & 0
\\
0 & 4 & 4\isqrt\,b_3
\\
0 & 4 & -4\isqrt\,b_3
\end{array}
\!\right\vert.
\]
means:
\[
b_3
\neq
0,
\]
and doing:
\[
w_3
\,\longmapsto\,
\frac{w_3}{b_3},
\]
one arrives at the announced form for the third graphing equation.
\endproof

\noindent{\bf Scholium.}
{\em In such elementarily normalized coordinates, one has
the diagonal normalization at the origin:}
\[
\aligned
\mathcal{L}\big\vert_0
&
=
\frac{\partial}{\partial z}
\bigg\vert_0,
\\
\overline{\mathcal{L}}\big\vert_0
&
=
\frac{\partial}{\partial\overline{z}}
\bigg\vert_0,
\\
\big[\mathcal{L},\overline{\mathcal{L}}\big]
\Big\vert_0
&
=
-\,2\,\isqrt\,
\frac{\partial}{\partial u_1}
\bigg\vert_0,
\\
\big[\mathcal{L},\,\big[\mathcal{L},\overline{\mathcal{L}}\big]\big]
\Big\vert_0
&
=
-\,4\,\isqrt\,
\frac{\partial}{\partial u_2}
\bigg\vert_0
+
4\,
\frac{\partial}{\partial u_3}
\bigg\vert_0,
\\
\big[\overline{\mathcal{L}},\,
\big[\mathcal{L},\overline{\mathcal{L}}\big]\big]
\Big\vert_0
&
=
-\,4\,\isqrt\,
\frac{\partial}{\partial u_2}
\bigg\vert_0
-
4\,
\frac{\partial}{\partial u_3}
\bigg\vert_0,
\endaligned
\]
{\em which conveniently fixes ideas 
when performing explicitly the Cartan equivalence
procedure.\qed}

\medskip

Taking zero remainders, one obtains Beloshapka's second cubic:
\[
\boxed{\,\,
{\text{\footnotesize\sf Model 
$\text{\sf (III)}_{\text{\sf 1}}$:}}
\ \ \ \ \ \ \ \ \ 
\aligned
v_1
&
=
z\overline{z},\,\,
\\
v_2
&
=
z^2\overline{z}+z\overline{z}^2,
\\
v_3
&
=
\isqrt\,\big(z^2\overline{z}-z\overline{z}^2\big).
\endaligned
\rule[-4pt]{0pt}{27pt}}
\]


\bigskip

\section{\sf General class $\text{\sf III}_{\text{\sf 2}}$}
\label{general-class-III-2}
\HEAD{\ref{general-class-III-2}.~General class 
$\text{\sf III}_{\text{\sf 2}}$}{
Jo\"el {\sc Merker} (Paris-Sud)}

\medskip

\medskip\noindent{\bf Proposition.}
{\em A local real analytic $5$-dimensional
CR-generic submanifold passing through the origin:}
\[
0
\,\in\,
M^5
\,\subset\,
\C^4
\]
{\em which belongs to the general Class $\text{\sf III}_{
\text{\sf 1}}$, namely such that,
for any local vector field generator $\mathcal{ L}$ of $T^{1, 0}M$, 
one has at every point:}
\[
\aligned
{\bf 3}
&
\,=\,
\rank_\C\Big(
\Big\{
\mathcal{L},\,\overline{\mathcal{L}},\,\,
\big[\mathcal{L},\overline{\mathcal{L}}\big]
\Big\}
\Big),
\\
{\bf 4}
&
\,=\,
\rank_\C\Big(
\Big\{
\mathcal{L},\,\overline{\mathcal{L}},\,\,
\big[\mathcal{L},\overline{\mathcal{L}}\big],\,\,
\big[\mathcal{L},\,\big[\mathcal{L},\overline{\mathcal{L}}\big]\big]
\Big\}
\Big),
\\
{\bf 4}
&
\,=\,
\rank_\C\Big(
\Big\{
\mathcal{L},\,\overline{\mathcal{L}},\,\,
\big[\mathcal{L},\overline{\mathcal{L}}\big],\,\,
\big[\mathcal{L},\,\big[\mathcal{L},\overline{\mathcal{L}}\big]\big],\,\,
\big[\overline{\mathcal{L}},\,\big[\mathcal{L},
\overline{\mathcal{L}}\big]\big]
\Big\}
\Big),
\\
{\bf 5}
&
\,=\,
\rank_\C\Big(
\Big\{
\mathcal{L},\,\overline{\mathcal{L}},\,\,
\big[\mathcal{L},\overline{\mathcal{L}}\big],\,\,
\big[\mathcal{L},\,\big[\mathcal{L},\overline{\mathcal{L}}\big]\big],\,\,
\big[\mathcal{L},\,\big[\mathcal{L},\,\big[\mathcal{L},
\overline{\mathcal{L}}\big]\big]\big]
\Big\}
\Big),
\endaligned
\]
{\em may always be represented, in suitable local holomorphic coordinates:}
\[
(z,w_1,w_2,w_3)
\in
\C^4
\]
{\em by three specific real analytic equations of the form:}
\[
\!\!\!\!\!\!\!\!\!\!\!\!\!\!\!\!\!\!\!\!
\boxed{\,\,
\underline{\text{\footnotesize\sf (III)}_{
\text{\tiny\sf 2}}:}
\ \ \ \ \
\aligned
v_1
&
=
z\overline{z}
+
c_1\,z^2\overline{z}^2
+
z\overline{z}\,{\rm O}_3\big(z,\overline{z}\big)
+
z\overline{z}\,u_1\,{\rm O}_1\big(z,\overline{z},u_1\big)
+
\\
&
\ \ \ \ \ \ \ \ \ \ \ \ \ \ \ \ \ \ \ \ \ \ 
+
z\overline{z}\,u_2\,{\rm O}_1\big(z,\overline{z},u_1,u_2\big)
+
z\overline{z}\,u_3\,{\rm O}_1\big(z,\overline{z},u_1,u_2,u_3\big),
\\
v_2
&
=
z^2\overline{z}
+
z\overline{z}^2
+
\beta_2\,z^2\overline{z}
+
\overline{\beta}_2\,z\overline{z}^2
+
c_2\,z^2\overline{z}^2
+
z\overline{z}\,{\rm O}_3\big(z,\overline{z}\big)
+
z\overline{z}\,u_1\,{\rm O}_1\big(z,\overline{z},u_1\big)
+
\\
&
\ \ \ \ \ \ \ \ \ \ \ \ \ \ \ \ \ \ \ \ \ \ 
+
z\overline{z}\,u_2\,{\rm O}_1\big(z,\overline{z},u_1,u_2\big)
+
z\overline{z}\,u_3\,{\rm O}_1\big(z,\overline{z},u_1,u_2,u_3\big),
\\
v_3
&
=
2\,z^3\overline{z}
+
2\,z\overline{z}^3
+
3\,z^2\overline{z}^2
+
z\overline{z}\,{\rm O}_3\big(z,\overline{z}\big)
+
z\overline{z}\,u_1\,{\rm O}_1\big(z,\overline{z},u_1\big)
+
\\
&
\ \ \ \ \ \ \ \ \ \ \ \ \ \ \ \ \ \ \ \ \ \ \ \ \ \ \ \ \ \ \ \ \ \ \ 
\ \ \ \ \ \ 
+
z\overline{z}\,u_2\,{\rm O}_1\big(z,\overline{z},u_1,u_2\big)
+
z\overline{z}\,u_3\,{\rm O}_1\big(z,\overline{z},u_1,u_2,u_3\big).
\endaligned
\rule[-5pt]{0pt}{27pt}}
\]

\proof
Choose first coordinates $(z, w_1, w_2, w_3)$ so that $M$ is:
\[
\aligned
v_1
&
=
\varphi_1\big(z,\overline{z},u_1,u_2,u_3\big)
=
a_1\,z\overline{z}
+
z\overline{z}\,{\rm O}_1\big(z,\overline{z},u_1,u_2,u_3\big),
\\
v_2
&
=
\varphi_2\big(z,\overline{z},u_1,u_2,u_3\big)
=
a_2\,z\overline{z}
+
z\overline{z}\,{\rm O}_1\big(z,\overline{z},u_1,u_2,u_3\big),
\\
v_3
&
=
\varphi_3\big(z,\overline{z},u_1,u_2,u_3\big)
=
a_3\,z\overline{z}
+
z\overline{z}\,{\rm O}_1\big(z,\overline{z},u_1,u_2,u_3\big),
\endaligned
\]
with:
\[
a_1,\,\,a_2,\,\,a_3
\,\,\in\,\R.
\]
As for the general class $\text{\sf III}_{\text{\sf 1}}$,
one sees that at least one of $a_1$, $a_2$, 
$a_3$ must be nonzero, say:
\[
a_1
\neq
0.
\]

Elementary linear transformations along coordinate axes yield:
\[
\aligned
v_1
&
=
z\overline{z}
+
z\overline{z}\,{\rm O}_1\big(z,\overline{z},u_1,u_2,u_3\big),
\\
v_2
&
=
\ \ \ \ \ \ \ \ \,
z\overline{z}\,{\rm O}_1\big(z,\overline{z},u_1,u_2,u_3\big),
\\
v_3
&
=
\ \ \ \ \ \ \ \ \,
z\overline{z}\,{\rm O}_1\big(z,\overline{z},u_1,u_2,u_3\big).
\endaligned
\]

To better organize these general remainders: 
\[
z\overline{z}\,{\rm O}_1\big(z,\overline{z},u_1,u_2,u_3\big),
\]
observe that every local converging power series function:
\[
F\big({\sf x}_1,{\sf x}_2,\dots,{\sf x}_\NN\big)
=
\sum_{\alpha_1,\alpha_2,\dots,\alpha_\NN\,\in\,\N}\,
F_{\alpha_1,\alpha_2,\dots,\alpha_\NN}\,
\big({\sf x}_1\big)^{\alpha_1}\,
\big({\sf x}_2\big)^{\alpha_2}\,
\cdots\,
\big({\sf x}_\NN\big)^{\alpha_\NN}
\]
can always be written under the progressively factorized form:
\[
F\big({\sf x}_1,{\sf x}_2,\dots,{\sf x}_\NN\big)
=
F_0
+
{\sf x}_1\,F_1\big({\sf x}_1\big)
+
{\sf x}_2\,F\big({\sf x}_1,{\sf x}_2\big)
+\cdots+
{\sf x}_\NN\,
F_\NN\big({\sf x}_1,{\sf x}_2,\dots,{\sf x}_\NN\big).
\]

Hence:
\[
\!\!\!\!\!\!\!\!\!\!\!\!\!\!\!\!\!\!\!\!
\aligned
v_1
&
=
z\overline{z}
+
z\overline{z}\,{\rm O}_1\big(z,\overline{z}\big)
+
z\overline{z}\,u_1\,{\rm O}_0\big(z,\overline{z},u_1\big)
+
z\overline{z}\,u_2\,{\rm O}_0\big(z,\overline{z},u_1,u_2\big)
+
z\overline{z}\,u_3\,{\rm O}_0\big(z,\overline{z},u_1,u_2,u_3\big),
\\
v_2
&
=
\ \ \ \ \ \ \ \ \,
z\overline{z}\,{\rm O}_1\big(z,\overline{z}\big)
+
z\overline{z}\,u_1\,{\rm O}_0\big(z,\overline{z},u_1\big)
+
z\overline{z}\,u_2\,{\rm O}_0\big(z,\overline{z},u_1,u_2\big)
+
z\overline{z}\,u_3\,{\rm O}_0\big(z,\overline{z},u_1,u_2,u_3\big),
\\
v_3
&
=
\ \ \ \ \ \ \ \ \,
z\overline{z}\,{\rm O}_1\big(z,\overline{z}\big)
+
z\overline{z}\,u_1\,{\rm O}_0\big(z,\overline{z},u_1\big)
+
z\overline{z}\,u_2\,{\rm O}_0\big(z,\overline{z},u_1,u_2\big)
+
z\overline{z}\,u_3\,{\rm O}_0\big(z,\overline{z},u_1,u_2,u_3\big).
\endaligned
\]

In the right-hand sides, let us emphasize three monomials
on each line:
\[
\aligned
&
d_1\,z\overline{z}\,u_1
+
e_1\,z\overline{z}\,u_2
+
f_1\,z\overline{z}\,u_3,
\\
&
d_2\,z\overline{z}\,u_1
+
e_2\,z\overline{z}\,u_2
+
f_2\,z\overline{z}\,u_3,
\\
&
d_3\,z\overline{z}\,u_1
+
e_3\,z\overline{z}\,u_2
+
f_3\,z\overline{z}\,u_3,
\endaligned
\]
so that:
\[
\aligned
v_1
&
=
z\overline{z}
+
z\overline{z}\,{\rm O}_1\big(z,\overline{z}\big)
+
d_1\,z\overline{z}\,u_1
+
z\overline{z}\,u_1\,{\rm O}_1\big(z,\overline{z},u_1\big)
+
\\
&
\ \ \ \ \ \ \ \ \ \ \ \ \ \ \ \ \ \ \ \ \ \ \ \ \ \ \ \ \ \ \ \ \
+
e_1\,z\overline{z}\,u_2
+
z\overline{z}\,u_2\,{\rm O}_1\big(z,\overline{z},u_1,u_2\big)
+
\\
&
\ \ \ \ \ \ \ \ \ \ \ \ \ \ \ \ \ \ \ \ \ \ \ \ \ \ \ \ \ \ \ \ \
+
f_1\,z\overline{z}\,u_3
+
z\overline{z}\,u_3\,{\rm O}_1\big(z,\overline{z},u_1,u_2,u_3\big),
\\
v_2
&
=
z\overline{z}
+
z\overline{z}\,{\rm O}_1\big(z,\overline{z}\big)
+
d_2\,z\overline{z}\,u_1
+
z\overline{z}\,u_1\,{\rm O}_1\big(z,\overline{z},u_1\big)
+
\\
&
\ \ \ \ \ \ \ \ \ \ \ \ \ \ \ \ \ \ \ \ \ \ \ \ \ \ \ \ \ \ \ \ \
+
e_2\,z\overline{z}\,u_2
+
z\overline{z}\,u_2\,{\rm O}_1\big(z,\overline{z},u_1,u_2\big)
+
\\
&
\ \ \ \ \ \ \ \ \ \ \ \ \ \ \ \ \ \ \ \ \ \ \ \ \ \ \ \ \ \ \ \ \
+
f_2\,z\overline{z}\,u_3
+
z\overline{z}\,u_3\,{\rm O}_1\big(z,\overline{z},u_1,u_2,u_3\big),
\\
v_3
&
=
z\overline{z}
+
z\overline{z}\,{\rm O}_1\big(z,\overline{z}\big)
+
d_3\,z\overline{z}\,u_1
+
z\overline{z}\,u_1\,{\rm O}_1\big(z,\overline{z},u_1\big)
+
\\
&
\ \ \ \ \ \ \ \ \ \ \ \ \ \ \ \ \ \ \ \ \ \ \ \ \ \ \ \ \ \ \ \ \
+
e_3\,z\overline{z}\,u_2
+
z\overline{z}\,u_2\,{\rm O}_1\big(z,\overline{z},u_1,u_2\big)
+
\\
&
\ \ \ \ \ \ \ \ \ \ \ \ \ \ \ \ \ \ \ \ \ \ \ \ \ \ \ \ \ \ \ \ \
+
f_3\,z\overline{z}\,u_3
+
z\overline{z}\,u_3\,{\rm O}_1\big(z,\overline{z},u_1,u_2,u_3\big).
\endaligned
\]

\medskip\noindent{\bf Assertion.}
{\em An appropriate choice of real constants in the biholomorphic
change of coordinates:}
\[
\aligned
z
&
=
z',
\\
w_1
&
=
w_1'
+
l_1'\,w_1'w_1'
+
m_1'\,w_1'w_2'
+
n_1'\,w_1'w_3',
\\
w_2
&
=
w_2'
+
l_2'\,w_1'w_1'
+
m_2'\,w_1'w_2'
+
n_2'\,w_1'w_3',
\\
w_3
&
=
w_3'
+
l_3'\,w_1'w_1'
+
m_3'\,w_1'w_2'
+
n_3'\,w_1'w_3',
\endaligned
\]
{\em annihilates the 9 new constants:}
\[
\aligned
0
&
=
d_1'
=
e_1'
=
f_1',
\\
0
&
=
d_2'
=
e_2'
=
f_2',
\\
0
&
=
d_3'
=
e_3'
=
f_3'.
\endaligned
\]

\proof
Indeed, the real and imaginary parts are:
\[
\aligned
u_1
&
=
u_1'
+
l_1'\big(u_1'u_1'-v_1'v_1'\big)
+
m_1'\big(u_1'u_2'-v_1'v_2'\big)
+
n_1'\big(u_1'u_3'-v_1'v_3'\big),
\\
v_1
&
=
v_1'
+
l_1'\big(u_1'v_1'+u_1'v_1'\big)
+
m_1'\big(u_1'v_2'+u_2'v_1'\big)
+
n_1'\big(u_1'v_3'+u_3'v_1'\big),
\endaligned
\]
\[
\aligned
u_2
&
=
u_2'
+
l_2'\big(u_1'u_1'-v_1'v_1'\big)
+
m_2'\big(u_1'u_2'-v_1'v_2'\big)
+
n_2'\big(u_1'u_3'-v_1'v_3'\big),
\\
v_2
&
=
v_2'
+
l_2'\big(u_1'v_1'+u_1'v_1'\big)
+
m_2'\big(u_1'v_2'+u_2'v_1'\big)
+
n_2'\big(u_1'v_3'+u_3'v_1'\big),
\endaligned
\]
\[
\aligned
u_3
&
=
u_3'
+
l_3'\big(u_1'u_1'-v_1'v_1'\big)
+
m_3'\big(u_1'u_2'-v_1'v_2'\big)
+
n_3'\big(u_1'u_3'-v_1'v_3'\big),
\\
v_3
&
=
v_3'
+
l_3'\big(u_1'v_1'+u_1'v_1'\big)
+
m_3'\big(u_1'v_2'+u_2'v_1'\big)
+
n_3'\big(u_1'v_3'+u_3'v_1'\big).
\endaligned
\]
Disregarding terms of order $4$ and higher, 
appropriate replacements
in the three equations of $M$ give:
\[
\aligned
&
v_1'
+
l_1'\big(u_1'v_1'+u_1'v_1'\big)
+
m_1'\big(u_1'v_2'+u_2'v_1'\big)
+
n_1'\big(u_1'v_3'+u_3'v_1'\big)
\equiv
\\
&
\ \ \ \ \ \ \ \
\equiv
z'\overline{z}'
+
z'\overline{z}'\,{\rm O}_1\big(z',\overline{z}'\big)
+
d_1\,z'\overline{z}'\,u_1'
+
e_1\,z'\overline{z}'\,u_2'
+
f_1\,z'\overline{z}'\,u_3'
\,\,\,
\mod\,{\rm O}(4),
\endaligned
\]
\[
\aligned
&
v_2'
+
l_2'\big(u_1'v_1'+u_1'v_1'\big)
+
m_2'\big(u_1'v_2'+u_2'v_1'\big)
+
n_2'\big(u_1'v_3'+u_3'v_1'\big)
\equiv
\\
&
\ \ \ \ \ \ \ \
\equiv
z'\overline{z}'
+
z'\overline{z}'\,{\rm O}_1\big(z',\overline{z}'\big)
+
d_2\,z'\overline{z}'\,u_1'
+
e_2\,z'\overline{z}'\,u_2'
+
f_2\,z'\overline{z}'\,u_3'
\,\,\,
\mod\,{\rm O}(4),
\endaligned
\]
\[
\aligned
&
v_3'
+
l_3'\big(u_1'v_1'+u_1'v_1'\big)
+
m_3'\big(u_1'v_2'+u_2'v_1'\big)
+
n_3'\big(u_1'v_3'+u_3'v_1'\big)
\equiv
\\
&
\ \ \ \ \ \ \ \
\equiv
z'\overline{z}'
+
z'\overline{z}'\,{\rm O}_1\big(z',\overline{z}'\big)
+
d_3\,z'\overline{z}'\,u_1'
+
e_3\,z'\overline{z}'\,u_2'
+
f_3\,z'\overline{z}'\,u_3'
\,\,\,
\mod\,{\rm O}(4).
\endaligned
\]

The equations: 
\[
\aligned
v_1'
&
=
\varphi_1'\big(z',\overline{z}',u_1',u_2',u_3'\big),
\\
v_2'
&
=
\varphi_2'\big(z',\overline{z}',u_1',u_2',u_3'\big),
\\
v_3'
&
=
\varphi_3'\big(z',\overline{z}',u_1',u_2',u_3'\big),
\endaligned
\]
of the transformed CR-generic submanifold:
\[
{M'}^5
\,\subset\,
{\C'}^4
\]
are obtained, up to order 4, by solving the above
three equations with respect to $v_1'$, $v_2'$, $v_3'$,
and one finds, modulo ${\rm O} ( 4)$:
\[
\aligned
v_1'
&
\equiv
z'\overline{z}'
+
z'\overline{z}'\,{\rm O}_1\big(z',\overline{z}'\big)
+
\underbrace{\big(d_1-2l_1'\big)}_{=:\,d_1'}\,z'\overline{z}'\,u_1'
+
\underbrace{\big(e_1-m_1'\big)}_{=:\,e_1'}\,z'\overline{z}'\,u_2'
+
\underbrace{\big(f_1-n_1'\big)}_{=:\,f_1'}\,z'\overline{z}'\,u_3',
\\
v_2'
&
\equiv
z'\overline{z}'
+
z'\overline{z}'\,{\rm O}_1\big(z',\overline{z}'\big)
+
\underbrace{\big(d_2-2l_2'\big)}_{=:\,d_2'}\,z'\overline{z}'\,u_1'
+
\underbrace{\big(e_2-m_2'\big)}_{=:\,e_2'}\,z'\overline{z}'\,u_2'
+
\underbrace{\big(f_2-n_2'\big)}_{=:\,f_2'}\,z'\overline{z}'\,u_3',
\\
v_3'
&
\equiv
z'\overline{z}'
+
z'\overline{z}'\,{\rm O}_1\big(z',\overline{z}'\big)
+
\underbrace{\big(d_3-2l_3'\big)}_{=:\,d_3'}\,z'\overline{z}'\,u_1'
+
\underbrace{\big(e_3-m_3'\big)}_{=:\,e_3'}\,z'\overline{z}'\,u_2'
+
\underbrace{\big(f_3-n_3'\big)}_{=:\,f_3'}\,z'\overline{z}'\,u_3'.
\endaligned
\]
Setting:
\[
\aligned
&
l_1':=\tfrac{1}{2}\,d_1,
\ \ \ \ \ \ \ \ \ \ \ \ \ \ \ \ \ \ \ \
m_1':=e_1,
\ \ \ \ \ \ \ \ \ \ \ \ \ \ \ \ \ \ \ \
n_1'
:=
f_1,
\\
&
l_2':=\tfrac{1}{2}\,d_2,
\ \ \ \ \ \ \ \ \ \ \ \ \ \ \ \ \ \ \ \
m_2':=e_2,
\ \ \ \ \ \ \ \ \ \ \ \ \ \ \ \ \ \ \ \
n_2'
:=
f_2,
\\
&
l_3':=\tfrac{1}{2}\,d_3,
\ \ \ \ \ \ \ \ \ \ \ \ \ \ \ \ \ \ \ \
m_3':=e_3,
\ \ \ \ \ \ \ \ \ \ \ \ \ \ \ \ \ \ \ \
n_3'
:=
f_3,
\endaligned
\]
one concludes.
\endproof

Emphasizing order $3$ terms in $(z, \overline{ z})$, 
plus order $4$ terms in the first line, 
one therefore arrives at:
\[
\!\!\!\!\!\!\!\!\!\!\!\!\!\!\!\!\!\!\!\!
\aligned
v_1
&
=
z\overline{z}
+
\alpha_1\,z^2\overline{z}
+
\overline{\alpha}_1\,z\overline{z}^2
+
\beta_1\,z^3\overline{z}
+
\overline{\beta}_1\,z\overline{z}^3
+
c_1\,z^2\overline{z}^2
+
z\overline{z}\,{\rm O}_3\big(z,\overline{z}\big)
+
z\overline{z}\,u_1\,{\rm O}_1\big(z,\overline{z},u_1\big)
+
\\
&
\ \ \ \ \ \ \ \ \ \ \ \ \ \ \ \ \ \ \ \ \ \ \ \ \ \ \ \ \ \ \ \ \ \ 
\ \ \ \ \ \ \ \ \ \ \ \ \ \ \ \ \ \ \ \ \ \ \ \ \ \ \ \ \ \ \
+
z\overline{z}\,u_2\,{\rm O}_1\big(z,\overline{z},u_1,u_2\big)
+
z\overline{z}\,u_3\,{\rm O}_1\big(z,\overline{z},u_1,u_2,u_3\big)
\\
v_2
&
=
\alpha_2\,z^2\overline{z}
+
\overline{\alpha}_2\,z\overline{z}^2
+
z\overline{z}\,{\rm O}_2\big(z,\overline{z}\big)
+
z\overline{z}\,u_1\,{\rm O}_1\big(z,\overline{z},u_1\big)
+
\\
&
\ \ \ \ \ \ \ \ \ \ \ \ \ \ \ \ \ \ \ \ \ \ \ \ \ \ \ \ \ \ \ 
+
z\overline{z}\,u_2\,{\rm O}_1\big(z,\overline{z},u_1,u_2\big)
+
z\overline{z}\,u_3\,{\rm O}_1\big(z,\overline{z},u_1,u_2,u_3\big)
\\
v_3
&
=
\alpha_3\,z^2\overline{z}
+
\overline{\alpha}_3\,z\overline{z}^2
+
z\overline{z}\,{\rm O}_2\big(z,\overline{z}\big)
+
z\overline{z}\,u_1\,{\rm O}_1\big(z,\overline{z},u_1\big)
+
\\
&
\ \ \ \ \ \ \ \ \ \ \ \ \ \ \ \ \ \ \ \ \ \ \ \ \ \ \ \ \ \ \ 
+
z\overline{z}\,u_2\,{\rm O}_1\big(z,\overline{z},u_1,u_2\big)
+
z\overline{z}\,u_3\,{\rm O}_1\big(z,\overline{z},u_1,u_2,u_3\big).
\endaligned
\]
In the first line, replacing:
\[
z
\,\longmapsto\,
z+\alpha_1\,z^2+\beta_1\,z^3,
\]
one makes:
\[
\alpha_1
=
0
=
\beta_1,
\]
without modifying the general form of the remainders. So: 
\[
\aligned
v_1
&
=
z\overline{z}
+
c_1\,z^2\overline{z}^2
+
z\overline{z}\,{\rm O}_3\big(z,\overline{z}\big)
+
z\overline{z}\,u_1\,{\rm O}_1\big(z,\overline{z},u_1\big)
+
\\
&
\ \ \ \ \ \ \ \ \ \ \ \ \ \ \ \ \ \ \ \ \ \ \ \ \ \ \ \ \ \ \ 
+
z\overline{z}\,u_2\,{\rm O}_1\big(z,\overline{z},u_1,u_2\big)
+
z\overline{z}\,u_3\,{\rm O}_1\big(z,\overline{z},u_1,u_2,u_3\big)
\\
v_2
&
=
\alpha_2\,z^2\overline{z}
+
\overline{\alpha}_2\,z\overline{z}^2
+
z\overline{z}\,{\rm O}_2\big(z,\overline{z}\big)
+
z\overline{z}\,u_1\,{\rm O}_1\big(z,\overline{z},u_1\big)
+
\\
&
\ \ \ \ \ \ \ \ \ \ \ \ \ \ \ \ \ \ \ \ \ \ \ \ \ \ \ \ \ \ \ 
+
z\overline{z}\,u_2\,{\rm O}_1\big(z,\overline{z},u_1,u_2\big)
+
z\overline{z}\,u_3\,{\rm O}_1\big(z,\overline{z},u_1,u_2,u_3\big)
\\
v_3
&
=
\alpha_3\,z^2\overline{z}
+
\overline{\alpha}_3\,z\overline{z}^2
+
z\overline{z}\,{\rm O}_2\big(z,\overline{z}\big)
+
z\overline{z}\,u_1\,{\rm O}_1\big(z,\overline{z},u_1\big)
+
\\
&
\ \ \ \ \ \ \ \ \ \ \ \ \ \ \ \ \ \ \ \ \ \ \ \ \ \ \ \ \ \ \ 
+
z\overline{z}\,u_2\,{\rm O}_1\big(z,\overline{z},u_1,u_2\big)
+
z\overline{z}\,u_3\,{\rm O}_1\big(z,\overline{z},u_1,u_2,u_3\big).
\endaligned
\]

It is now time to recall (\cite{ Merker-Pocchiola-Sabzevari-5-CR-II})
that the coefficients of the natural local generator for
$T^{1, 0}M$:
\[
\mathcal{L}
=
\frac{\partial}{\partial z}
+
A_1\,\frac{\partial}{\partial u_1}
+
A_2\,\frac{\partial}{\partial u_2}
+
A_3\,\frac{\partial}{\partial u_3}
\]
are given by the formulas:
\[
A_1
=
\frac{1}{\Delta}\,
\left\vert\!
\begin{array}{ccc}
-\varphi_{1,z} & \varphi_{1,u_2} & \varphi_{1,u_3}
\\
-\varphi_{2,z} & \isqrt+\varphi_{2,u_2} & \varphi_{2,u_3}
\\
-\varphi_{3,z} & \varphi_{3,u_2} & \isqrt+\varphi_{3,u_3}
\end{array}
\!\right\vert,
\]
\[
A_2
=
\frac{1}{\Delta}\,
\left\vert\!
\begin{array}{ccc}
\isqrt+\varphi_{1,u_1} & -\varphi_{1,z} & \varphi_{1,u_3}
\\
\varphi_{2,u_1} & -\varphi_{2,z} & \varphi_{2,u_3}
\\
\varphi_{3,u_1} & -\varphi_{3,z} & \isqrt+\varphi_{3,u_3}
\end{array}
\!\right\vert,
\]
\[
A_3
=
\frac{1}{\Delta}\,
\left\vert\!
\begin{array}{ccc}
\isqrt+\varphi_{1,u_1} & \varphi_{1,u_2} & -\varphi_{1,z}
\\
\varphi_{2,u_1} & \isqrt+\varphi_{2,u_2} & -\varphi_{2,z}
\\
\varphi_{3,u_1} & \varphi_{3,u_2} & -\varphi_{3,z}
\end{array}
\!\right\vert,
\]
with common denominator:
\[
\Delta
\,:=\,
\left\vert\!
\begin{array}{ccc}
\isqrt+\varphi_{1,u_1} & \varphi_{1,u_2} & \varphi_{1,u_3}
\\
\varphi_{2,u_1} & \isqrt+\varphi_{2,u_2} & \varphi_{2,u_3}
\\
\varphi_{3,u_1} & \varphi_{3,u_2} & \isqrt+\varphi_{3,u_3}
\end{array}
\!\right\vert.
\]

One notices:
\[
\varphi_{j,u_l}
=
z\overline{z}\,{\rm O}(1)
=
{\rm O}(3)
\eqno
{\scriptstyle{(1\,\leqslant\,j,\,l\,\leqslant\,3)}},
\]
whence:
\[
\aligned
\Delta
&
=
\left\vert\!
\begin{array}{ccc}
\isqrt+{\rm O}(3) & {\rm O}(3) & {\rm O}(3) 
\\
{\rm O}(3) & \isqrt+{\rm O}(3) & {\rm O}(3) 
\\
{\rm O}(3) & {\rm O}(3) & \isqrt+{\rm O}(3) 
\end{array}
\!\right\vert
\\
&
=
\isqrt^3
+
z\overline{z}\,{\rm O}(1)
=
\isqrt^3
+
{\rm O}(3).
\endaligned
\]

Concerning numerators, that of $A_1$ is:
\[
A_1^{\sf num}
=
\left\vert\!
\begin{array}{ccc}
-\overline{z}+{\rm O}(3) & {\rm O}(3) & {\rm O}(3)
\\
-2\alpha_2z\overline{z}-\overline{\alpha}_2\overline{z}^2+{\rm O}(3)
& \isqrt+{\rm O}(3) & {\rm O}(3)
\\
-2\alpha_3z\overline{z}-\overline{\alpha}_3\overline{z}^2+{\rm O}(3)
& {\rm O}(3) & \isqrt+{\rm O}(3)
\end{array}
\!\right\vert,
\]
so that:
\[
\aligned
A_1
&
=
\frac{-\,\isqrt^2\,\overline{z}+{\rm O}(3)}{
\isqrt^3+{\rm O}(3)}
\\
&
=
\isqrt\,\overline{z}
+
{\rm O}(3).
\endaligned
\]

Next:
\[
A_2^{\sf num}
=
\left\vert\!
\begin{array}{ccc}
\isqrt+{\rm O}(3) &
-\overline{z}+{\rm O}(3) & {\rm O}(3)
\\ 
{\rm O}(3) &
-2\alpha_2z\overline{z}-\overline{\alpha}_2\overline{z}^2+{\rm O}(3)
& {\rm O}(3)
\\
{\rm O}(3) &
-2\alpha_3z\overline{z}-\overline{\alpha}_3\overline{z}^2+{\rm O}(3)
 & \isqrt+{\rm O}(3)
\end{array}
\!\right\vert,
\]
so that:
\[
\aligned
A_1
&
=
\frac{-\,\isqrt^2\,2\alpha_2\,z\overline{z}
-\isqrt^2\,\overline{\alpha}_2\,\overline{z}^2+{\rm O}(3)}{
\isqrt^3+{\rm O}(3)}
\\
&
=
\isqrt\,2\alpha_2\,z\overline{z}
+
\isqrt\,\overline{\alpha}_2\,\overline{z}^2
+
{\rm O}(3).
\endaligned
\]

Quite similarly (mental exercise):
\[
A_3
=
\isqrt\,2\alpha_3\,z\overline{z}
+
\isqrt\,\overline{\alpha}_3\,\overline{z}^2
+
{\rm O}(3).
\]

Thus:
\[
\aligned
\mathcal{L}
&
=
\frac{\partial}{\partial z}
+
\isqrt\,\Big(
\overline{z}+{\rm O}(3)
\Big)\,\frac{\partial}{\partial u_1}
+
\isqrt\,\Big(
2\alpha_2\,z\overline{z}+\overline{\alpha}_2\,\overline{z}^2+{\rm O}(3)
\Big)\,\frac{\partial}{\partial u_2}
+
\\
&
\ \ \ \ \ \ \ \ \ \ \ \ \ \ \ \ \ \ \ \ \ \ \ \ \ \ \ \ \ \ \ \ \ \ 
\ \ \ \ \ \ \ \ \ \ \ \ \ \ \ 
+
\isqrt\,\Big(
2\alpha_3\,z\overline{z}+\overline{\alpha}_3\,\overline{z}^2+{\rm O}(3)
\Big)\,\frac{\partial}{\partial u_3},
\endaligned
\]
\[
\aligned
\overline{\mathcal{L}}
&
=
\frac{\partial}{\partial\overline{z}}
-
\isqrt\,\Big(
z+{\rm O}(3)
\Big)\,\frac{\partial}{\partial u_1}
-
\isqrt\,\Big(
2\overline{\alpha}_2\,z\overline{z}+\alpha_2\,z^2+{\rm O}(3)
\Big)\,\frac{\partial}{\partial u_2}
+
\\
&
\ \ \ \ \ \ \ \ \ \ \ \ \ \ \ \ \ \ \ \ \ \ \ \ \ \ \ \ \ \ \ \ \ \ 
\ \ \ \ \ \ \ \ \ \ \ \ \ \ \ 
+
\isqrt\,\Big(
2\overline{\alpha}_3\,z\overline{z}
+
\alpha_3\,z^2+{\rm O}(3)
\Big)\,\frac{\partial}{\partial u_3},
\endaligned
\]
whence:
\[
\aligned
\big[\mathcal{L},\overline{\mathcal{L}}\big]
=
-\,\isqrt\,
\Big(
2+{\rm O}(2)
\Big)\,
\frac{\partial}{\partial u_1}
&
-
\isqrt\,
\Big(
4\alpha_2\,z+4\overline{\alpha}_2\,\overline{z}
+
{\rm O}(2)
\Big)\,
\frac{\partial}{\partial u_2}
-
\\
&
-
\isqrt\,
\Big(
4\alpha_3\,z+4\overline{\alpha}_3\,\overline{z}
+
{\rm O}(2)
\Big)\,
\frac{\partial}{\partial u_3},
\endaligned
\]
and lastly:
\[
\aligned
\big[\mathcal{L},\,\big[\mathcal{L},\overline{\mathcal{L}}\big]\big]
&
=
{\rm O}(1)\,\frac{\partial}{\partial u_1}
-
\isqrt\,\Big(
4\alpha_2+{\rm O}(1)
\Big)\,
\frac{\partial}{\partial u_2}
-
\isqrt\,\Big(
4\alpha_3+{\rm O}(1)
\Big)\,
\frac{\partial}{\partial u_3},
\\
\big[\overline{\mathcal{L}},\,
\big[\mathcal{L},\overline{\mathcal{L}}\big]\big]
&
=
{\rm O}(1)\,\frac{\partial}{\partial u_1}
-
\isqrt\,\Big(
4\overline{\alpha}_2+{\rm O}(1)
\Big)\,
\frac{\partial}{\partial u_2}
-
\isqrt\,\Big(
4\overline{\alpha}_3+{\rm O}(1)
\Big)\,
\frac{\partial}{\partial u_3}.
\endaligned
\]

At the origin, visibly:
\[
\aligned
\mathcal{L}\big\vert_0
&
=
\frac{\partial}{\partial z}
\bigg\vert_0,
\\
\overline{\mathcal{L}}\big\vert_0
&
=
\frac{\partial}{\partial\overline{z}}
\bigg\vert_0,
\\
\big[\mathcal{L},\overline{\mathcal{L}}\big]
\Big\vert_0
&
=
-\,2\,\isqrt\,
\frac{\partial}{\partial u_1}
\bigg\vert_0,
\\
\big[\mathcal{L},\,
\big[\mathcal{L},\overline{\mathcal{L}}\big]\big]
\Big\vert_0
&
=
-\,4\,\isqrt\,\alpha_2\,
\frac{\partial}{\partial u_2}
\bigg\vert_0
-
4\,\isqrt\,\alpha_3\,
\frac{\partial}{\partial u_3}
\bigg\vert_0,
\\
\big[\overline{\mathcal{L}},\,
\big[\mathcal{L},\overline{\mathcal{L}}\big]\big]
\Big\vert_0
&
=
-\,4\,\isqrt\,\overline{\alpha}_2\,
\frac{\partial}{\partial u_2}
\bigg\vert_0
-
4\,\isqrt\,\overline{\alpha}_3\,
\frac{\partial}{\partial u_3}
\bigg\vert_0,
\endaligned
\]
By hypothesis, these five vectors should {\em not}
constitute a basis for:
\[
\C\otimes_\R T_0M
\,=\,
\C\frac{\partial}{\partial z}
\bigg\vert_0
\oplus
\C\frac{\partial}{\partial\overline{z}}
\bigg\vert_0
\oplus
\C\frac{\partial}{\partial u_1}
\bigg\vert_0
\oplus
\C\frac{\partial}{\partial u_2}
\bigg\vert_0
\oplus
\C\frac{\partial}{\partial u_3}
\bigg\vert_0,
\] 
of dimension {\bf 5}, while the first four should be independent,
hence:
\[
{\bf 4}
=
\rank_\C
\left(\!
\begin{array}{ccccc}
1 & 0 & 0 & 0 & 0
\\
0 & 1 & 0 & 0 & 0
\\
0 & 0 & -2\isqrt & 0 & 0
\\
0 & 0 & 0 & -4\isqrt\alpha_2 & -4\isqrt\alpha_3
\\
0 & 0 & 0 & -4\isqrt\overline{\alpha_2} & -4\isqrt\overline{\alpha}_3
\end{array}
\!\right).
\]

Without loss of generality:
\[
\alpha_2
\neq
0,
\]
and doing:
\[
z\,\longmapsto\,
\lambda\,z,
\]
with:
\[
1
=
\lambda^2\overline{\lambda}\,\alpha_2,
\]
one makes:
\[
\alpha_2
=
1.
\]

Thus:
\[
0
=
\left\vert\!
\begin{array}{cc}
1 & \alpha_3
\\
1 & \overline{\alpha}_3
\end{array}
\!\right\vert,
\]
that is to say:
\[
\alpha_3
=:
a_3
\,\in\,R.
\]

Lastly, doing:
\[
w_3
\,\longmapsto\,
w_3-a_3\,w_2,
\]
one makes:
\[
\alpha_3
=
0.
\]

\medskip

Emphasize now fourth order terms in all lines:
\[
\aligned
v_1
&
=
z\overline{z}
+
c_1\,z^2\overline{z}^2
+
z\overline{z}\,{\rm O}_3\big(z,\overline{z}\big)
+
z\overline{z}\,u_1\,{\rm O}_1\big(z,\overline{z},u_1\big)
+
\\
&
\ \ \ \ \ \ \ \ \ \ \ \ \ \ \ \ \ \ \ \ \ \ \ \ \ \ \ \ \ \ \ 
+
z\overline{z}\,u_2\,{\rm O}_1\big(z,\overline{z},u_1,u_2\big)
+
z\overline{z}\,u_3\,{\rm O}_1\big(z,\overline{z},u_1,u_2,u_3\big),
\\
v_2
&
=
z^2\overline{z}
+
z\overline{z}^2
+
\beta_2\,z^3\overline{z}
+
\overline{\beta}_2\,z\overline{z}^3
+
c_2\,z^2\overline{z}^2
+
z\overline{z}\,{\rm O}_3\big(z,\overline{z}\big)
+
z\overline{z}\,u_1\,{\rm O}_1\big(z,\overline{z},u_1\big)
+
\\
&
\ \ \ \ \ \ \ \ \ \ \ \ \ \ \ \ \ \ \ \ \ \ \ \ \ \ \ \ \ \ \ 
+
z\overline{z}\,u_2\,{\rm O}_1\big(z,\overline{z},u_1,u_2\big)
+
z\overline{z}\,u_3\,{\rm O}_1\big(z,\overline{z},u_1,u_2,u_3\big),
\\
v_3
&
=
\beta_3\,z^3\overline{z}
+
\overline{\beta}_3\,z\overline{z}^3
+
c_3\,z^2\overline{z}^2
+
z\overline{z}\,{\rm O}_3\big(z,\overline{z}\big)
+
z\overline{z}\,u_1\,{\rm O}_1\big(z,\overline{z},u_1\big)
+
\\
&
\ \ \ \ \ \ \ \ \ \ \ \ \ \ \ \ \ \ \ \ \ \ \ \ \ \ \ \ \ \ \ 
+
z\overline{z}\,u_2\,{\rm O}_1\big(z,\overline{z},u_1,u_2\big)
+
z\overline{z}\,u_3\,{\rm O}_1\big(z,\overline{z},u_1,u_2,u_3\big).
\endaligned
\]
Again:
\[
\Delta
=
\isqrt^3
+
z\overline{z}\,{\rm O}(0).
\]

Next:
\[
\!\!\!\!\!\!\!\!\!\!\!\!\!\!\!\!\!\!\!\!
\footnotesize
\aligned
A_1^{\sf num}
&
=
\left\vert\!
\def\arraystretch{1.25}
\begin{array}{ccc}
\substack{
-\overline{z}-2c_1z\overline{z}^2
+\overline{z}\,{\rm O}_3(z,\overline{z})
+\overline{z}u_1{\rm O}_1(z,\overline{z},u_1)+
\\
+\overline{z}u_2{\rm O}_1(z,\overline{z},u_1,u_2)
+\overline{z}u_3{\rm O}_1(z,\overline{z},u_1,u_2,u_3)}
\bigskip
&
z\overline{z}{\rm O}_1(z,\overline{z},u_1,u_2,u_3)
&
z\overline{z}{\rm O}_1(z,\overline{z},u_1,u_2,u_3)
\\
\substack{
-2z\overline{z}-\overline{z}^2-3\beta_2z^2\overline{z}
-\overline{\beta}_2\overline{z}^3-2c_2z\overline{z}^2+
\\
+\overline{z}\,{\rm O}_3(z,\overline{z})
+\overline{z}u_1{\rm O}_1(z,\overline{z},u_1)+
\\
+\overline{z}u_2{\rm O}_1(z,\overline{z},u_1,u_2)
+\overline{z}u_3{\rm O}_1(z,\overline{z},u_1,u_2,u_3)}
\bigskip
&
\isqrt+z\overline{z}{\rm O}_1(z,\overline{z},u_1,u_2,u_3)
&
z\overline{z}{\rm O}_1(z,\overline{z},u_1,u_2,u_3)
\\
\substack{
-3\beta_3z^2\overline{z}
-\overline{\beta}_3\overline{z}^3-2c_3z\overline{z}^2+
\\
+\overline{z}\,{\rm O}_3(z,\overline{z})
+\overline{z}u_1{\rm O}_1(z,\overline{z},u_1)+
\\
+\overline{z}u_2{\rm O}_1(z,\overline{z},u_1,u_2)
+\overline{z}u_3{\rm O}_1(z,\overline{z},u_1,u_2,u_3)}
&
z\overline{z}{\rm O}_1(z,\overline{z},u_1,u_2,u_3)
&
\isqrt+z\overline{z}{\rm O}_1(z,\overline{z},u_1,u_2,u_3)
\end{array}
\!\right\vert
\\
&
=
-\,\isqrt^2\,
\Big(
\overline{z}+2c_1z\overline{z}^2
+
\overline{z}\,{\rm O}_3(z,\overline{z})
+
\overline{z}\,u_1\,{\rm O}_1\big(z,\overline{z},u_1\big)
+
\overline{z}\,u_2\,{\rm O}_1\big(z,\overline{z},u_1,u_2\big)
+
\overline{z}\,u_3\,{\rm O}_1\big(z,\overline{z},u_1,u_2,u_3\big)
\Big),
\endaligned
\]
so that:
\[
\!\!\!\!\!\!\!\!\!\!\!\!\!\!\!\!\!\!\!\!
\footnotesize
\aligned
A_1
&
=
\frac{
-\,\isqrt^2\,
\Big(
\overline{z}+2c_1z\overline{z}^2
+
\overline{z}\,{\rm O}_3(z,\overline{z})
+
\overline{z}\,u_1\,{\rm O}_1\big(z,\overline{z},u_1\big)
+
\overline{z}\,u_2\,{\rm O}_1\big(z,\overline{z},u_1,u_2\big)
+
\overline{z}\,u_3\,{\rm O}_1\big(z,\overline{z},u_1,u_2,u_3\big)
\Big)}{
\isqrt^3+z\overline{z}\,{\rm O}(1)}
\\
&
=
\isqrt\,
\Big(
\overline{z}+2c_1z\overline{z}^2
+
\overline{z}\,{\rm O}_3(z,\overline{z})
+
\overline{z}\,u_1\,{\rm O}_1\big(z,\overline{z},u_1\big)
+
\overline{z}\,u_2\,{\rm O}_1\big(z,\overline{z},u_1,u_2\big)
+
\overline{z}\,u_3\,{\rm O}_1\big(z,\overline{z},u_1,u_2,u_3\big)
\Big).
\endaligned
\]

Secondly:
\[
\footnotesize
\aligned
A_2^{\sf num}
&
=
\left\vert\!
\def\arraystretch{1.25}
\begin{array}{ccc}
\isqrt+z\overline{z}{\rm O}_1(z,\overline{z},u_1,u_2,u_3)
&
\substack{
-\overline{z}-2c_1z\overline{z}^2
+\overline{z}\,{\rm O}_3(z,\overline{z})
+\overline{z}u_1{\rm O}_1(z,\overline{z},u_1)+
\\
+\overline{z}u_2{\rm O}_1(z,\overline{z},u_1,u_2)
+\overline{z}u_3{\rm O}_1(z,\overline{z},u_1,u_2,u_3)}
\bigskip
&
z\overline{z}{\rm O}_1(z,\overline{z},u_1,u_2,u_3)
\\
z\overline{z}{\rm O}_1(z,\overline{z},u_1,u_2,u_3)
&
\substack{
-2z\overline{z}-\overline{z}^2-3\beta_2z^2\overline{z}
-\overline{\beta}_2\overline{z}^3-2c_2z\overline{z}^2+
\\
+\overline{z}\,{\rm O}_3(z,\overline{z})
+\overline{z}u_1{\rm O}_1(z,\overline{z},u_1)+
\\
+\overline{z}u_2{\rm O}_1(z,\overline{z},u_1,u_2)
+\overline{z}u_3{\rm O}_1(z,\overline{z},u_1,u_2,u_3)}
\bigskip
&
z\overline{z}{\rm O}_1(z,\overline{z},u_1,u_2,u_3)
\\
z\overline{z}{\rm O}_1(z,\overline{z},u_1,u_2,u_3)
&
\substack{
-3\beta_3z^2\overline{z}
-\overline{\beta}_3\overline{z}^3-2c_3z\overline{z}^2+
\\
+\overline{z}\,{\rm O}_3(z,\overline{z})
+\overline{z}u_1{\rm O}_1(z,\overline{z},u_1)+
\\
+\overline{z}u_2{\rm O}_1(z,\overline{z},u_1,u_2)
+\overline{z}u_3{\rm O}_1(z,\overline{z},u_1,u_2,u_3)}
&
\isqrt+z\overline{z}{\rm O}_1(z,\overline{z},u_1,u_2,u_3)
\end{array}
\!\right\vert
\\
&
=
-\,\isqrt^2\,
\Big(
2\,z\overline{z}
+
\overline{z}^2
+
3\beta_2\,z^2\overline{z}
+
\overline{\beta}_2\,\overline{z}^3
+
2c_2\,z\overline{z}^2
+
\\
&
\ \ \ \ \ \ \ \ \ \ \ \ \ \ \ \ \ \ \ \ \ \ \ \ 
+
\overline{z}\,{\rm O}_3(z,\overline{z})
+
\overline{z}\,u_1\,{\rm O}_1\big(z,\overline{z},u_1\big)
+
\overline{z}\,u_2\,{\rm O}_1\big(z,\overline{z},u_1,u_2\big)
+
\overline{z}\,u_3\,{\rm O}_1\big(z,\overline{z},u_1,u_2,u_3\big)
\Big),
\endaligned
\]
so that:
\[
\footnotesize
\aligned
A_2
&
=
\frac{A_2^{\sf num}}{\isqrt^3+z\overline{z}\,{\rm O}(1)}
\\
&
=
\isqrt\,
\Big(
2\,z\overline{z}
+
\overline{z}^2
+
3\beta_2\,z^2\overline{z}
+
\overline{\beta}_2\,\overline{z}^3
+
2c_2\,z\overline{z}^2
+
\\
&
\ \ \ \ \ \ \ \ \ \ \ \ \ \ \ \ \ \ \ \ \ \ \ \ 
+
\overline{z}\,{\rm O}_3(z,\overline{z})
+
\overline{z}\,u_1\,{\rm O}_1\big(z,\overline{z},u_1\big)
+
\overline{z}\,u_2\,{\rm O}_1\big(z,\overline{z},u_1,u_2\big)
+
\overline{z}\,u_3\,{\rm O}_1\big(z,\overline{z},u_1,u_2,u_3\big)
\Big).
\endaligned
\]

Quite similarly:
\[
\footnotesize
\aligned
A_3
&
=
\isqrt\,
\Big(
3\beta_3\,z^2\overline{z}
+
\overline{\beta}_3\,\overline{z}^3
+
2c_3\,z\overline{z}^2
+
\\
&
\ \ \ \ \ \ \ \ \ \ \ \ \ \ \ \ \ \ \ \ \ \ \ \ 
+
\overline{z}\,{\rm O}_3(z,\overline{z})
+
\overline{z}\,u_1\,{\rm O}_1\big(z,\overline{z},u_1\big)
+
z\overline{z}\,u_2\,{\rm O}_0\big(z,\overline{z},u_1,u_2\big)
+
z\overline{z}\,u_3\,{\rm O}_0\big(z,\overline{z},u_1,u_2,u_3\big)
\Big).
\endaligned
\]

To compactify remainders, attribute weights:
\[
\aligned
\weight(z)
&
:=
1,
\\
\weight(w_1)
&
:=
2,
\\
\weight(w_2)
&
:=
3,
\\
\weight(w_3)
&
:\geqslant
4.
\endaligned
\]

Hence:
\[
\aligned
\mathcal{L}
&
=
\frac{\partial}{\partial z}
+
\isqrt\,
\Big(
\overline{z}
+
2c_1\,z\overline{z}^2
+
{\rm O}_{\sf weighted}(4)
\Big)\,
\frac{\partial}{\partial u_1}
+
\\
&
\ \ \ \ \ \ \ \ \ \ \ 
+
\isqrt\Big(
2\,z\overline{z}+\overline{z}^2
+
3\beta_2\,z^2\overline{z}
+
\overline{\beta}_2\,\overline{z}^3
+
2c_2\,z\overline{z}^2
+
{\rm O}_{\sf weighted}(4)
\Big)\,
\frac{\partial}{\partial u_2}
+
\\
&
\ \ \ \ \ \ \ \ \ \ \ 
+
\isqrt\Big(
3\beta_3\,z^2\overline{z}
+
\overline{\beta}_3\,\overline{z}^3
+
2c_3\,z\overline{z}^2
+
{\rm O}_{\sf weighted}(4)
\Big)\,
\frac{\partial}{\partial u_3},
\endaligned
\]
\[
\aligned
\overline{\mathcal{L}}
&
=
\frac{\partial}{\partial\overline{z}}
-
\isqrt\,
\Big(
z
+
2c_1\,z^2\overline{z}
+
{\rm O}_{\sf weighted}(4)
\Big)\,
\frac{\partial}{\partial u_1}
-
\\
&
\ \ \ \ \ \ \ \ \ \ \ 
-
\isqrt\Big(
2\,z\overline{z}+z^2
+
3\overline{\beta}_2\,z\overline{z}^2
+
\beta_2\,z^3
+
2c_2\,z^2\overline{z}
+
{\rm O}_{\sf weighted}(4)
\Big)\,
\frac{\partial}{\partial u_2}
-
\\
&
\ \ \ \ \ \ \ \ \ \ \ 
-
\isqrt\Big(
3\overline{\beta}_3\,z\overline{z}^2
+
\beta_3\,z^3
+
2c_3\,z^2\overline{z}
+
{\rm O}_{\sf weighted}(4)
\Big)\,
\frac{\partial}{\partial u_3}.
\endaligned
\]

When computing brackets up to length three, remainders of weighted
order $\geqslant 4$ do not contribute (exercise):
\[
\aligned
\big[\mathcal{L},\overline{\mathcal{L}}\big]
&
=
-\,2\isqrt\,
\Big(
1+4c_1\,z\overline{z}
+
{\rm O}_{\sf weighted}(3)
\Big)\,
\frac{\partial}{\partial u_1}
+
\\
&
\ \ \ \ \ 
-\,2\,\isqrt\,
\Big(
2\,z+2\,\overline{z}
+
3\beta_2\,z^2
+
3\overline{\beta}_2\,\overline{z}^2
+
4c_2\,z\overline{z}
+
{\rm O}_{\sf weighted}(3)
\Big)\,
\frac{\partial}{\partial u_2}
+
\\
&
\ \ \ \ \ 
-\,2\,\isqrt\,
\Big(
3\beta_3\,z^2
+
3\overline{\beta}_3\,\overline{z}^2
+
4c_3\,z\overline{z}
+
{\rm O}_{\sf weighted}(3)
\Big)\,
\frac{\partial}{\partial u_3},
\endaligned
\]
\[
\aligned
\big[\mathcal{L},\,
\big[\mathcal{L},\overline{\mathcal{L}}\big]\big]
&
=
-\,2\isqrt\,
\Big(
4c_1\,\overline{z}
+
{\rm O}_{\sf weighted}(2)
\Big)\,
\frac{\partial}{\partial u_1}
+
\\
&
\ \ \ \ \ 
-\,2\,\isqrt\,
\Big(
2
+
6\beta_2\,z
+
4c_2\,\overline{z}
+
{\rm O}_{\sf weighted}(2)
\Big)\,
\frac{\partial}{\partial u_2}
+
\\
&
\ \ \ \ \ 
-\,2\,\isqrt\,
\Big(
6\beta_3\,z
+
4c_3\,\overline{z}
+
{\rm O}_{\sf weighted}(2)
\Big)\,
\frac{\partial}{\partial u_3},
\endaligned
\]
\[
\aligned
\big[\overline{\mathcal{L}},\,
\big[\mathcal{L},\overline{\mathcal{L}}\big]\big]
&
=
-\,2\isqrt\,
\Big(
4c_1z
+
{\rm O}_{\sf weighted}(2)
\Big)\,
\frac{\partial}{\partial u_1}
+
\\
&
\ \ \ \ \ 
-\,2\,\isqrt\,
\Big(
2
+
6\overline{\beta}_2\,\overline{z}
+
4c_2\,z
+
{\rm O}_{\sf weighted}(2)
\Big)\,
\frac{\partial}{\partial u_2}
+
\\
&
\ \ \ \ \ 
-\,2\,\isqrt\,
\Big(
6\overline{\beta}_3\,\overline{z}
+
4c_3\,z
+
{\rm O}_{\sf weighted}(2)
\Big)\,
\frac{\partial}{\partial u_3}.
\endaligned
\]

The hypothesis:
\[
{\bf 4}
\,=\,
\rank_\C\Big(
\Big\{
\mathcal{L},\,\overline{\mathcal{L}},\,\,
\big[\mathcal{L},\overline{\mathcal{L}}\big],\,\,
\big[\mathcal{L},\,\big[\mathcal{L},\overline{\mathcal{L}}\big]\big],\,\,
\big[\overline{\mathcal{L}},\,\big[\mathcal{L},
\overline{\mathcal{L}}\big]\big]
\Big\}
\Big)
\]
reads:
\[
0
\equiv
\left\vert\!
\begin{array}{ccc}
\substack{1+4c_1z\overline{z}+{\rm O}_{\sf weighted}(3)}
&
\substack{4c_1\overline{z}+{\rm O}_{\sf weighted}(2)}
&
\substack{4c_1z+{\rm O}_{\sf weighted}(2)}\bigskip
\\
\substack{
2z+2\overline{z}+3\beta_2z^2+\overline{\beta}_2\overline{z}^2+\\
+4c_2z\overline{z}+{\rm O}_{\sf weighted}(3)}
&
\substack{2+6\beta_2z+4c_2\overline{z}+\\
+{\rm O}_{\sf weighted}(2)}
&
\substack{2+6\overline{\beta}_2\overline{z}+4c_2z+\\
+{\rm O}_{\sf weighted}(2)}\bigskip
\\
\substack{3\beta_3z^2+3\overline{\beta}_3\overline{z}^2
+4c_3z\overline{z}+\\
+{\rm O}_{\sf weighted}(3)}
&
\substack{6\beta_3z+4c_3\overline{z}+\\
+{\rm O}_{\sf weighted}(2)}
&
\substack{6\overline{\beta}_3\overline{z}
+4c_3z+\\
+{\rm O}_{\sf weighted}(2)}
\end{array}
\!\right\vert,
\]
and picking only order $1$ terms:
\[
0
\equiv
12\,\overline{\beta}_3\,\overline{z}
-
12\,\beta_3\,z
+
8\,c_3\,z
-
8\,c_3\,\overline{z},
\]
which is:
\[
3\,\beta_3
=
2\,c_3,
\]
so that:
\[
\boxed{\,\text{\rm $\beta_3$ is real}.\,}
\]

Further, the hypothesis:
\[
{\bf 5}
\,=\,
\rank_\C\Big(
\Big\{
\mathcal{L},\,\overline{\mathcal{L}},\,\,
\big[\mathcal{L},\overline{\mathcal{L}}\big],\,\,
\big[\mathcal{L},\,\big[\mathcal{L},\overline{\mathcal{L}}\big]\big],\,\,
\big[\mathcal{L},\,\big[\mathcal{L},\,\big[\mathcal{L},
\overline{\mathcal{L}}\big]\big]\big]
\Big\}
\Big),
\]
requires (exercise):
\[
c_3
\neq
0.
\]

Naturally, a real dilation of $z$ makes:
\[
\beta_3
:=
2,
\ \ \ \ \ \ \ \ \ \ \ \ \ \ \ \ \ \ \ \ \ \ \ \
c_3
:=
3,
\]
so that order $4$ terms on the third line are
the ones announced:
\[
2\,z^3\overline{z}
+
2\,z\overline{z}^3
+
3\,z^2\overline{z}^2,
\]
while a simultaneous correcting dilation along the $w_1$-axis and
along the $w_2$-axis keeps unchanged the preceding normalizations,
which concludes.
\endproof

\noindent{\bf Scholium.}
{\em In such elementarily normalized coordinates, one has
the diagonal normalization at the origin:}
\[
\aligned
\mathcal{L}\big\vert_0
&
=
\frac{\partial}{\partial z}
\bigg\vert_0,
\\
\overline{\mathcal{L}}\big\vert_0
&
=
\frac{\partial}{\partial\overline{z}}
\bigg\vert_0,
\\
\big[\mathcal{L},\overline{\mathcal{L}}\big]
\Big\vert_0
&
=
-\,2\,\isqrt\,
\frac{\partial}{\partial u_1}
\bigg\vert_0,
\\
\big[\mathcal{L},\,\big[\mathcal{L},\overline{\mathcal{L}}\big]\big]
\Big\vert_0
&
=
-\,4\,\isqrt\,
\frac{\partial}{\partial u_2}
\bigg\vert_0,
\\
\big[\mathcal{L},\,
\big[\mathcal{L},\,\big[\mathcal{L},\overline{\mathcal{L}}\big]\big]\big]
\Big\vert_0
&
=
-\,24\,\isqrt\,
\frac{\partial}{\partial u_3}
\bigg\vert_0,
\endaligned
\]
{\em which conveniently fixes ideas 
when performing explicitly the Cartan equivalence
procedure.\qed}

\medskip

Moreover, by further analysis, one can
make $c_1 = 0$, $\beta_2 = 0$, $c_2 = 0$
after some appropriate biholomorphism,
keeping untouched all the preceding normalizations.
In any case, neither $c_1$, nor $\beta_2$, nor
$c_2$ entered the preceding invariant calculations
concerning iterated Lie brackets
between local vector field sections of $T^{1,0}M$ and 
of $T^{0,1}M$, which confirms that
$c_1$, $\beta_2$, $c_2$ have no invariant character. 

\medskip

Taking zero remainders, one obtains:
\[
\boxed{\,\,
{\text{\footnotesize\sf Model
$\text{\sf (III)}_{\text\sf 2}$:}}
\ \ \ \ \ \ \ \ \ 
\aligned
v_1
&
=
z\overline{z},\,\,
\\
v_2
&
=
z^2\overline{z}+z\overline{z}^2,
\\
v_3
&
=
2\,z^3\overline{z}+2\,z\overline{z}^3
+
3\,z^2\overline{z}^2.
\endaligned
\rule[-4pt]{0pt}{36pt}}
\]

\medskip

To conclude, one must express the constraint:
\[
{\bf 4}
\,=\,
\rank_\C\Big(
\Big\{
\mathcal{L},\,\overline{\mathcal{L}},\,\,
\big[\mathcal{L},\overline{\mathcal{L}}\big],\,\,
\big[\mathcal{L},\,\big[\mathcal{L},\overline{\mathcal{L}}\big]\big],\,\,
\big[\overline{\mathcal{L}},\,\big[\mathcal{L},
\overline{\mathcal{L}}\big]\big]
\Big\}
\Big).
\]

Starting from:
\[
\aligned
\mathcal{L}
&
=
\frac{\partial}{\partial z}
+
A_1\,\frac{\partial}{\partial u_1}
+
A_2\,\frac{\partial}{\partial u_2}
+
A_3\,\frac{\partial}{\partial u_3},
\\
\overline{\mathcal{L}}
&
=
\frac{\partial}{\partial\overline{z}}
+
\overline{A}_1\,\frac{\partial}{\partial u_1}
+
\overline{A}_2\,\frac{\partial}{\partial u_2}
+
\overline{A}_3\,\frac{\partial}{\partial u_3},
\endaligned
\]
firstly, secondly, one has thirdly:
\[
\small
\aligned
\big[\mathcal{L},\overline{\mathcal{L}}\big]
&
=
\Big(
\mathcal{L}\big(\overline{A}_1\big)
-
\overline{\mathcal{L}}\big(A_1\big)
\Big)\,
\frac{\partial}{\partial u_1}
+
\Big(
\mathcal{L}\big(\overline{A}_2\big)
-
\overline{\mathcal{L}}\big(A_2\big)
\Big)\,
\frac{\partial}{\partial u_2}
+
\Big(
\mathcal{L}\big(\overline{A}_3\big)
-
\overline{\mathcal{L}}\big(A_3\big)
\Big)\,
\frac{\partial}{\partial u_3},
\endaligned
\]
fourthly:
\[
\aligned
\big[\mathcal{L},\,
\big[\mathcal{L},\overline{\mathcal{L}}\big]\big]
&
=
\Big(
\mathcal{L}\big(\mathcal{L}\big(\overline{A}_1\big)\big)
-
2\,\mathcal{L}\big(\overline{\mathcal{L}}\big(A_1\big)\big)
+
\overline{\mathcal{L}}\big(\mathcal{L}\big(A_1\big)\big)
\Big)\,
\frac{\partial}{\partial u_1}
+
\\
&
+
\Big(
\mathcal{L}\big(\mathcal{L}\big(\overline{A}_2\big)\big)
-
2\,\mathcal{L}\big(\overline{\mathcal{L}}\big(A_2\big)\big)
+
\overline{\mathcal{L}}\big(\mathcal{L}\big(A_2\big)\big)
\Big)\,
\frac{\partial}{\partial u_2}
+
\\
&
+
\Big(
\mathcal{L}\big(\mathcal{L}\big(\overline{A}_3\big)\big)
-
2\,\mathcal{L}\big(\overline{\mathcal{L}}\big(A_3\big)\big)
+
\overline{\mathcal{L}}\big(\mathcal{L}\big(A_3\big)\big)
\Big)\,
\frac{\partial}{\partial u_3},
\endaligned
\]
and fifthly:
\[
\aligned
\big[\overline{\mathcal{L}},\,
\big[\mathcal{L},\overline{\mathcal{L}}\big]\big]
&
=
\Big(
-\,
\overline{\mathcal{L}}\big(\overline{\mathcal{L}}\big(A_1\big)\big)
+
2\,\overline{\mathcal{L}}\big(\mathcal{L}\big(\overline{A}_1\big)\big)
-
\mathcal{L}\big(\overline{\mathcal{L}}\big(\overline{A}_1\big)\big)
\Big)\,
\frac{\partial}{\partial u_1}
+
\\
&
+
\Big(
-\,
\overline{\mathcal{L}}\big(\overline{\mathcal{L}}\big(A_2\big)\big)
+
2\,\overline{\mathcal{L}}\big(\mathcal{L}\big(\overline{A}_2\big)\big)
-
\mathcal{L}\big(\overline{\mathcal{L}}\big(\overline{A}_2\big)\big)
\Big)\,
\frac{\partial}{\partial u_2}
+
\\
&
+
\Big(
-\,
\overline{\mathcal{L}}\big(\overline{\mathcal{L}}\big(A_3\big)\big)
+
2\,\overline{\mathcal{L}}\big(\mathcal{L}\big(\overline{A}_3\big)\big)
-
\mathcal{L}\big(\overline{\mathcal{L}}\big(\overline{A}_3\big)\big)
\Big)\,
\frac{\partial}{\partial u_3},
\endaligned
\]

So the hypothesis means the identical vanishing of the $3 \times 3$
determinant:
\[
\aligned
0
\equiv
\left\vert\!
\begin{array}{ccc}
\substack{
\mathcal{L}(\overline{A}_1)-\overline{\mathcal{L}}(A_1)}
&
\substack{
\mathcal{L}(\overline{A}_2)-\overline{\mathcal{L}}(A_2)}
&
\substack{
\mathcal{L}(\overline{A}_3)-\overline{\mathcal{L}}(A_3)}\bigskip
\\
\substack{
\mathcal{L}(\mathcal{L}(\overline{A}_1))
-2\mathcal{L}(\overline{\mathcal{L}}(A_1))+\\
+\overline{\mathcal{L}}(\mathcal{L}(A_1))}
&
\substack{
\mathcal{L}(\mathcal{L}(\overline{A}_2))
-2\mathcal{L}(\overline{\mathcal{L}}(A_2))+\\
+\overline{\mathcal{L}}(\mathcal{L}(A_2))}
&
\substack{
\mathcal{L}(\mathcal{L}(\overline{A}_3))
-2\mathcal{L}(\overline{\mathcal{L}}(A_3))+\\
+\overline{\mathcal{L}}(\mathcal{L}(A_3))}\bigskip
\\
\substack{
-\overline{\mathcal{L}}(\overline{\mathcal{L}}(A_1))
+2\overline{\mathcal{L}}(\mathcal{L}(\overline{A}_1))-\\
-\mathcal{L}(\overline{\mathcal{L}}(\overline{A}_1))}
&
\substack{
-\overline{\mathcal{L}}(\overline{\mathcal{L}}(A_2))
+2\overline{\mathcal{L}}(\mathcal{L}(\overline{A}_2))-\\
-\mathcal{L}(\overline{\mathcal{L}}(\overline{A}_2))}
&
\substack{
-\overline{\mathcal{L}}(\overline{\mathcal{L}}(A_3))
+2\overline{\mathcal{L}}(\mathcal{L}(\overline{A}_3))-\\
-\mathcal{L}(\overline{\mathcal{L}}(\overline{A}_3))}
\end{array}
\!\right\vert.
\endaligned
\]

\medskip

When expressed back in terms of the three graphing functions:
\[
\aligned
&
\varphi_1\big(z,\overline{z},u_1,u_2,u_3\big),
\\
&
\varphi_2\big(z,\overline{z},u_1,u_2,u_3\big),
\\
&
\varphi_3\big(z,\overline{z},u_1,u_2,u_3\big),
\endaligned
\]
one obtains a rational expression in the third-order
jet of the three graphing functions
$\varphi_1$, $\varphi_2$, $\varphi_3$ whose numerator 
contains hundreds of lines.


\bigskip

\section{\sf General class $\text{\sf IV}_{\text{\sf 1}}$}
\label{general-class-IV-1}
\HEAD{\ref{general-class-IV-1}.~General class 
$\text{\sf IV}_{\text{\sf 1}}$}{
Jo\"el {\sc Merker} (Paris-Sud)}

\medskip

\medskip\noindent{\bf Proposition.}
{\em A local real analytic hypersurface passing through the origin:}
\[
0
\,\in\,
M^5
\,\subset\,
\C^3
\]
{\em which belongs to the general Class $\text{\sf IV}_{
\text{\sf 1}}$, namely such that,
for any two local generators $\big\{
\mathcal{ L}_1, \mathcal{ L}_2 \big\}$ for $T^{1, 0}M$:}
\[
\Big\{
\mathcal{L}_1,\,\mathcal{L}_2,\,
\overline{\mathcal{L}}_1,\,\overline{\mathcal{L}}_2,\,\,
\big[\mathcal{L}_1,\overline{\mathcal{L}}_1\big]\Big\}
\]
{\em constitute a frame for $\C\otimes_\R TM$, 
and such that in addition, the Levi form:}
\[
\text{\footnotesize\sf Levi-Form}^M(p)
\]
{\em if of rank $2$ at every point, 
may always be represented, in suitable local holomorphic coordinates:}
\[
(z_1,z_2,w)
\in
\C^3
\]
{\em by a specific real analytic equation of the form:}
\[
\boxed{\,\,
\underline{\text{\footnotesize\sf 
$\text{\sf (IV)}_{\text{\sf 1}}$:}}
\ \ \ \ \
\aligned
v
=
z_1\overline{z}_1
\pm
z_2\overline{z}_2
+
{\rm O}_3\big(z_1,z_2,\overline{z}_1,\overline{z}_2,u\big),\,\,
\rule[-5pt]{0pt}{20pt}
\endaligned}
\]
{\em with remainder satisfying:}
\[
0
\equiv
\remainder
\big(0,0,\overline{z}_1,\overline{z}_2,u\big)
\equiv
\remainder
\big(z_1,z_2,0,0,u\big).
\]

\proof
One starts with:
\[
\aligned
v
&
=
\varphi\big(z_1,z_2,\overline{z}_1,\overline{z}_2,u)
\\
&
=
{\rm O}_2\big(z_1,z_2,\overline{z}_1,\overline{z}_2,u\big),
\endaligned
\]
assuming:
\[
0
\equiv
\varphi\big(0,0,\overline{z}_1,\overline{z}_2,u\big)
\equiv
\varphi\big(z_1,z_2,0,0,u\big).
\]

Hence:
\[
\aligned
v
&
=
P_2\big(z_1,z_2,\overline{z}_1,\overline{z}_2\big)
+
{\rm O}_3\big(z_1,z_2,\overline{z}_1,\overline{z}_2,u\big)
\\
&
=
a\,z_1\overline{z}_1
+
\beta\,z_2\overline{z}_1
+
\overline{\beta}\,z_1\overline{z}_2
+
c\,z_2\overline{z}_2
+
{\rm O}_3\big(z_1,z_2,\overline{z}_1,\overline{z}_2,u\big),
\endaligned
\]
with:
\[
a\in\R,
\ \ \ \ \ \ \ \ \ \ \ \ \ \ \ \ \ \ \ \ \ \
\beta\in\C,
\ \ \ \ \ \ \ \ \ \ \ \ \ \ \ \ \ \ \ \ \ \
c\in\R.
\]

The first natural local generator of $T^{1, 0}M$ is:
\[
\aligned
\mathcal{L}_1
&
=
\frac{\partial}{\partial z_1}
-
\frac{\varphi_{z_1}}{\isqrt+\varphi_u}\,
\frac{\partial}{\partial u}
\\
&
=
\frac{\partial}{\partial z_1}
-
\bigg(
\frac{a\,\overline{z}_1+\overline{\beta}\,\overline{z}_2+{\rm O}(2)}{
\isqrt+{\rm O}(2)}
\bigg)\,
\frac{\partial}{\partial u},
\\
&
=
\frac{\partial}{\partial z_1}
+
\isqrt\,\Big(
a\,\overline{z}_1+\overline{\beta}\,\overline{z}_2+{\rm O}(2)
\Big)\,
\frac{\partial}{\partial u},
\endaligned
\]
and similarly, the second is:
\[
\mathcal{L}_2
=
\frac{\partial}{\partial z_2}
+
\isqrt\,
\Big(
\beta\,\overline{z}_1+c\,\overline{z}_2+{\rm O}(2)
\Big)\,
\frac{\partial}{\partial u}.
\]

Conjugates are:
\[
\aligned
\overline{\mathcal{L}}_1
&
=
\frac{\partial}{\partial\overline{z}_1}
-
\isqrt\,
\Big(
a\,z_1+\beta\,z_2+{\rm O}(2)
\Big)\,
\frac{\partial}{\partial u},
\\
\overline{\mathcal{L}}_2
&
=
\frac{\partial}{\partial\overline{z}_2}
-
\isqrt\,
\Big(
\overline{\beta}\,z_1+c\,z_2+{\rm O}(2)
\Big)\,
\frac{\partial}{\partial u}.
\endaligned
\]

Taking:
\[
\rho_0
=
du
-
A_1\,dz_1
-
A_2\,dz_2
-
\overline{A}_1\,d\overline{z}_1
-
\overline{A}_2\,d\overline{z}_2,
\]
the Levi matrix at the origin is:
\[
\!\!\!\!\!\!\!\!\!\!\!\!\!\!\!\!\!\!\!\!
\!\!\!\!\!\!\!\!\!\!\!\!\!\!\!\!\!\!\!\!
\aligned
\left(\!
\begin{array}{cc}
\rho_0\big(\isqrt\big[\mathcal{L}_1,\overline{\mathcal{L}}_1\big]\big)
&
\rho_0\big(\isqrt\big[\mathcal{L}_2,\overline{\mathcal{L}}_1\big]\big)
\\
\rho_0\big(\isqrt\big[\mathcal{L}_1,\overline{\mathcal{L}}_2\big]\big)
&
\rho_0\big(\isqrt\big[\mathcal{L}_2,\overline{\mathcal{L}}_2\big]\big)
\end{array}
\!\right)(0)
&
=
\isqrt\,
\left(\!
\begin{array}{cc}
\mathcal{L}_1\big(\overline{A}_1\big)
-
\overline{\mathcal{L}}_1\big(A_1\big)
&
\mathcal{L}_2\big(\overline{A}_1\big)
-
\overline{\mathcal{L}}_1\big(A_2\big)
\\
\mathcal{L}_1\big(\overline{A}_2\big)
-
\overline{\mathcal{L}}_2\big(A_1\big)
&
\mathcal{L}_2\big(\overline{A}_2\big)
-
\overline{\mathcal{L}}_2\big(A_2\big)
\end{array}
\!\right)(0)
\\
&
=
\left(\!
\begin{array}{cc}
2a & 2\beta
\\
2\overline{\beta} & 2c
\end{array}
\!\right),
\endaligned
\]
and its determinant should be nonzero.

A linear transformation:
\[
\left(\!\!
\begin{array}{c}
z_1
\\
z_2
\end{array}
\!\!\right)
\,\longmapsto\,
\underbrace{
\left(\!\begin{array}{cc}
\lambda & \mu
\\
\nu & \chi
\end{array}
\!\right)}_{\in\,{\sf GL}_2(\C)}
\left(\!\!
\begin{array}{c}
z_1
\\
z_2
\end{array}
\!\!\right)
\]
normalizes this $2 \times 2$ Hermitian matrix to:
\[
2\,
\left(\!
\begin{array}{cc}
1 & 0
\\
0 & \pm 1
\end{array}
\!\right),
\]
which concludes.
\endproof

\noindent{\bf Scholium.}
{\em In such elementarily normalized coordinates, one has
the diagonal normalizations at the origin of the fields:}
\[
\aligned
\mathcal{L}_1\big\vert_0
&
=
\frac{\partial}{\partial z_1}
\bigg\vert_0,
\\
\mathcal{L}_2\big\vert_0
&
=
\frac{\partial}{\partial z_2}
\bigg\vert_0,
\\
\overline{\mathcal{L}}_1\big\vert_0
&
=
\frac{\partial}{\partial\overline{z}_1}
\bigg\vert_0,
\\
\overline{\mathcal{L}}_2\big\vert_0
&
=
\frac{\partial}{\partial\overline{z}_2}
\bigg\vert_0,
\\
\big[\mathcal{L}_1,\overline{\mathcal{L}}_1\big]
\Big\vert_0
&
=
-\,2\,\isqrt\,\,
\frac{\partial}{\partial u}
\bigg\vert_0,
\endaligned
\]
{\em and of the Levi Matrix:}
\[
\aligned
\left(\!
\begin{array}{cc}
\rho_0\big(\isqrt\big[\mathcal{L}_1,\overline{\mathcal{L}}_1\big]\big)
&
\rho_0\big(\isqrt\big[\mathcal{L}_2,\overline{\mathcal{L}}_1\big]\big)
\\
\rho_0\big(\isqrt\big[\mathcal{L}_1,\overline{\mathcal{L}}_2\big]\big)
&
\rho_0\big(\isqrt\big[\mathcal{L}_2,\overline{\mathcal{L}}_2\big]\big)
\end{array}
\!\right)(0)
&
=
2\,
\left(\!
\begin{array}{cc}
1 & 0
\\
0 & \pm 1
\end{array}
\!\right),
\endaligned
\]
{\em which conveniently fixes ideas when performing explicitly 
the Cartan equivalence procedure.\qed}

\medskip

Taking zero remainders, one obtains (two):
\[
\boxed{\,\,
{\text{\footnotesize\sf Model(s) 
$\text{\sf (IV)}_{\text{\sf 1}}$
:}}
\ \ \ \ \ \ \ \ \ 
v
=
z_1\overline{z}_1
\pm
z_1\overline{z}_2.\,\,
\rule[-4pt]{0pt}{15pt}}
\]


\bigskip

\section{\sf General class $\text{\sf IV}_{\text{\sf 2}}$}
\label{general-class-IV-2}
\HEAD{\ref{general-class-IV-2}.~General class 
$\text{\sf IV}_{\text{\sf 2}}$}{
Jo\"el {\sc Merker} (Paris-Sud)}

\medskip

\medskip\noindent{\bf Proposition.}
{\em A local real analytic hypersurface passing through the origin:}
\[
0
\,\in\,
M^5
\,\subset\,
\C^3
\]
{\em which belongs to the general Class $\text{\sf IV}_{
\text{\sf 2}}$, namely such that,
for any two local generators $\big\{
\mathcal{ L}_1, \mathcal{ L}_2 \big\}$ for $T^{1, 0}M$:}
\[
\Big\{
\mathcal{L}_1,\,\mathcal{L}_2,\,
\overline{\mathcal{L}}_1,\,\overline{\mathcal{L}}_2,\,\,
\big[\mathcal{L}_1,\overline{\mathcal{L}}_1\big]\Big\}
\]
{\em constitute a frame for $\C\otimes_\R TM$, 
such that in addition, the Levi form:}
\[
\text{\footnotesize\sf Levi-Form}^M(p)
\]
{\em if of rank $1$ at every point, and such that
lastly, the Freeman form:}
\[
\text{\sf Freeman-Form}^M(p)
\] 
{\em is nondegenerate at every point, 
may always be represented, in suitable local holomorphic coordinates:}
\[
(z_1,z_2,w)
\in
\C^3
\]
{\em by a specific real analytic equation of the form:}
\[
\boxed{\,\,
\underline{\text{\footnotesize\sf 
$\text{\sf (IV)}_{\text{\sf 2}}$:}}
\ \ \ \ \
\aligned
v
=
z_1\overline{z}_1
+
{\textstyle{\frac{1}{2}}}\,
z_1z_1\overline{z}_2
+
{\textstyle{\frac{1}{2}}}\,
z_2\overline{z}_1\overline{z}_1
+
{\rm O}_4\big(z_1,z_2,\overline{z}_1,\overline{z}_2\big)
+
u\,
{\rm O}_2\big(z_1,z_2,\overline{z}_1,\overline{z}_2,u\big),\,\,
\rule[-5pt]{0pt}{20pt}
\endaligned}
\]
{\em with remainders both satisfying:}
\[
0
\equiv
\remainder
\big(0,0,\overline{z}_1,\overline{z}_2,u\big)
\equiv
\remainder
\big(z_1,z_2,0,0,u\big).
\]

\proof
One starts with:
\[
\aligned
v
&
=
\varphi\big(z_1,z_2,\overline{z}_1,\overline{z}_2,u)
\\
&
=
{\rm O}_2\big(z_1,z_2,\overline{z}_1,\overline{z}_2,u\big),
\endaligned
\]
assuming:
\[
0
\equiv
\varphi\big(0,0,\overline{z}_1,\overline{z}_2,u\big)
\equiv
\varphi\big(z_1,z_2,0,0,u\big).
\]

Emphasize second and third order terms:
\[
\aligned
v
&
=
P_2\big(z_1,z_2,\overline{z}_1,\overline{z}_2\big)
+
P_3\big(z_1,z_2,\overline{z}_1,\overline{z}_2\big)
+
\\
&
\ \ \ \ \
+
u\,Q_2\big(z_1,z_2,\overline{z}_1,\overline{z}_2\big)
+
{\rm O}_4\big(z_1,z_2,\overline{z}_1,\overline{z}_2,u\big).
\endaligned
\]

Of course:
\[
\aligned
0
&
\equiv
P_3\big(0,0,\overline{z}_1,\overline{z}_2\big)
\equiv
P_3\big(z_1,z_2,0,0\big),
\\
0
&
\equiv
Q_2\big(0,0,\overline{z}_1,\overline{z}_2\big)
\equiv
Q_2\big(z_1,z_2,0,0\big),
\\
0
&
\equiv
{\rm O}_4\big(0,0,\overline{z}_1,\overline{z}_2,u\big)
\equiv
{\rm O}_4\big(z_1,z_2,0,0,u\big).
\endaligned
\]

As in what precedes:
\[
P_2\big(z_1,z_2,\overline{z}_1,\overline{z}_2\big)
=
a\,z_1\overline{z}_1
+
\beta\,z_2\overline{z}_1
+
\overline{\beta}\,z_1\overline{z}_2
+
c\,z_2\overline{z}_2,
\]
with:
\[
a\in\R,
\ \ \ \ \ \ \ \ \ \ \ \ \ \ \ \ \ \ \ \ \ \
\beta\in\C,
\ \ \ \ \ \ \ \ \ \ \ \ \ \ \ \ \ \ \ \ \ \
c\in\R.
\]

Taking:
\[
\rho_0
=
du
-
A_1\,dz_1
-
A_2\,dz_2
-
\overline{A}_1\,d\overline{z}_1
-
\overline{A}_2\,d\overline{z}_2,
\]
what was done for the General Class
$\text{\sf (IV)}_{\text\sf 1}$ yields that
the Levi matrix at the origin is:
\[
\!\!\!\!\!\!\!\!\!\!\!\!\!\!\!\!\!\!\!\!
\!\!\!\!\!\!\!\!\!\!\!\!\!\!\!\!\!\!\!\!
\aligned
\left(\!
\begin{array}{cc}
\rho_0\big(\isqrt\big[\mathcal{L}_1,\overline{\mathcal{L}}_1\big]\big)
&
\rho_0\big(\isqrt\big[\mathcal{L}_2,\overline{\mathcal{L}}_1\big]\big)
\\
\rho_0\big(\isqrt\big[\mathcal{L}_1,\overline{\mathcal{L}}_2\big]\big)
&
\rho_0\big(\isqrt\big[\mathcal{L}_2,\overline{\mathcal{L}}_2\big]\big)
\end{array}
\!\right)(0)
&
=
\isqrt\,
\left(\!
\begin{array}{cc}
\mathcal{L}_1\big(\overline{A}_1\big)
-
\overline{\mathcal{L}}_1\big(A_1\big)
&
\mathcal{L}_2\big(\overline{A}_1\big)
-
\overline{\mathcal{L}}_1\big(A_2\big)
\\
\mathcal{L}_1\big(\overline{A}_2\big)
-
\overline{\mathcal{L}}_2\big(A_1\big)
&
\mathcal{L}_2\big(\overline{A}_2\big)
-
\overline{\mathcal{L}}_2\big(A_2\big)
\end{array}
\!\right)(0)
\\
&
=
\left(\!
\begin{array}{cc}
2a & 2\beta
\\
2\overline{\beta} & 2c
\end{array}
\!\right),
\endaligned
\]
and it should be of rank $1$.

A linear transformation:
\[
\left(\!\!
\begin{array}{c}
z_1
\\
z_2
\end{array}
\!\!\right)
\,\longmapsto\,
\underbrace{
\left(\!\begin{array}{cc}
\lambda & \mu
\\
\nu & \chi
\end{array}
\!\right)}_{\in\,{\sf GL}_2(\C)}
\left(\!\!
\begin{array}{c}
z_1
\\
z_2
\end{array}
\!\!\right)
\]
normalizes it to:
\[
2\,
\left(\!
\begin{array}{cc}
1 & 0
\\
0 & 0
\end{array}
\!\right),
\]
whence:
\[
\aligned
v
&
=
z_1\overline{z}_1
+
P_3\big(z_1,z_2,\overline{z}_1,\overline{z}_2\big)
+
\\
&
\ \ \ \ \ \ \ \ \ \ \ \ \,
+
u\,Q_2\big(z_1,z_2,\overline{z}_1,\overline{z}_2\big)
+
{\rm O}_4\big(z_1,z_2,\overline{z}_1,\overline{z}_2,u\big).
\endaligned
\]

To compute:
\[
\overline{\mathcal{L}}_1\big(A_1\big),
\]
observing:
\[
\aligned
\varphi_u
&
=
{\rm O}(2),
\\
\varphi_{z_1}
&
=
{\rm O}(1),
\endaligned
\]
start with:
\[
\aligned
A_1
&
=
-\,\frac{\varphi_{z_1}}{\isqrt+\varphi_u}
\\
&
=
\isqrt\,
\bigg(
\frac{\varphi_{z_1}}{1-\isqrt\,\varphi_u}
\bigg)
\\
&
=
\isqrt\,\varphi_{z_1}\,
\Big(
1
+
\isqrt\,\varphi_u
+
\big(\isqrt\,\varphi_u\big)^2
+\cdots
\Big)
\\
&
=
\isqrt\,\varphi_{z_1}
+
{\rm O}(3)
\\
&
=
\isqrt\big(
\overline{z}_1
+
P_{3,z_1}
+
u\,Q_{2,z_1}
+
{\rm O}(3)
\big).
\endaligned
\]

Next:
\[
\aligned
\overline{\mathcal{L}}_1
&
=
\frac{\partial}{\partial\overline{z}_1}
+
\overline{A}_1\,
\frac{\partial}{\partial u}
\\
&
=
\frac{\partial}{\partial\overline{z}_1}
+
{\rm O}(1)\,\frac{\partial}{\partial u}.
\endaligned
\]
Observing that:
\[
\aligned
\varphi_{z_1u}
&
=
Q_{2,z_1}
+
{\rm O}(3)
\\
&
=
{\rm O}(1),
\endaligned
\]
one computes:
\[
\aligned
\overline{\mathcal{L}}_1\big(A_1\big)
&
=
\bigg(
\frac{\partial}{\partial\overline{z}_1}
+
{\rm O}(1)\,\frac{\partial}{\partial u}
\bigg)
\big[
\isqrt\,\varphi_{z_1}+{\rm O}(3)
\big]
\\
&
=
\isqrt\,\varphi_{z_1\overline{z}_1}
+
{\rm O}(1)\,\varphi_{z_1u}
+
{\rm O}(2)
\\
&
=
\isqrt\,\varphi_{z_1\overline{z}_1}
+
{\rm O}(2)
\\
&
=
\isqrt\,
\big(
1+P_{3,z_1\overline{z}_1}
+
u\,Q_{2,z_1\overline{z}_1}
+
{\rm O}(2)
\big).
\endaligned
\]

One treats similarly the other
\[
\overline{\mathcal{L}}_1\big(A_2\big),
\ \ \ \ \ \ \ \ \ \ \ \ \ \ \ \ \ \ \
\overline{\mathcal{L}}_2\big(A_1\big),
\ \ \ \ \ \ \ \ \ \ \ \ \ \ \ \ \ \ \
\overline{\mathcal{L}}_2\big(A_2\big),
\]
plus all conjugates, and obtains:
\[
\mathmotsf{Levi-Matrix}
\,=\,
2\,
\left(\!
\begin{array}{cc}
\substack{
1+P_{3,z_1\overline{z}_1}+\\
+uQ_{2,z_1\overline{z}_1}+{\rm O}(2)}
&
\substack{
P_{3,z_2\overline{z}_1}+\\
uQ_{2,z_2\overline{z}_1}+{\rm O}(2)}\bigskip
\\
\substack{
P_{3,z_1\overline{z}_2}+\\
+uQ_{2,z_1\overline{z}_2}+{\rm O}(2)}
&
\substack{
P_{3,z_2\overline{z}_2}+\\
uQ_{2,z_2\overline{z}_2}+{\rm O}(2)}
\end{array}
\!\right).
\]

The vanishing of its determinant modulo ${\rm O}(2)$-terms:
\[
0
\equiv
P_{3,z_2\overline{z}_2}
+
u\,Q_{2,z_2\overline{z}_2}
+
{\rm O}(2)
\]
yields:
\[
\aligned
0
&
\equiv
P_{3,z_2\overline{z}_2}
\\
0
&
\equiv
Q_{2,z_2\overline{z}_2}.
\endaligned
\]

The graphing equation then has the form:
\[
\aligned
v
&
=
z_1\overline{z}_1
+
\alpha\,z_1z_1\overline{z}_1
+
\overline{\alpha}_1\,z_1\overline{z}_1\overline{z}_1
+
\beta\,z_1z_2\overline{z}_2
+
\overline{\beta}\,z_2\overline{z}_1\overline{z}_2
+
\\
&
\ \ \ \ \ \ \ \ \ \ \ \ \
+
\gamma\,z_1z_2\overline{z}_1
+
\overline{\gamma}\,z_1\overline{z}_1\overline{z}_2
+
\delta\,z_2z_2\overline{z}_1
+
\overline{\delta}\,z_1\overline{z}_2\overline{z}_2
+
\\
&
\ \ \ \ \ \ \ \ \ \ \ \ \
+
a\,u\,z_1\overline{z}_1
+
\varepsilon\,u\,z_1\overline{z}_2
+
\overline{\varepsilon}\,u\,\overline{z}_1z_2
+
{\rm O}_4
\big(z_1,z_2,\overline{z}_1,\overline{z}_2,u\big).
\endaligned
\]

Doing:
\[
z_1
\,\longmapsto\,
z_1
+
\alpha\,z_1z_1
+
\gamma\,z_1z_2
+
\delta\,z_2z_2,
\]
one makes:
\[
\alpha
=
\gamma
=
\delta
=
0.
\]

It remains:

\[
\aligned
v
&
=
z_1\overline{z}_1
+
\beta\,z_1z_2\overline{z}_2
+
\overline{\beta}\,z_2\overline{z}_1\overline{z}_2
+
\\
&
\ \ \ \ \ \ \ \ \ \ \ \ \
+
a\,u\,z_1\overline{z}_1
+
\varepsilon\,u\,z_1\overline{z}_2
+
\overline{\varepsilon}\,u\,\overline{z}_1z_2
+
{\rm O}_4
\big(z_1,z_2,\overline{z}_1,\overline{z}_2,u\big).
\endaligned
\]

The function $k$ of Section~9, page~81 of~\cite{
Merker-Pocchiola-Sabzevari-5-CR-II} is:
\[
\aligned
k
&
=
\frac{
-\,\mathcal{L}_2\big(\overline{A}_1\big)
+
\overline{\mathcal{L}}_1\big(A_2\big)}{
\mathcal{L}_1\big(\overline{A}_1\big)
-
\overline{\mathcal{L}}_1\big(A_1\big)}
\\
&
=
\frac{2\isqrt\,\varphi_{z_2\overline{z}_1}+{\rm O}(1)}{
-\,2\,\isqrt+{\rm O}(1)}
\\
&
=
-\,\varphi_{z_2\overline{z}_1}
+
{\rm O}(2)
\\
&
=
-\,\overline{\beta}\,\overline{z}_1
-
\overline{\varepsilon}\,u
+
{\rm O}(2).
\endaligned
\]

Finally, the nondegeneracy of the Freeman form:
\[
0
\neq
\overline{\mathcal{L}}_1(k)
\]
reads:
\[
\aligned
\overline{\mathcal{L}}_1(k)
&
=
\bigg(
\frac{\partial}{\partial\overline{z}_1}
+
{\rm O}(1)\,\frac{\partial}{\partial u}
\bigg)
\Big[
-\,\overline{\beta}\,\overline{z}_1
-
\overline{\varepsilon}\,u
+
{\rm O}(2)
\Big]
\\
&
=
-\,\overline{\beta}
+
{\rm O}(1),
\endaligned
\]
which necessitates:
\[
\boxed{\,\beta\neq 0.\,}
\]

A complex dilation in the $z_2$-axis then terminates.
\endproof

\noindent{\bf Scholium.}
{\em In such elementarily normalized coordinates, one has
the diagonal normalizations at the origin of the fields:}
\[
\aligned
\mathcal{L}_1\big\vert_0
&
=
\frac{\partial}{\partial z_1}
\bigg\vert_0,
\\
\mathcal{L}_2\big\vert_0
&
=
\frac{\partial}{\partial z_2}
\bigg\vert_0,
\\
\overline{\mathcal{L}}_1\big\vert_0
&
=
\frac{\partial}{\partial\overline{z}_1}
\bigg\vert_0,
\\
\overline{\mathcal{L}}_2\big\vert_0
&
=
\frac{\partial}{\partial\overline{z}_2}
\bigg\vert_0,
\\
\big[\mathcal{L}_1,\overline{\mathcal{L}}_1\big]
\Big\vert_0
&
=
-\,2\,\isqrt\,\,
\frac{\partial}{\partial u}
\bigg\vert_0,
\endaligned
\]
{\em and of the Levi Matrix:}
\[
\aligned
\left(\!
\begin{array}{cc}
\rho_0\big(\isqrt\big[\mathcal{L}_1,\overline{\mathcal{L}}_1\big]\big)
&
\rho_0\big(\isqrt\big[\mathcal{L}_2,\overline{\mathcal{L}}_1\big]\big)
\\
\rho_0\big(\isqrt\big[\mathcal{L}_1,\overline{\mathcal{L}}_2\big]\big)
&
\rho_0\big(\isqrt\big[\mathcal{L}_2,\overline{\mathcal{L}}_2\big]\big)
\end{array}
\!\right)(0)
&
=
2\,
\left(\!
\begin{array}{cc}
1 & 0
\\
0 & 0
\end{array}
\!\right),
\endaligned
\]
{\em and furthermore, one has the normalization of the
function $k$:}
\[
\boxed{\,
k
=
-\,\overline{z}_1+{\rm O}(1),\,}
\]
{\em whence the nondegeneracy of the Freeman form reads:}
\[
\boxed{\,
\overline{\mathcal{L}}_1(k)
\big\vert_0
=
-\,1,\,}
\]
{\em which conveniently fixes ideas when performing explicitly 
the Cartan equivalence procedure.\qed}

\medskip

Taking, not zero remainders, but adapted ones
(\cite{ Gaussier-Merker-2003}):
\[
\boxed{\,\,
{\text{\footnotesize\sf Model 
$\text{\sf (III)}_{\text{\sf 1}}$:}}
\ \ \ \ \ \ \ \ \ 
v
=
\frac{z_1\overline{z}_1+
\frac{1}{2}\,z_1z_1\overline{z}_2
+
\frac{1}{2}\,z_2\overline{z}_1\overline{z}_1}{
1-z_2\overline{z}_2}.\,\,
\rule[-4pt]{0pt}{15pt}}
\]


\bigskip

\section{\sf Smooth categories}
\label{smooth-categories}
\HEAD{\ref{smooth-categories}.~General class 
$\text{\sf IV}_{\text{\sf 2}}$}{
Jo\"el {\sc Merker} (Paris-Sud)}

\medskip

Since the power series expansions have been used only up to a certain
order in all reasonings, the hypothesis of real analyticity is not
necessary. For instance, existence of preliminary coordinates in which
all pluriharmonic terms are removed are (well) known to also exist
when $M^{ 2n + c} \subset \C^{ n+c}$ is of class $\mathcal{
C}^\kappa$, with $\kappa \in \N$, $\kappa \geqslant 2$, or $\kappa =
\infty$, provided one truncates the requirements to a certain order
governed by the size of $\kappa$. Adapting the elementary proofs
provided here, one may obtain (exercise) easy generalizations of the
six propositions.


\bigskip

\section{\sf Graphing ontologies}
\label{graphing-ontologies}
\HEAD{\ref{graphing-ontologies}.~Graphing ontologies}{
Jo\"el {\sc Merker} (Paris-Sud)}

\medskip

\noindent{\bf General class $\text{\sf I}$.}
Before finer biholomorphic equivalence considerations, 
there is a one-to one correspondence between:
\[
\Big(
M^3
\,\subset\,
\C^2
\Big)
\,\,\in\,\,
\text{\sf General Class $\text{\sf I}$},
\]
and graphing functions: 
\[
\aligned
v
&
=
\varphi(x,y,u\big)
\\
&
=
z\overline{z}
+
z\overline{z}\,{\rm O}_1\big(z,\overline{z}\big)
+
z\overline{z}\,{\rm O}_1(u),
\endaligned
\]
of smoothness $\mathcal{ C}^\kappa$ ($\kappa \geqslant 3$), 
or $\mathcal{ C}^\infty$, or $\mathcal{ C}^\omega$, 
{\em without any further restriction}.

\medskip\noindent{\bf General class $\text{\sf II}$.}
Before finer biholomorphic equivalence considerations, 
there is a one-to one correspondence between:
\[
\Big(
M^4
\,\subset\,
\C^3
\Big)
\,\,\in\,\,
\text{\sf General Class $\text{\sf II}$},
\]
and pairs of graphing functions: 
\[
\aligned
v_1
&
=
\varphi_1\big(x,y,u_1,u_2\big)
\\
&
=
z\overline{z} 
+
z\overline{z}\,{\rm O}_2\big(z,\overline{z}\big)
+
z\overline{z}\,{\rm O}_1(u_1)
+
z\overline{z}\,{\rm O}_1(u_2),\,\,
\\
v_2
&
=
\varphi_2\big(x,y,u_1,u_2\big)
\\
&
=
z^2\overline{z}
+
z\overline{z}^2
+
z\overline{z}\,{\rm O}_2\big(z,\overline{z}\big)
+
z\overline{z}\,{\rm O}_1(u_1)
+
z\overline{z}\,{\rm O}_1(u_2).\,\,
\endaligned
\]
of smoothness $\mathcal{ C}^\kappa$ ($\kappa \geqslant 3$), 
or $\mathcal{ C}^\infty$, or $\mathcal{ C}^\omega$, 
again {\em without any further restriction}.

\medskip\noindent{\bf General class $\text{\sf III}_{\text{\sf 1}}$.}
Before finer biholomorphic equivalence considerations, 
there is a one-to one correspondence between:
\[
\Big(
M^5
\,\subset\,
\C^4
\Big)
\,\,\in\,\,
\text{\sf General Class $\text{\sf III}_{\text{\sf 1}}$},
\]
and triples of graphing functions: 
\[
\aligned
v_1
&
=
\varphi_1\big(x,y,u_1,u_2,u_3\big)
\\
&
=
z\overline{z}
+
z\overline{z}\,{\rm O}_2\big(z,\overline{z}\big)
+
z\overline{z}\,{\rm O}_1(u_1)
+
z\overline{z}\,{\rm O}_1(u_2)
+
z\overline{z}\,{\rm O}_1(u_3),
\endaligned
\]
\[
\aligned
v_2
&
=
\varphi_2\big(x,y,u_1,u_2,u_3\big)
\\
&
=
z^2\overline{z}
+
z\overline{z}^2
+
z\overline{z}\,{\rm O}_2\big(z,\overline{z}\big)
+
z\overline{z}\,{\rm O}_1(u_1)
+
z\overline{z}\,{\rm O}_1(u_2)
+
z\overline{z}\,{\rm O}_1(u_3),
\endaligned
\]
\[
\aligned
v_3
&
=
\varphi_3\big(x,y,u_1,u_2,u_3\big)
\\
&
=
\isqrt\,\big(z^2\overline{z}
-
z\overline{z}^2\big)
+
z\overline{z}\,{\rm O}_2\big(z,\overline{z}\big)
+
z\overline{z}\,{\rm O}_1(u_1)
+
z\overline{z}\,{\rm O}_1(u_2)
+
z\overline{z}\,{\rm O}_1(u_3),
\endaligned
\]
of smoothness $\mathcal{ C}^\kappa$ ($\kappa \geqslant 4$), 
or $\mathcal{ C}^\infty$, or $\mathcal{ C}^\omega$, 
still {\em without any further restriction}.

\medskip\noindent{\bf General class $\text{\sf III}_{\text{\sf 2}}$.}
Before finer biholomorphic equivalence considerations, 
there is a one-to one correspondence between:
\[
\Big(
M^5
\,\subset\,
\C^4
\Big)
\,\,\in\,\,
\text{\sf General Class $\text{\sf III}_{\text{\sf 2}}$},
\]
and triples of graphing functions: 
\[
\aligned
v_1
&
=
\varphi_1\big(x,y,u_1,u_2,u_3\big)
\\
&
=
z\overline{z}
+
z\overline{z}\,{\rm O}_3\big(z,\overline{z}\big)
+
z\overline{z}\,u_1\,{\rm O}_1\big(z,\overline{z},u_1\big)
+
\\
&
\ \ \ \ \ \ \ \ \ \ \ \ \ \ \ \ \ \ \ \ \ \ 
+
z\overline{z}\,u_2\,{\rm O}_1\big(z,\overline{z},u_1,u_2\big)
+
z\overline{z}\,u_3\,{\rm O}_1\big(z,\overline{z},u_1,u_2,u_3\big),
\endaligned
\]
\[
\aligned
v_2
&
=
\varphi_2\big(x,y,u_1,u_2,u_3\big)
\\
&
=
z^2\overline{z}
+
z\overline{z}^2
+
z\overline{z}\,{\rm O}_3\big(z,\overline{z}\big)
+
z\overline{z}\,u_1\,{\rm O}_1\big(z,\overline{z},u_1\big)
+
\\
&
\ \ \ \ \ \ \ \ \ \ \ \ \ \ \ \ \ \ \ \ \ \ 
+
z\overline{z}\,u_2\,{\rm O}_1\big(z,\overline{z},u_1,u_2\big)
+
z\overline{z}\,u_3\,{\rm O}_1\big(z,\overline{z},u_1,u_2,u_3\big),
\endaligned
\]
\[
\aligned
v_3
&
=
\varphi_3\big(x,y,u_1,u_2,u_3\big)
\\
&
=
2\,z^3\overline{z}
+
2\,z\overline{z}^3
+
3\,z^2\overline{z}^2
+
z\overline{z}\,{\rm O}_3\big(z,\overline{z}\big)
+
z\overline{z}\,u_1\,{\rm O}_1\big(z,\overline{z},u_1\big)
+
\\
&
\ \ \ \ \ \ \ \ \ \ \ \ \ \ \ \ \ \ \ \ \ \ \ \ \ \ \ \ \ \ \ \ \ \ \ 
\ \ \ \ \ \ 
+
z\overline{z}\,u_2\,{\rm O}_1\big(z,\overline{z},u_1,u_2\big)
+
z\overline{z}\,u_3\,{\rm O}_1\big(z,\overline{z},u_1,u_2,u_3\big),
\endaligned
\]
of smoothness $\mathcal{ C}^\kappa$ ($\kappa \geqslant 5$), 
or $\mathcal{ C}^\infty$, or $\mathcal{ C}^\omega$, 
{\em however with the restriction that:}
\[
\aligned
0
\equiv
\left\vert\!
\begin{array}{ccc}
\substack{
\mathcal{L}(\overline{A}_1)-\overline{\mathcal{L}}(A_1)}
&
\substack{
\mathcal{L}(\overline{A}_2)-\overline{\mathcal{L}}(A_2)}
&
\substack{
\mathcal{L}(\overline{A}_3)-\overline{\mathcal{L}}(A_3)}\bigskip
\\
\substack{
\mathcal{L}(\mathcal{L}(\overline{A}_1))
-2\mathcal{L}(\overline{\mathcal{L}}(A_1))+\\
+\overline{\mathcal{L}}(\mathcal{L}(A_1))}
&
\substack{
\mathcal{L}(\mathcal{L}(\overline{A}_2))
-2\mathcal{L}(\overline{\mathcal{L}}(A_2))+\\
+\overline{\mathcal{L}}(\mathcal{L}(A_2))}
&
\substack{
\mathcal{L}(\mathcal{L}(\overline{A}_3))
-2\mathcal{L}(\overline{\mathcal{L}}(A_3))+\\
+\overline{\mathcal{L}}(\mathcal{L}(A_3))}\bigskip
\\
\substack{
-\overline{\mathcal{L}}(\overline{\mathcal{L}}(A_1))
+2\overline{\mathcal{L}}(\mathcal{L}(\overline{A}_1))-\\
-\mathcal{L}(\overline{\mathcal{L}}(\overline{A}_1))}
&
\substack{
-\overline{\mathcal{L}}(\overline{\mathcal{L}}(A_2))
+2\overline{\mathcal{L}}(\mathcal{L}(\overline{A}_2))-\\
-\mathcal{L}(\overline{\mathcal{L}}(\overline{A}_2))}
&
\substack{
-\overline{\mathcal{L}}(\overline{\mathcal{L}}(A_3))
+2\overline{\mathcal{L}}(\mathcal{L}(\overline{A}_3))-\\
-\mathcal{L}(\overline{\mathcal{L}}(\overline{A}_3))}
\end{array}
\!\right\vert.
\endaligned
\]
a possibly subtle and complicated
differential relation\big/condition one should take account of 
when launching the Cartan equivalence method.

\medskip\noindent{\bf General class $\text{\sf IV}_{\text{\sf 1}}$.}
Before finer biholomorphic equivalence considerations, 
there is a one-to one correspondence between:
\[
\Big(
M^5
\,\subset\,
\C^3
\Big)
\,\,\in\,\,
\text{\sf General Class $\text{\sf IV}_{\text{\sf 1}}$},
\]
and graphing functions: 
\[
\aligned
v
&
=
\varphi\big(x_1,x_2,y_1,y_2,u\big)
\\
&
=
z_1\overline{z}_1
\pm
z_2\overline{z}_2
+
{\rm O}_3\big(z_1,z_2,\overline{z}_1,\overline{z}_2,u\big),\,\,
\rule[-5pt]{0pt}{20pt}
\endaligned
\]
of smoothness $\mathcal{ C}^\kappa$ ($\kappa \geqslant 3$), 
or $\mathcal{ C}^\infty$, or $\mathcal{ C}^\omega$, 
{\em without any further restriction} for the fourth time.

\medskip\noindent{\bf General class $\text{\sf IV}_{\text{\sf 2}}$.}
Before finer biholomorphic equivalence considerations, 
there is a one-to one correspondence between:
\[
\Big(
M^5
\,\subset\,
\C^3
\Big)
\,\,\in\,\,
\text{\sf General Class $\text{\sf IV}_{\text{\sf 1}}$},
\]
and graphing functions: 
\[
\aligned
v
&
=
\varphi\big(x_1,x_2,y_1,y_2,u\big)
\\
&
=
z_1\overline{z}_1
+
z_1z_1\overline{z}_2
+
z_2\overline{z}_1\overline{z}_1
+
{\rm O}_4\big(z_1,z_2,\overline{z}_1,\overline{z}_2\big)
+
u\,
{\rm O}_2\big(z_1,z_2,\overline{z}_1,\overline{z}_2,u\big),
\endaligned
\]
of smoothness $\mathcal{ C}^\kappa$ ($\kappa \geqslant 4$), 
or $\mathcal{ C}^\infty$, or $\mathcal{ C}^\omega$, 
{\em however with the restriction of identical vanishing:}
\[
\aligned
0
\equiv\,\,
&
\mathmotsf{Levi-Determinant}
=
\frac{4}{(\isqrt+\varphi_u)^3(-\isqrt+\varphi_u)^3}
\bigg\{
\varphi_{z_2\overline{z}_2}\varphi_{z_1\overline{z}_1}
-
\varphi_{z_2\overline{z}_1}\varphi_{z_1\overline{z}_2}
+
\\
&
+
\varphi_{z_2\overline{z}_1}\varphi_{\overline{z}_2}\varphi_{z_1u}
\varphi_u
-
\varphi_{z_2\overline{z}_1}\varphi_{\overline{z}_2}\varphi_{z_1}
\varphi_{uu}
-
\varphi_{\overline{z}_1}\varphi_{z_2u}\varphi_{z_1}
\varphi_{\overline{z}_2u}
+
\varphi_{\overline{z}_1}\varphi_{z_2u}\varphi_u
\varphi_{z_1\overline{z}_2}
-
\\
&
-\,
\varphi_{z_2}\varphi_{\overline{z}_1u}\varphi_{\overline{z}_2}
\varphi_{z_1u}
-
\varphi_{z_2}\varphi_{\overline{z}_1}\varphi_{uu}
\varphi_{z_1\overline{z}_2}
+
\varphi_{z_2}\varphi_{\overline{z}_1u}
\varphi_u\varphi_{z_1\overline{z}_2}
-
\varphi_{z_2\overline{z}_2}\varphi_{\overline{z}_1}
\varphi_{z_1u}\varphi_u
+
\\
&
+
\varphi_{z_2\overline{z}_2}\varphi_{z_1}
\varphi_{\overline{z}_1}\varphi_{uu}
-
\varphi_{z_2\overline{z}_2}\varphi_{z_1}
\varphi_{\overline{z}_1u}\varphi_u
+
\varphi_{z_2\overline{z}_1}\varphi_{z_1}
\varphi_{\overline{z}_2u}\varphi_u
+
\varphi_{z_2}\varphi_{\overline{z}_2u}
\varphi_{\overline{z}_1}\varphi_{z_1u}
-
\\
&
-\,
\varphi_{z_2}\varphi_{\overline{z}_2u}
\varphi_{z_1\overline{z}_1}\varphi_u
+
\varphi_{\overline{z}_2}\varphi_{z_2u}
\varphi_{z_1}\varphi_{\overline{z}_1u}
-
\varphi_{\overline{z}_2}\varphi_{z_2u}
\varphi_u\varphi_{z_1\overline{z}_1}
+
\varphi_{\overline{z}_2}\varphi_{z_2}
\varphi_{uu}\varphi_{z_1\overline{z}_1}
+
\\
&
+
\isqrt
\big(
\varphi_{z_2\overline{z}_2}\varphi_{z_1}
\varphi_{\overline{z}_1u}
+
\varphi_{\overline{z}_1}\varphi_{z_2u}
\varphi_{z_1\overline{z}_2}
+
\varphi_{z_2\overline{z}_1}\varphi_{\overline{z}_2}
\varphi_{z_1u}
+
\varphi_{z_2}\varphi_{\overline{z}_2u}
\varphi_{z_1\overline{z}_1}
\big)
-
\\
&
-\,\isqrt\big(
\varphi_{\overline{z}_2}\varphi_{z_2u}
\varphi_{z_1\overline{z}_1}
+
\varphi_{z_2\overline{z}_1}\varphi_{z_1}
\varphi_{\overline{z}_2u}
+
\varphi_{z_2}\varphi_{\overline{z}_1u}
\varphi_{z_1\overline{z}_2}
+
\varphi_{z_2\overline{z}_2}\varphi_{\overline{z}_1}
\varphi_{z_1u}
\big)
-
\\
&
-\,
\varphi_{z_2\overline{z}_1}\varphi_{z_1\overline{z}_2}
\varphi_u\varphi_u
+
\varphi_{z_2\overline{z}_2}\varphi_{z_1\overline{z}_1}
\varphi_u\varphi_u
\bigg\},\,\,
\endaligned
\]
a possibly subtle and complicated differential relation\big/condition
that Samuel Pocchiola took systematically account of while performing
{\em completely explicitly} the Cartan equivalence method in~\cite{
Pocchiola-2013}.


\vfill\end{document}